\documentclass[twoside,12pt]{article}
\usepackage[dvips]{graphicx}

\begin{document}

\newtheorem{lem}{Lemma}[section]
\newtheorem{pro}[lem]{Proposition}
\newtheorem{theo}[lem]{Theorem}
\newtheorem{cor}[lem]{Corollary}
\newtheorem{ass}[lem]{Assumption}
\newcommand{\qed}{\hfill \hbox{\rule{4pt}{8pt}}}

\renewcommand{\textfraction}{.01}
\renewcommand{\topfraction}{.99}

\title{{\bf 
Floer's chain complexes \\
for Lagrangian submanifolds \\ 
in symplectic manifolds \\
with concave ends
}}

\author{Manabu AKAHO 
\\ Department of Mathematics
\\ Tokyo Metropolitan University 
\\ {\tt akaho@math.metro-u.ac.jp}}

\date{}

\maketitle

%%%%%%%%%%%%%%%%%%%%%%%%%%%%%%%%%%%%%%%%%%%%%%%%%%%%%%%%%%%%%%%%%%
\section{Introduction}%%%%%%%%%%%%%%%%%%%%%%%%%%%%%%%%%%%%%%%%%%%%
%%%%%%%%%%%%%%%%%%%%%%%%%%%%%%%%%%%%%%%%%%%%%%%%%%%%%%%%%%%%%%%%%%

We usually define Floer homology for Lagrangian submanifolds
in symplectic manifolds,
but there are many important {\it non-smooth}
guys such as algebraic varieties and Lagrangian subvarieties.
Our first plan to construct Floer theory for such non-smooth
objects is to do for something like open strata of them. 
The open strata are non-compact,
and we will start with {\it concave} ends.

Floer's chain complexes for Lagrangian submanifolds
in closed symplectic manifolds are generated by
intersection points of Lagrangian submanifolds
and whose differentials count pseudo-holomorphic strips 
with Lagrangian boundary conditions.
In this paper we will propose
Floer's chain complexes for Lagrangian submanifolds 
in symplectic manifolds with {\it concave ends}.
\\

A {\it symplectic form} $\omega $
on a smooth manifold $X$ is a non-degenerate closed $2$-form.
The non-degeneracy induces the existence of 
almost complex structures $J$ such that
$\omega (\cdot ,J\cdot )$ is a metric on $X$.
In particular, we consider 
time-dependent almost complex structures $J_t, t\in [0,1]$.

A {\it Lagrangian submanifold} $L$ 
is an $n$-dimensional submanifold in $X^{2n}$
such that $\omega |_{TL}=0$.
Here we assume the following conditions
for Lagrangian submanifolds $L_0$ and $L_1$.
\begin{ass}(Nondegeneracy of intersections)\label{ass1.1}%%%%%%%%%
$L_0$ and $L_1$ intersect transversally.
\end{ass}%%%%%%%%%%%%%%%%%%%%%%%%%%%%%%%%%%%%%%%%%%%%%%%%%%%%%%%%%
\begin{ass}(Admissibility)\label{ass1.2}%%%%%%%%%%%%%%%%%%%%%%%%%%
If $u:S^1\times [0,1]\to X$ is a map such that 
$u(S^1,0)\subset L_0$ and $u(S^1,1)\subset L_1$,
then $\int_{S^1\times [0,1]}u^*\omega =0$.
\end{ass}%%%%%%%%%%%%%%%%%%%%%%%%%%%%%%%%%%%%%%%%%%%%%%%%%%%%%%%%%

Let $(X,J)$ be an almost complex manifold 
and $(\Sigma ,j)$ a Riemann surface with a complex structure $j$.
A pseudo-holomorphic curve is a map $u:\Sigma \to X$ such that
\[
\overline{\partial }_Ju:=\frac{1}{2}(du+J\circ du \circ j)=0.
\]
Our Riemann surface is ${\bf R}\times [0,1]$ 
with the natural complex structure $i$,
and our almost complex structures on $X$ are time-dependent,
then we consider the following 
elliptic partial differential equation
\[
\overline{\partial }_{J_t}u(\tau ,t)
:=\frac{\partial u(\tau ,t)}{\partial \tau }+
J_t(u(\tau ,t))\frac{\partial u(\tau ,t)}{\partial t}=0,
\]
where $(\tau ,t)\in {\bf R}\times [0,1]$.
Define $\mathcal{M}(p,q)$, for $p$ and $q\in L_0\cap L_1$,
to be the set of maps $u:{\bf R}\times [0,1]\to X$ such that 
\begin{itemize}
\item $u({\bf R},0)\subset L_0$ and $u({\bf R},1)\subset L_1$
\item $\lim_{\tau \to -\infty }u(\tau ,[0,1])=p$ and 
      $\lim_{\tau \to \infty }u(\tau ,[0,1])=q$
\item $\overline{\partial }_{J_t}u(\tau , t)=0$.
\end{itemize}
We call a map satisfying the above conditions a
{\it pseudo-holomorphic strip}
and $\mathcal{M}(p,q)$ 
the {\it moduli space of pseudo-holomorphic strips}.
Note that ${\bf R}$ acts on $\mathcal{M}(p,q)$
by $(a*u)(\tau ,t):=u(\tau -a,t)$ with $a\in {\bf R}$.
We denote the quotient by $\hat{\mathcal{M}}(p,q)$.
Then Floer proved the following theorem \cite{f1}.
\begin{theo}\label{theo1.3}%%%%%%%%%%%%%%%%%%%%%%%%%%%%%%%%%%%%%%%
If $J_t$ is generic, then $\mathcal{M}(p,q)$ 
is a smooth finite dimensional manifold.
\end{theo}%%%%%%%%%%%%%%%%%%%%%%%%%%%%%%%%%%%%%%%%%%%%%%%%%%%%%%%%
The generic means that $J_t$ is an element 
of a Baire category set in a certain Banach space
of time-dependent almost complex structures. 

For the compactification of the moduli spaces in Floer's context,
we need the following condition.
\begin{ass}($\pi _2$-condition)\label{ass1.4}%%%%%%%%%%%%%%%%%%%%%
If $u:D^2\to X$ is a map such that $u(\partial D^2)\subset L$,
then $\int _{D^2}u^*\omega =0$.
\end{ass}%%%%%%%%%%%%%%%%%%%%%%%%%%%%%%%%%%%%%%%%%%%%%%%%%%%%%%%%%
Assume that $L_0$ and $L_1$ satisfy the $\pi _2$-condition.
\begin{theo}\label{theo1.5}%%%%%%%%%%%%%%%%%%%%%%%%%%%%%%%%%%%%%%%
Let $X$ be a closed symplectic manifold.
Then $\hat{\mathcal{M}}(p,q)$ can be compactified
(with respect to the topology of uniform convergence with all
derivatives on compact set).
\end{theo}%%%%%%%%%%%%%%%%%%%%%%%%%%%%%%%%%%%%%%%%%%%%%%%%%%%%%%%%
In fact,
if the dimension of $\hat{\mathcal{M}}(p,q)$ is equal to $0$,
then $\hat{\mathcal{M}}(p,q)$ is compact, 
and if the dimension of $\hat{\mathcal{M}}(p,q)$ is equal to $1$,
then it can be compactified so that the boundary is
\begin{eqnarray}
\bigcup_{\dim \hat{\mathcal{M}}(p,r)=
\dim \hat{\mathcal{M}}(r,q)=0}
\hat{\mathcal{M}}(p,r)\times \hat{\mathcal{M}}(r,q).
\label{1}
\end{eqnarray}
(To show $(\ref{1})$ we need also the gluing theorem \cite{f1}.)

Let $C$ be the free ${\bf Z}_2$-vector space 
over the elements of $L_0\cap L_1$.
We define a linear map $\partial :C\to C$ 
in terms of the canonical bases
\[
\partial p
:=\sum _{q\in L_0\cap L_1}
\sharp \hat{\mathcal{M}}(p,q)q,
\]
where the sum ranges over all $q\in L_0\cap L_1$
such that $\dim \hat{\mathcal{M}}(p,q)=0$
and where $\sharp \hat{\mathcal{M}}(p,q)$ is 
the modulo 2 number of the elements of $\hat{\mathcal{M}}(p,q)$.
Then Floer proved the following theorem \cite{f1}.
\begin{theo}\label{theo1.6}%%%%%%%%%%%%%%%%%%%%%%%%%%%%%%%%%%%%%%%
$\partial \circ \partial =0$.
\end{theo}%%%%%%%%%%%%%%%%%%%%%%%%%%%%%%%%%%%%%%%%%%%%%%%%%%%%%%%%
The idea of Theorem \ref{theo1.6} is very important for us, 
hence we adopt the proof.
For $p\in L_0\cap L_1$
\begin{eqnarray*}
\partial \partial p
&=&
\partial \sum _{q\in L_0\cap L_1} \sharp \hat{\mathcal{M}}(p,q)q
\\
&=&
\sum _{r\in L_0\cap L_1}\sum _{q\in L_0\cap L_1}
\sharp \hat{\mathcal{M}}(p,q)\sharp \hat{\mathcal{M}}(q,r)r.
\end{eqnarray*}
$\sum _{q\in L_0\cap L_1}
\sharp \hat{\mathcal{M}}(p,q)\sharp \hat{\mathcal{M}}(q,r)$
is nothing but the number of the boundary components of
the compactification of 
a $1$-dimensional smooth manifold $\hat{\mathcal{M}}(p,r)$,
and then even.

We call the chain complex $(C,\partial )$ the 
{\it Floer's chain complex} for $L_0$ and $L_1$ in $X$.
\\

By using the universal Novikov ring as in \cite{fooo}, 
we can remove the admissibility.
On the other hand,
many persons made effort to weaken the $\pi _2$-condition 
\cite{oh} and \cite{o},
and it grows into an obstruction theory of the boundary operators
\cite{fooo}. 
\\

Next we will consider $(\ref{1})$ from another angle.
Let $\{u_i\}_{i=1,2,\ldots }$ be a sequence of elements in
$\hat{\mathcal{M}}(p,q)$ which converges to an element 
$(v,w)\in \hat{\mathcal{M}(p,r)}\times \hat{\mathcal{M}}(r,q)$,
see Figure 1.
%%%%%%%%%%%%%%%%%%%%%%%%%%%%%%%%%%%%%%%%%%%%%%%%%%%%%%%%%%%%%%%%%%
\begin{figure}[h]
\setlength{\unitlength}{1mm}
\begin{picture}(0,60)(0,0)
%%%% waku
%\put(0,0){\framebox(130,60){}}
%\put(65,0){\line(0,1){60}}
%\put(32.5,0){\line(0,1){60}}
%\put(97.5,0){\line(0,1){60}}
%\put(0,5){\line(1,0){130}}
%\put(0,30){\line(1,0){130}}
%\put(0,55){\line(1,0){130}}
%%%% u_i tachi
\put(15,45){\makebox(5,5){$u_i$}}
\bezier{230}(32.5,55)(17.5,42.5)(32.5,30)
\bezier{230}(32.5,30)(47.5,17.5)(32.5,5)
\bezier{230}(32.5,55)(47.5,42.5)(34.5,32)
\bezier{230}(30.5,28)(17.5,17.5)(32.5,5)
%%%% kouten 
\put(30,0){\makebox(5,5){$q$}}
\put(30,55){\makebox(5,5){$p$}}
%%%% yajirushi
\put(60,25){\makebox(10,10){$\Longrightarrow$}}
%%%% v to w
\put(95,40){\makebox(5,5){$v$}}
\bezier{230}(97.5,55)(82.5,42.5)(97.5,30)
\bezier{230}(97.5,30)(112.5,17.5)(97.5,5)
\put(95,15){\makebox(5,5){$w$}}
\bezier{230}(97.5,55)(112.5,42.5)(97.5,30)
\bezier{230}(97.5,30)(82.5,17.5)(97.5,5)
%%%% kouten
\put(95,0){\makebox(5,5){$q$}}
\put(95,55){\makebox(5,5){$p$}}
\put(100,27.5){\makebox(5,5){$r$}}
%%%%
\end{picture}
\begin{center}
Figure 1
\end{center}
\end{figure}
%%%%%%%%%%%%%%%%%%%%%%%%%%%%%%%%%%%%%%%%%%%%%%%%%%%%%%%%%%%%%%%%%%
\\
This phenomenon implies that,
{\it at the limit of the sequence, a generator} $r$ 
{\it of the Floer's chain complex appears 
at the intersection point} $r$.
For simplicity, we assume that
$X$, around $r$, is locally isomorphic to ${\bf C}^n$,
where $r$ corresponds to the origin,
and $L_0$ and $L_1$ are locally isomorphic to ${\bf R}^n$
and $(\sqrt{-1}{\bf R})^n$, respectively.
Consider ${\bf C}^n\setminus \{0\}$ 
to be $(0,\infty )\times S^{2n-1}$
through the polar coordinate, 
and moreover transform it into ${\bf R}\times S^{2n-1}$ by 
\[
(0,\infty )\times S^{2n-1}\to {\bf R}\times S^{2n-1},
\
(\rho , x)\mapsto (\alpha ,x):=(\log \rho ,x),
\]
%%%%%%%%%%%%%%%%%%%%%%%%%%%%%%%%%%%%%%%%%%%%%%%%%%%%%%%%%%%%%%%%%%
\begin{figure}[h]
\setlength{\unitlength}{1mm}
\begin{picture}(0,140)(0,0)
%%%% waku
%\put(0,0){\framebox(130,140){}}
%\put(65,0){\line(0,1){140}}
%\put(0,60){\line(1,0){130}}
%\put(0,70){\line(1,0){130}}
%\put(0,80){\line(1,0){130}}
%%%% hidari no tsutsu
\put(17.5,40){\line(0,1){55}}
\bezier{50}(17.45,94.98)(17.45,100)(16,100)
\put(47.5,40){\line(0,1){55}}
\bezier{50}(47.4,94.98)(47.4,100)(49,100)
%%%% hidari no ue wakka
\bezier{170}(17.5,85)(17.5,80)(32.5,80)
\bezier{170}(32.5,80)(47.5,80)(47.5,85)
\bezier{20}(17.5,85)(17.5,90)(32.5,90)
\bezier{20}(32.5,90)(47.5,90)(47.5,85)
%%%% hidari no shita wakka
\bezier{170}(17.5,55)(17.5,50)(32.5,50)
\bezier{170}(32.5,50)(47.5,50)(47.5,55)
\bezier{20}(17.5,55)(17.5,60)(32.5,60)
\bezier{20}(32.5,60)(47.5,60)(47.5,55)
%%%% migiue no tsutsu
\put(82.5,80){\line(0,1){55}}
\bezier{50}(82.45,134.98)(82.45,140)(81,140)
\put(112.5,80){\line(0,1){55}}
\bezier{50}(112.4,134.98)(112.4,140)(114,140)
%%%% migiue no ue wakka
\bezier{170}(82.5,125)(82.5,120)(97.5,120)
\bezier{170}(97.5,120)(112.5,120)(112.5,125)
\bezier{20}(82.5,125)(82.5,130)(97.5,130)
\bezier{20}(97.5,130)(112.5,130)(112.5,125)
%%%% migiue no shita wakka
\bezier{170}(82.5,95)(82.5,90)(97.5,90)
\bezier{170}(97.5,90)(112.5,90)(112.5,95)
\bezier{20}(82.5,95)(82.5,100)(97.5,100)
\bezier{20}(97.5,100)(112.5,100)(112.5,95)
%%%% migi no wakka
\bezier{170}(82.5,70)(82.5,65)(97.5,65)
\bezier{170}(97.5,65)(112.5,65)(112.5,70)
\bezier{170}(82.5,70)(82.5,75)(97.5,75)
\bezier{170}(97.5,75)(112.5,75)(112.5,70)
%%%% migishita no tsutsu
\put(82.5,0){\line(0,1){60}}
\put(112.5,0){\line(0,1){60}}
%%%% migiue no ue wakka
\bezier{170}(82.5,45)(82.5,40)(97.5,40)
\bezier{170}(97.5,40)(112.5,40)(112.5,45)
\bezier{20}(82.5,45)(82.5,50)(97.5,50)
\bezier{20}(97.5,50)(112.5,50)(112.5,45)
%%%% migiue no shita wakka
\bezier{170}(82.5,15)(82.5,10)(97.5,10)
\bezier{170}(97.5,10)(112.5,10)(112.5,15)
\bezier{20}(82.5,15)(82.5,20)(97.5,20)
\bezier{20}(97.5,20)(112.5,20)(112.5,15)
%%%% yajirushi
\put(60,65){\makebox(10,10){$\Longrightarrow$}}
%%%% hidari ue \Lambda tachi
\put(21,84){\line(2,1){9}}
\put(44,86){\line(-2,-1){9}}
%%%% hidari shita \Lambda tachi
\put(21,54){\line(2,1){9}}
\put(44,56){\line(-2,-1){9}}
%%%% migi no \Lambda tachi with Reeb chords
\put(86,69){\line(2,1){9}}
\put(109,71){\line(-2,-1){9}}
\put(88,70){\line(5,-1){13.53}}
\put(107,70){\line(-5,1){13.53}}
%%%% migiue ue no \Lambda tachi
\put(86,124){\line(2,1){9}}
\put(109,126){\line(-2,-1){9}}
%%%% migiue shita no \Lambda tachi
\put(86,94){\line(2,1){9}}
\put(109,96){\line(-2,-1){9}}
%%%% migishita ue no \Lambda tachi
\put(86,44){\line(2,1){9}}
\put(109,46){\line(-2,-1){9}}
%%%% migishita shita no \Lambda tachi
\put(86,14){\line(2,1){9}}
\put(109,16){\line(-2,-1){9}}
%%%%
\put(75.5,68){\makebox(10,10){$\Lambda _0$}}
\put(109.5,62){\makebox(10,10){$\Lambda _1$}}
%%%% hidari no strip
\bezier{400}(23,95)(23,45)(26.25,45)
\bezier{50}(28.5,97)(28.5,45)(26.25,45)
\bezier{400}(36,92.4)(36,42.4)(39.25,42.4)
\bezier{400}(41.5,94.4)(41.5,42.4)(39.25,42.4)
\put(26.25,45){\line(5,-1){13}}
%%%% migishita no strip
\bezier{400}(88,55)(88,5)(91.25,5)
\bezier{50}(93.5,57)(93.5,5)(91.25,5)
\bezier{400}(101,52.4)(101,2.4)(104.25,2.4)
\bezier{400}(106.5,54.4)(106.5,2.4)(104.25,2.4)
\put(91.25,5){\line(5,-1){13}}
%%%% migiue no strips
\put(88,85){\line(0,1){50}}
\put(101,82.4){\line(0,1){50}}
\put(107,82.4){\line(0,1){50}}
\multiput(94,85)(0,2){25}{\line(0,1){0.5}}
%%%%% migi no tensen
\multiput(88,55)(0,2){15}{\line(0,1){0.5}}
\multiput(101,52.4)(0,2){15}{\line(0,1){0.5}}
\multiput(94,57)(0,2){15}{\line(0,1){0.5}}
\multiput(107,54.4)(0,2){15}{\line(0,1){0.5}}
%%%% \infty toka
\put(49.5,38){\makebox(5,5){$-\infty $}}
\put(114.5,-2){\makebox(5,5){$-\infty $}}
\put(114.5,57){\makebox(5,5){$\infty $}}
\put(114.5,78){\makebox(5,5){$-\infty $}}
%%%% strips no namae
\put(30,100){\makebox(5,5){$u_i$}}
\put(91,134){\makebox(5,5){$v$}}
\put(99,136){\makebox(5,5){$w$}}
\put(95,-2.5){\makebox(5,5){$z$}}
%%%%
\end{picture}
\begin{center}
Figure 2
\end{center}
\end{figure}
%%%%%%%%%%%%%%%%%%%%%%%%%%%%%%%%%%%%%%%%%%%%%%%%%%%%%%%%%%%%%%%%%%
\\
so that $L_0$ and $L_1$ are locally diffeomorphic to
${\bf R}\times \Lambda _0$ and ${\bf R}\times \Lambda _1$,
respectively,
where $\Lambda _0$ and $\Lambda _1$ are (Legendrian) submanifolds
in $S^{2n-1}$ (with the standard contact form).
Then the limit of $\{u_i\}_{i=1,2,\ldots }$ 
looks like the following:
$u_i$'s grow toward $-\infty$, and at the limit there appear
three pseudo-holomorphic strips $v$, $w$ and $z$ 
as in Figure 2.
(The segments between $\Lambda _0$ and $\Lambda _1$ are 
Reeb chords.)
\\

In comparison with the intersection point $r$,
this phenomenon implies that, 
{\it at the limit of the sequence, a generator} $z$
{\it of the Floer's chain complex for} 
$L_0\setminus \{r\}$ {\it and} $L_1\setminus \{r\}$ 
{\it in} $X\setminus \{r\}$ 
{\it appears at} $-\infty $ {\it of the (concave) end}.
\\

The content of our paper is as follows:
Section 2 defines concave/convex ends 
of non-compact symplectic manifolds 
and pseudo-holomorphic strips.
Section 3 proposes Floer's chain complexes 
for Lagrangian submanifolds 
in symplectic manifolds with concave/convex ends.
Section 4 proves some gluing arguments 
for pseudo-holomorphic strips and 
Section 5 observes the bubbling off phenomenon 
for pseudo-holomorphic curves. 
%Section 6 considers the transversality problem in our cases.
\\
\\
{\it Sign convention}; through this paper,
$L^p_1$ is the space of functions/sections $f$
such that $\int |f|^p +\int |Df|^p< \infty $ {\it for} $p> 2$.

%%%%%%%%%%%%%%%%%%%%%%%%%%%%%%%%%%%%%%%%%%%%%%%%%%%%%%%%%%%%%%%%%%
\section{Pseudo-holomorphic strips %%%%%%%%%%%%%%%%%%%%%%%%%%%%%%%
in symplectic manifolds with concave ends}%%%%%%%%%%%%%%%%%%%%%%%%
%%%%%%%%%%%%%%%%%%%%%%%%%%%%%%%%%%%%%%%%%%%%%%%%%%%%%%%%%%%%%%%%%%

Let $M$ be a smooth oriented manifold of dimension $2n+1$.
A {\it contact form} on $M$ is a 1-form $\lambda $
such that $\lambda \wedge (d\lambda )^n$ is a volume form on $M$.
Then the 2-form $d(e^{\alpha }\lambda )$ is a symplectic structure
in ${\bf R}\times M$, where $\alpha \in {\bf R}$.
We call $({\bf R}\times M,d(e^{\alpha }\lambda ))$ 
the {\it symplectization} of $(M,\lambda )$,
and there is a natural projection $\pi _M:{\bf R}\times M\to M$.

Let $X$ be a non-compact symplectic manifold
which, out side a compact set, is diffeomorphic to cylinders.
If a cylinder is  symplectically isomorphic to 
$((-\infty ,R_-)\times M_-,d(e^{\alpha }\lambda _-))$
for $(M_-,\lambda _-)$ and  $R_-\in {\bf R}$,
then we call the cylinder a {\it concave end}.
Similarly, if a cylinder is  symplectically isomorphic to 
$((R_+,\infty )\times M_+,d(e^{\alpha }\lambda _+))$
for $(M_+,\lambda _+)$ and $R_+\in {\bf R}$,
then we call the cylinder a {\it convex end}.
\\

A {\it Legendrian submanifold} $\Lambda $ 
is an $n$-dimensional submanifold in $(M^{2n+1},\lambda )$ 
which satisfies $T\Lambda \subset \mbox{Ker} \lambda $.
Then ${\bf R}\times \Lambda $ is a Lagrangian submanifold
in the symplectization. 

Let $L$ be a Lagrangian submanifold in a non-compact symplectic
manifold whose ends are concave or convex. 
If $L$ is non-compact, then we assume that 
$L$ satisfies the following.
\begin{ass}(Cone-condition)\label{ass2.1}%%%%%%%%%%%%%%%%%%%%%%%%%
$L|_{(-\infty ,R_-)\times M_-}=(-\infty ,R_-)\times \Lambda _-$
and
$L|_{(R_+,\infty )\times M_+}=(R_+,\infty )\times \Lambda _+$,
where $\Lambda _-$ and $\Lambda _+$ are Legendrian submanifolds.
\end{ass}%%%%%%%%%%%%%%%%%%%%%%%%%%%%%%%%%%%%%%%%%%%%%%%%%%%%%%%%%

Associated to $\lambda $ there are two important structures.
First of all the so-called Reeb vector field $X=X_{\lambda }$ 
defined by $i_X\lambda \equiv 1$ and $i_Xd\lambda \equiv 0$,
and secondly the contact structure $\xi =\xi _{\lambda }$ 
given by $\xi :=\mbox{Ker} \lambda \subset TM$.
$TM$ splits into ${\bf R}X_{\lambda }\oplus \xi$,
and we have a natural projection $\pi _{\xi }:TM\to \xi$.
(We shall use the same notation $\xi $ to denote $\pi _M^*\xi$.)
Moreover, on the symplectization, 
$T({\bf R}\times M)$ splits into $E\oplus \xi$, 
where $E:={\bf R}\frac{\partial }{\partial \theta}
\oplus \pi _M^*{\bf R}X_{\lambda}$,
and we will use a natural projection 
$\pi _E:T({\bf R}\times M)\to E$.

Let $\{\varphi _t\}_{0\leq t\leq T}$ be the isotopy on $M$
such that
\[
\left\{
\begin{array}{rcl}
\frac{d}{dt}\varphi _t&=&X_{\lambda }\circ \varphi _t ,\\
\varphi _0&=& \mbox{id}.
\end{array}
\right.
\]
From the definition, we can conclude that
$\frac{d}{dt}(\varphi _t^*\lambda )=0$ and
$\frac{d}{dt}(\varphi _t^*d\lambda )=0$,
and then 
$d\varphi _T(\xi _p)=\xi _{\varphi _T(p)}$ and 
$d\varphi _TX_{\lambda }(p)=X_{\lambda }(\varphi _T(p))$
for $p\in M$.

We call a map $\gamma :{\bf R}/T{\bf Z}\to M$ such that
$\dot{\gamma }=X_{\lambda }(\gamma )$ 
a {\it closed characteristic} of period $T$.
Similarly, we call a map $\gamma :[0,T]\to M$ such that
$\dot{\gamma }=X_{\lambda }(\gamma )$ with 
$\gamma (0)\in \Lambda _0$ and $\gamma (T)\in \Lambda _1$,
where $\Lambda _0$ and $\Lambda _1$ are Legendrian submanifolds,
a {\it Reeb chord} from $\Lambda _0 $ to $\Lambda _1$ 
of length $T$.
\\

The restriction of $d\lambda $ on $\xi $ is non-degenerate,
hence $d\lambda |_{\xi }$ induces complex structures $I$ on
$\xi $ such that the bilinear form
\[
d\lambda (x)(h,I(x)k),\ h,k\in \xi _x
\]
is a positive definite inner product, and then
\[
g_M(h,k):=\lambda (h)\lambda (k)
+d\lambda (\pi _{\xi }(h),I\pi _{\xi }(k)),\ h,k\in TM,
\]
gives a metric on $M$.
For $a,b\in {\bf R}$ and $k\in \xi$,
\[
\overline{I}
(a\frac{\partial }{\partial \theta }+bX_{\lambda }+k)
:=
-b\frac{\partial }{\partial \theta }+aX_{\lambda }+Ik
\]
is an almost complex structure on ${\bf R}\times M$,
and the equation of pseudo-holomorphic curve
for $\overline{u}:=(\alpha ,u):\Sigma \to {\bf R}\times M$
turns out to be
\[
\left\{
\begin{array}{rcl}
\pi _{\xi }\circ du+I\pi _{\xi }\circ du \circ j&=& 0,
\\
(u^*\lambda )\circ j&=& d\alpha.
\end{array}
\right.
\]
For pseudo-holomorphic strips 
we consider time-dependent complex structures $I_t$ on $\xi$, 
and if $X$ is a symplectic manifolds with concave or convex ends,
then we suppose that time-dependent 
almost complex structures $J_t$ on $X$, 
out side a compact set, are of the form $\overline{I_t}$. 
\\

Let $L_0$ and $L_1$ be Lagrangian submanifolds in $X$ such that
$L_0|_{(-\infty ,R_-]\times M_-}
=(-\infty ,R_-]\times \Lambda _0^-$
and 
$L_0|_{[R_+,\infty )\times M_+}
=[R_+,\infty )\times \Lambda _0^+$,
and also
$L_1|_{(-\infty ,R_-]\times M_-}
=(-\infty ,R_-]\times \Lambda _1^-$
and 
$L_1|_{[R_+,\infty )\times M_+}
=[R_+,\infty )\times \Lambda _1^+$.
We assume that $L_0$ and $L_1$ satisfy Assumption \ref{ass1.1}
(, and hence 
$\Lambda _0^-\cap \Lambda _1^-=\emptyset$ 
and $\Lambda _0^+\cap \Lambda _1^+=\emptyset$).
Moreover we assume the following condition
for $\Lambda _0^-$ and $\Lambda _1^-$.
\begin{ass}\label{ass2.2}%%%%%%%%%%%%%%%%%%%%%%%%%%%%%%%%%%%%%%%%%
(Nondegeneracy of Reeb chords)
Let $\gamma :[0,T]\to M$ be a Reeb chord from 
$\Lambda _i^-$ to $\Lambda _j^-$, $i\neq j$.
Then $d\varphi _T(T_{\gamma (0)}\Lambda _i^-)$ and
$T_{\gamma (T)}\Lambda _j^-$ 
intersect transversally in $\xi _{\gamma (T)}$,
where $\{\varphi _t\}_{0\leq t\leq T}$ is the isotopy
generated by the Reeb vector field.
\end{ass}%%%%%%%%%%%%%%%%%%%%%%%%%%%%%%%%%%%%%%%%%%%%%%%%%%%%%%%%%
From the above assumption we can conclude that 
Reeb chords are isolated.
\\

Now we will consider non-constant pseudo-holomorphic strips,
i.e., maps $\overline{u} :\Sigma ={\bf R}\times [0,1]\to X$
such that
\begin{eqnarray*}
&\overline{\partial }_{J_t}\overline{u}(\tau ,t)=0
\\
&\overline{u}({\bf R},0)\subset L_0 
\mbox{ and } \overline{u}({\bf R},1)\subset L_1,
\end{eqnarray*}
with the following asymptotic conditions.
(We are interested in concave ends,
hence for simplicity we shall use notation $M$ instead of $M_-$.)
In the following 
we denote $\overline{u}$ on concave ends $(-\infty ,R_-]\times M$
by $(\alpha ,u)$:
\\
\\
${\bf (I)}$
$\lim_{\tau \to -\infty}\overline{u}(\tau ,[0,1])=p_-$
and 
$\lim_{\tau \to \infty}\overline{u}(\tau ,[0,1])=p_+$
for $p_-$ and $p_+\in L_0\cap L_1$.
\\
\\
${\bf (II)}$
$\lim_{\tau \to \infty}\overline{u}(\tau ,[0,1])=p_+
\in L_0\cap L_1$,
and
there is a Reeb chord $\gamma _-:[0,T_-]\to M$
from $\Lambda _0^-$ to $\Lambda _1^-$
such that 
$\lim_{\tau \to -\infty }\alpha (\tau ,t)=-\infty $ and 
$\lim_{\tau \to -\infty }u(\tau ,t)=\gamma _-(T_-t)$.
\\
\\
${\bf (II')}$
$\lim_{\tau \to -\infty}\overline{u}(\tau ,[0,1])=p_-
\in L_0\cap L_1$,
and
there is a Reeb chord $\gamma _+:[0,T_+]\to M$
from $\Lambda _1^-$ to $\Lambda _0^-$
such that 
$\lim_{\tau \to \infty }\alpha (\tau ,t)=-\infty $ and 
$\lim_{\tau \to \infty }u(\tau ,t)=\gamma _+(T_+(1-t))$.
\\
\\
${\bf (III)}$
There is a Reeb chord $\gamma _-:[0,T_-]\to M$
from $\Lambda _0^-$ to $\Lambda _1^-$
such that 
$\lim_{\tau \to -\infty }\alpha (\tau ,t)=-\infty $ and 
$\lim_{\tau \to -\infty }u(\tau ,t)=\gamma _-(T_-t)$,
and 
a Reeb chord $\gamma _+:[0,T_+]\to M$
from $\Lambda _1^-$ to $\Lambda _0^-$
such that 
$\lim_{\tau \to \infty }\alpha (\tau ,t)=-\infty $ and 
$\lim_{\tau \to \infty }u(\tau ,t)=\gamma _+(T_+(1-t))$.
\\
\\
Moreover we will consider the following 
extra pseudo-holomorphic strips in the symplectization of $M$, 
i.e.,
$\overline{u}:=(\alpha ,u)
:\Sigma ={\bf R}\times [0,1]\to {\bf R}\times M$
such that
\begin{eqnarray*}
&\overline{\partial }_{\overline{I_t}}\overline{u}(\tau ,t)=0
\\
&\overline{u}({\bf R},0)\subset {\bf R}\times \Lambda _0^- 
\mbox{ and } 
\overline{u}({\bf R},1)\subset {\bf R}\times \Lambda _1^-
\end{eqnarray*}
with the asymptotic conditions:
\\
\\
${\bf (IV)}$
There is a Reeb chord $\gamma _-:[0,T_-]\to M$
from $\Lambda _1^-$ to $\Lambda _0^-$
such that 
$\lim_{\tau \to -\infty }\alpha (\tau ,t)=\infty $ and 
$\lim_{\tau \to -\infty }u(\tau ,t)=\gamma _-(T_-(1-t))$, and
a Reeb chord $\gamma _+:[0,T_+]\to M$
from $\Lambda _0^-$ to $\Lambda _1^-$
such that 
$\lim_{\tau \to \infty }\alpha (\tau ,t)=\infty $ and 
$\lim_{\tau \to \infty }u(\tau ,t)=\gamma _+(T_+t)$.
\\
\\
${\bf (V)}$
There is a Reeb chord $\gamma _-:[0,T_-]\to M$
from $\Lambda _0^-$ to $\Lambda _1^-$
such that 
$\lim_{\tau \to -\infty }\alpha (\tau ,t)=-\infty $ and 
$\lim_{\tau \to -\infty }u(\tau ,t)=\gamma _-(T_-t)$, and
a Reeb chord $\gamma _+:[0,T_+]\to M$
from $\Lambda _0^-$ to $\Lambda _1^-$
such that 
$\lim_{\tau \to \infty }\alpha (\tau ,t)=\infty $ and 
$\lim_{\tau \to \infty }u(\tau ,t)=\gamma _+(T_+t)$.
\\
\\
${\bf (V')}$
There is a Reeb chord $\gamma _-:[0,T_-]\to M$
from $\Lambda _1^-$ to $\Lambda _0^-$
such that 
$\lim_{\tau \to -\infty }\alpha (\tau ,t)=\infty $ and 
$\lim_{\tau \to -\infty }u(\tau ,t)=\gamma _-(T_-(1-t))$, and
a Reeb chord $\gamma _+:[0,T_+]\to M$
from $\Lambda _1^-$ to $\Lambda _0^-$
such that 
$\lim_{\tau \to \infty }\alpha (\tau ,t)=-\infty $ and 
$\lim_{\tau \to \infty }u(\tau ,t)=\gamma _+(T_+(1-t))$.
\\
\\
If $\gamma _-=\gamma _+=\gamma $ of length $T$, 
then there is a trivial pseudo-holomorphic strip 
$\overline{u}:=(T\tau,\gamma (Tt))$ of the type ${\bf (V)}$, 
and also $\overline{u}:=(-T\tau ,\gamma (T(1-t)))$
of the type ${\bf (V')}$. 

For the each asymptotic,
we put the following exponential decay condition:
\\
\\
${\bf (1)}$
If $\lim _{\tau \to -\infty}\overline{u}(\tau, t)
=p_-\in L_0\cap L_1$,
then there are some constants 
$\rho _-$, $C_{\beta }^-$ and $d_->0$, and a function 
$\xi_- (\tau ,t)\in C^{\infty }({\bf R}\times [0,1],T_{p_-}X)$
such that
\begin{itemize}
\item $\overline{u}(\tau ,t)
=\exp _{p_-}\xi _-(\tau ,t)$ for $\tau <\rho _-$,
\item $|\partial ^{\beta }\xi _-(\tau ,t)|
<C_{\beta}^-e^{d_-\tau }$.
\end{itemize}
${\bf (1')}$
If $\lim _{\tau \to \infty}\overline{u}(\tau, t)
=p_+\in L_0\cap L_1$,
then there are some constants 
$\rho _+$, $C_{\beta }^+$ and $d_+>0$, and a function 
$\xi_+ (\tau ,t)\in C^{\infty }({\bf R}\times [0,1],T_{p_+}X)$
such that
\begin{itemize}
\item $\overline{u}(\tau ,t)
=\exp _{p_+}\xi (\tau ,t)$ for $\rho _+<\tau $,
\item $|\partial ^{\beta }\xi _+(\tau ,t)|
<C_{\beta }^+e^{-d_+\tau }$.
\end{itemize}
Let 
$\iota :[0,1]\times \{z\in {\bf R}^{2n},|z|<\varepsilon \}\to M$
be an immersion 
such that 
$\iota (t ,0)=\gamma (Tt)$, where $\gamma $ 
is a Reeb chord of length $T$:
\\
\\
${\bf (2)}$
If $\lim_{\tau \to -\infty }\alpha (\tau ,t)=-\infty $ and 
$\lim_{\tau \to -\infty }u(\tau ,t)=\gamma _-(T_-t)$,
then there are some $\tau _-$ such that
$u(\tau ,t)\in \mbox{Im }\iota _-$ for $\tau <\tau _-$,
where $\iota _-$ is the immersion with respect to $\gamma _-$.
If we denote the $\iota _-$ pull-back of $(\alpha ,u)$ by  
\[
(\alpha ,\theta_-, z_-):(-\infty ,\tau _-]\times [0,1]
\to {\bf R}\times [0,1]\times 
\{z\in {\bf R}^{2n},|z|<\varepsilon \},\]
then there are some constants
$c_-$, $C_{\beta }^-$ and $d_->0$ such that 
\begin{eqnarray*}
|\partial ^{\beta }[\alpha (\tau ,t)-(T_-\tau +c_-)]|
&\leq &C_{\beta }^-e^{d_-\tau},
\\
|\partial ^{\beta }[\theta _-(\tau ,t)-t]|
&\leq & C_{\beta }^-e^{d_-\tau},
\\
|\partial ^{\beta }z_-(\tau ,t)|
&\leq & C_{\beta }^-e^{d_-\tau}.
\end{eqnarray*}
${\bf (2')}$
If $\lim_{\tau \to \infty }\alpha (\tau ,t)=-\infty $ and 
$\lim_{\tau \to \infty }u(\tau ,t)=\gamma _+(T_+(1-t))$,
then there are some $\tau _+$ such that
$u(\tau ,t)\in \mbox{Im }\iota _+$ for $\tau _+<\tau $,
where $\iota _+$ is the immersion with respect to $\gamma _+$.
If we denote the $\iota _+$ pull-back of $(\alpha ,u)$ by  
\[
(\alpha ,\theta_+, z_+):
[\tau _+,\infty )\times [0,1]\to 
{\bf R}\times [0,1]\times \{z\in {\bf R}^{2n},|z|<\varepsilon \},
\]
then there are some constants $c_+$, $C_{\beta }^+$ and $d_+>0$ 
such that 
\begin{eqnarray*}
|\partial ^{\beta }[\alpha (\tau ,t)-(-T_+\tau +c_+)]|
&\leq &C_{\beta }^+e^{-d_+\tau},
\\
|\partial ^{\beta }[\theta _+(\tau ,t)-(1-t)]|
&\leq & C_{\beta }^+e^{-d_+\tau},
\\
|\partial ^{\beta }z_+(\tau ,t)|
&\leq & C_{\beta }^+e^{-d_+\tau}.
\end{eqnarray*}
${\bf (3)}$
If $\lim_{\tau \to -\infty }\alpha (\tau ,t)=\infty $ and 
$\lim_{\tau \to -\infty }u(\tau ,t)=\gamma _-(T_-(1-t))$,
then there are some $\tau _-$ such that
$u(\tau ,t)\in \mbox{Im }\iota _-$ for $\tau <\tau _-$,
where $\iota _-$ is the immersion with respect to $\gamma _-$.
If we denote the $\iota _-$ pull-back of $(\alpha ,u)$ by  
\[
(\alpha ,\theta_-, z_-):
(-\infty ,\tau _-]\times [0,1]\to 
{\bf R}\times [0,1]\times \{z\in {\bf R}^{2n},|z|<\varepsilon \},
\]
then there are some constants 
$c_-$, $C_{\beta }^-$ and $d_->0$ such that 
\begin{eqnarray*}
|\partial ^{\beta }[\alpha (\tau ,t)-(-T_-\tau +c_-)]|
&\leq &C_{\beta }^-e^{d_-\tau},
\\
|\partial ^{\beta }[\theta _-(\tau ,t)-(1-t)]|
&\leq & C_{\beta }^-e^{d_-\tau},
\\
|\partial ^{\beta }z_-(\tau ,t)|
&\leq & C_{\beta }^-e^{d_-\tau}.
\end{eqnarray*}
${\bf (3')}$
If $\lim_{\tau \to \infty }\alpha (\tau ,t)=\infty $ and 
$\lim_{\tau \to \infty }u(\tau ,t)=\gamma _+(T_+t)$,
then there are some $\tau _+$ such that
$u(\tau ,t)\in \mbox{Im }\iota _+$ for $\tau _+<\tau $,
where $\iota _+$ is the immersion with respect to $\gamma _+$.
If we denote the $\iota _+$ pull-back of $(\alpha ,u)$ by  
\[
(\alpha ,\theta_+, z_+):
[\tau _+,\infty )\times [0,1]\to 
{\bf R}\times [0,1]\times \{z\in {\bf R}^{2n},|z|<\varepsilon \},
\]
then there are some constants
$c_+$, $C_{\beta }^+$ and $d_+>0$ such that 
\begin{eqnarray*}
|\partial ^{\beta }[\alpha (\tau ,t)-(T_+\tau +c_+)]|
&\leq &C_{\beta }^+e^{-d_+\tau},
\\
|\partial ^{\beta }[\theta _+(\tau ,t)-t]|
&\leq & C_{\beta }^+e^{-d_+\tau},
\\
|\partial ^{\beta }z_+(\tau ,t)|
&\leq & C_{\beta }^+e^{-d_+\tau}.
\end{eqnarray*}

Define $\mathcal{M}_I(p,q)$ to be the set of
pseudo-holomorphic strips 
of the form ${\bf (I)}$ with ${\bf (1)}$ and ${\bf (1')}$,
also $\mathcal{M}_{II}(\gamma _-,p_+)$
of the form ${\bf (II)}$ with ${\bf (2)}$ and ${\bf (1')}$,
$\mathcal{M}_{II'}(p_-,\gamma _+)$
of the form ${\bf (II')}$ with ${\bf (1)}$ and ${\bf (2')}$,
$\mathcal{M}_{III}(\gamma _-,\gamma _+)$
of the form ${\bf (III)}$ with ${\bf (2)}$ and ${\bf (2')}$,
$\mathcal{M}_{IV}(\gamma _-,\gamma _+)$
of the form ${\bf (IV)}$ with ${\bf (3)}$ and ${\bf (3')}$,
$\mathcal{M}_{V}(\gamma _-,\gamma _+)$
of the form ${\bf (V)}$ with ${\bf (2)}$ and ${\bf (3')}$
and finally
$\mathcal{M}_{V'}(\gamma _-,\gamma _+)$
of the form ${\bf (V')}$ with ${\bf (3)}$ and ${\bf (2')}$.
Note that ${\bf R}$ acts on the moduli spaces  of the type 
${\bf (I)}$, ${\bf (II)}$, ${\bf (II')}$ and ${\bf (III)}$
by $\frac{d}{da}(a*\overline{u})$, $a\in {\bf R}$.
On the other hand, for the moduli spaces of the type 
${\bf (IV)}$, ${\bf (V)}$ and ${\bf (V')}$, 
${\bf R}^2$ acts on them by
$\frac{d}{da}(a*\overline{u})$ for $(a,0)\in {\bf R}^2$ and
$\frac{d}{db}(b\sharp \overline{u})$ for $(0,b)\in {\bf R}^2$,
where $b\sharp \overline{u}:=(\alpha -b,u)$ for
$\overline{u}:=(\alpha ,u)$.
We shall denote these quotients by $\hat{\mathcal{M}}_*(*,*)$.
\\

First we recall the index for strips of the type ${\bf (I)}$
with ${\bf (1)}$ and ${\bf (1')}$.
(This part is so standard, 
the reader may skip to the next content.)
Choose a trivialization
$\{e^-_1,e^-_2,\ldots ,e^-_{2n+1},e^-_{2n+2}\}$ of $T_{p_-}X$
such that 
$\omega _{p_-}(e^-_{2i-1},e^-_{2j})=\delta _{ij}$,
and similarly 
$\{e^+_1,e^+_2,\ldots ,e^+_{2n+1},e^+_{2n+2}\}$ of $T_{p_+}X$
such that $\omega _{p_+}(e^+_{2i-1},e^+_{2j})=\delta _{ij}$.
We assume that our time-dependent almost complex structures $J_t$
satisfy the following condition.
\begin{ass}\label{ass2.3}%%%%%%%%%%%%%%%%%%%%%%%%%%%%%%%%%%%%%%%%%
A time-dependent almost complex structure $J_t$, $t\in [0,1]$,
satisfies that $J_t$ on $T_{p_-}X$ is standard with respect to
$\{e^-_1,e^-_2,\ldots ,e^-_{2n+1},e^-_{2n+2}\}$,
i.e.,
$J_te^-_{2i-1}=e^-_{2i}$ and $J_te^-_{2i}=-e^-_{2i-1}$,
and similarly
$J_t$ on $T_{p_+}X$ is standard with respect to
$\{e^+_1,e^+_2,\ldots ,e^+_{2n+1},e^+_{2n+2}\}$.
\end{ass}%%%%%%%%%%%%%%%%%%%%%%%%%%%%%%%%%%%%%%%%%%%%%%%%%%%%%%%%%
(We can always choose such almost complex structures.)
Let $g(t)$, $t\in [0,1]$, be a time-dependent metric on $X$ 
such that 
$L_0$ is totally geodesic with respect to $g(0)$
and similarly $L_1$ is totally geodesic with respect to $g(1)$.
We denote by $\exp :=\exp (t) :TX\to X$ 
the exponential map.
Let $\overline{u}\in L^p_1({\bf R}\times [0,1];X,L_0,L_1)$ 
be a map satisfying the Lagrangian boundary conditions and 
${\bf (I)}$ with the decay conditions ${\bf (1)}$ 
and ${\bf (1')}$.
($\overline{u}$ need not be pseudo-holomorphic.)
For $\eta \in L^p_1(\overline{u}^*TX, \overline{u}^*TL_0,
\overline{u}^*TL_1)$, a section of $\overline{u}^*TX$
with $\eta (\tau ,0)\in T_{\overline{u}(\tau ,0)}L_0$
and $\eta (\tau ,1)\in T_{\overline{u}(\tau ,1)}L_1$,
we define a map 
$f_{\overline{u}}:
L^p_1(\overline{u}^*TX, \overline{u}^*TL_0,\overline{u}^*TL_1)\to
L^p(\overline{u}^*TX\otimes 
\wedge ^{0,1}T^*({\bf R}\times [0,1]))$
by
\begin{eqnarray}
\label{2}
f_{\overline{u}}(\eta ):=
\Phi_{\overline{u}}(\eta )^{-1}
\overline{\partial }_{J_t}(\exp _{\overline{u}}(\eta )), 
\end{eqnarray}
where 
$\Phi_{\overline{u}}(\eta ):T_{\overline{u}}X\to 
T_{\exp _{\overline{u}}}X$ 
denotes parallel transport of a connection
along the geodesic $t\to \exp _{\overline{u}}(t\eta )$. 
The differential $Df_{\overline{u}}(0)$ is
\[
Df_{\overline{u}}(0)\eta =\nabla _{\partial /\partial \tau}\eta
+J_t(\overline{u}(\tau ,t))\nabla _{\partial /\partial t}\eta
+(\nabla _{\eta }J_t)\partial _t\overline{u}.
\]
As $\tau \to -\infty $, the right hand side is 
\begin{equation}\label{3}
\nabla _{\partial /\partial \tau}\eta
+J_t(p_-)\nabla _{\partial /\partial t}\eta ,
\end{equation}
and if we denote $\eta =\sum_{i=1}^{2n+2}\eta _ie^-_i$, 
then an equation $(\ref{3})=0$ turns out to be
\[
\frac{\partial }{\partial \tau}
\left[
\begin{array}{c}
\eta _1
\\
\vdots
\\ 
\eta _{2n+2}
\end{array}
\right]
+
J_0
\frac{\partial }{\partial t}
\left[
\begin{array}{c}
\eta _1
\\
\vdots
\\ \eta _{2n+2}
\end{array}
\right]
=0.
\]
For simplicity we shall use $E_{\overline{u}}$ 
to denote $Df_{\overline{u}}(0)$ 
and put $E_{\overline{u}}=\frac{d}{d\tau }-Q_t$.
Then we conclude that 
$Q_{-\infty }$ has no eigenvectors of eigenvalue $0$
from Assumption \ref{ass1.1}.
Similarly,
as $\tau \to \infty $, $Q_{\infty }$ also
has no eigenvectors of eigenvalue $0$.
Since $Q_{-\infty}$ and $Q_{\infty }$ are invertible,
$E_{\overline{u}}:
L^p_1(\overline{u}^*TX, \overline{u}^*TL_0,\overline{u}^*TL_1)\to
L^p(\overline{u}^*TX\otimes 
\wedge ^{0,1}T^*({\bf R}\times [0,1]))$ 
is Fredholm.
$\mbox{Ind}E_{\overline{u}}$ denotes 
the index of $E_{\overline{u}}$.
\\

Similarly we will introduce an index for strips of the type 
${\bf (II)}$ with ${\bf (2)}$ and ${\bf (1')}$.
For Reeb chords $\gamma :[0,T]\to M$
we define $\overline{\gamma }(t):=\gamma (Tt)$.
Consider the pull-back $\overline{\gamma }^*\xi $ over 
$\mathcal{I}:=[0,1]$
and choose a trivialization
$\{e_1,e_2,\ldots ,e_{2n}\}$
such that 
$e_i(t)=d\varphi _{Tt}e_i(0)$, $t\in [0,1]$,
and 
$\overline{\gamma }^*d\lambda (e_{2i-1}(0),e_{2j}(0))
=\delta _{ij}$.
Take a time-dependent connection 
$\nabla ^{\lambda }:=\nabla ^{\lambda }(t)$,
$t\in [0,1]$, on $\xi $ 
so that 
the holonomy of $\overline{\gamma }^*\nabla ^{\lambda }(t)$
agrees with $d\varphi _{Tt}$ along the Reeb chords,
i.e.,
$\overline{\gamma }^*\nabla ^{\lambda }(t)e_i(t)=0$.
Let $\{f_1,f_2,\ldots ,f_{2n}\}$ be another trivialization 
such that $\overline{\gamma }^*I_t(t)$ is 
the standard complex structure $J_0$
with respect to $\{f_1(t),f_2(t),\ldots ,f_{2n}(t)\}$,
i.e.,
$\overline{\gamma }^*I_t(t)f_{2i-1}(t)=f_{2i}(t)$ 
and $\overline{\gamma }^*I_t(t)f_{2i-1}(t)=-f_{2i}(t)$.
If $e_i(t)=\sum_{j=1}^{2n}a_{ij}(t)f_j(t)$ and $A:=[a_{ij}(t)]$,
then
\[
\overline{\gamma }^*I_t\overline{\gamma }^*
\nabla ^{\lambda }_{\frac{\partial }{\partial t}}
\sum_{i=1}^{2n}\eta _i(t)f_i(t)
=
[f_1 f_2 \cdots f_{2n}]
\left\{
J_0\frac{\partial }{\partial t}
\left[
\begin{array}{c}
\eta _1
\\
\vdots 
\\
\eta _{2n}
\end{array}
\right]
-J_0\frac{\partial A}{\partial t}
A^{-1}
\left[
\begin{array}{c}
\eta _1
\\
\vdots 
\\
\eta _{2n}
\end{array}
\right]
\right\}.
\]
We conclude that $-J_0\frac{\partial A}{\partial t}A^{-1}$
is symmetric 
if the complex structures satisfy the following condition.
\begin{ass}\label{ass2.4}%%%%%%%%%%%%%%%%%%%%%%%%%%%%%%%%%%%%%%%%%
A time-dependent complex structure $I_{t}$, $t\in [0,1]$,
on $\xi $ satisfies that
$\overline{\gamma }^*d\lambda $ is standard with respect to 
$\{f_1(t),f_2(t),\ldots ,f_{2n}(t)\}$, i.e.,
$\overline{\gamma }^*d\lambda (f_{2i-1}(t),f_{2j}(t))
=\delta _{ij}.$ 
\end{ass}%%%%%%%%%%%%%%%%%%%%%%%%%%%%%%%%%%%%%%%%%%%%%%%%%%%%%%%%%
(We can always choose 
$\nabla ^{\lambda }$ and $I_t$ as above.
In fact, if we choose them so that $e_i(t)=f_i(t)$, 
then $a_{ij}(t)=\delta _{ij}$.)
The double $\mathcal{I}'$ of $\mathcal{I}$ is the circle 
$\mathcal{I}'=\mathcal{I}\cup _{\partial \mathcal{I}}
\overline{\mathcal{I}}$
where $\overline{\mathcal{I}}$ is the mirror image of 
$\mathcal{I}$ and corresponding points 
on the boundary are identified.
We denote $\kappa $ a natural involution. 
The doubling of $\overline{\gamma }^*\xi $ is completely similar.
Let $\overline{\overline{\gamma }^*\xi }$ denote the vector bundle
over $\overline{\mathcal{I}}$ whose fiber 
$(\overline{\overline{\gamma }^*\xi })_{\kappa (t)}$ is 
$(\overline{\gamma }^*\xi )_t$
with the complex structure $-\overline{\gamma }^*I_t(t)$. 
Then the double $\overline{\gamma }^*\xi '$ over $\mathcal{I}'$
is obtained by gluing $\overline{\gamma }^*\xi $ 
and $\overline{\overline{\gamma }^*\xi}$ 
along $\partial \mathcal{I}$.
The identification is realized by the automorphism $s$
of $\overline{\gamma }^*\xi|_{\partial \mathcal{I}}$ 
which is the anti-complex reflection through 
$T_{\gamma (0)}\Lambda _0$ over $\gamma (0)$
and through $T_{\gamma (T)}\Lambda _1$ over $\gamma (T)$.
We denote $\overline{\gamma }^*I_t'$
the double of the complex structure on $\overline{\gamma }^*\xi '$
which is $\overline{\gamma }^*I_t$ on $\overline{\gamma }^*\xi $
and $-\overline{\gamma }^*I_t$ on 
$\overline{\overline{\gamma }^*\xi }$.
An element $\eta '\in L^p_1(\overline{\gamma }^*\xi')$ 
is defined by 
$\eta '|_{\mathcal{I}}=\eta _1\in L^p_1(\overline{\gamma }^*\xi)$
and 
$\tilde{\kappa }\circ (\eta '|_{\overline{\mathcal{I}}})
\circ \kappa  =\eta _2\in L^p_1 (\overline{\gamma }^*\xi)$,
where $\tilde{\kappa }$ is the natural involution 
lifting $\kappa$.
These satisfy 
$\eta _2|_{\partial \mathcal{I}}=s\circ (\eta _1
|_{\partial \mathcal{I}})$
which is equivalent to
\begin{eqnarray*}
\eta _1+\eta _2
& \in & 
L^p_1(\overline{\gamma }^*\xi ;T_{\gamma (0)}\Lambda _0,
T_{\gamma (T)}\Lambda _1),
\\
\overline{\gamma }^*I_t(\eta _1-\eta _2)
& \in & 
L^p_1(\overline{\gamma }^*\xi ;T_{\gamma (0)}\Lambda _0,
T_{\gamma (T)}\Lambda _1).
\end{eqnarray*}
Conversely, any couple $(\eta _1,\eta _2)\in 
L^p_1 (\overline{\gamma }^*\xi)\times 
L^p_1 (\overline{\gamma }^*\xi)$
satisfying the above conditions defines 
$\eta '\in L^p_1(\overline{\gamma }^*\xi')$.
We can now define the double operator 
$\overline{\gamma }^*\nabla ^{\lambda }(t)'$ by
$\overline{\gamma }^*\nabla ^{\lambda }(t)$
on $\overline{\gamma }^*\xi$
and 
$\tilde{\kappa }\circ \overline{\gamma }^*\nabla ^{\lambda }(t)
\circ \kappa $ on $\overline{\overline{\gamma }^*\xi}$.
We conclude that the holonomy around the circle $\mathcal{I}'$
is equal to $-\mbox{id}$ from Assumption \ref{ass2.2},
and then the equation for sections $\eta '(t)\in 
L^p_1(\overline{\gamma}^*\xi ')$
\[
\overline{\gamma }^*I_t'\overline{\gamma }^*
\nabla ^{\lambda }_{\frac{\partial }{\partial t}}(t)'\eta '(t)=0
\]
has no eigen functions of eigen value $0$.
\\

Let $\nabla :=\nabla (t)$ be a time-dependent connection 
on $TX$ which, out side a compact set, is of the form:
the restriction on $E$ is trivial, i.e.,
$\nabla (a\frac{\partial }{\partial \theta }+bX_{\lambda })
=da\otimes \frac{\partial }{\partial \theta }
+db\otimes X_{\lambda }$,
and the restriction on $\xi $
is the pull-back of $\nabla ^{\lambda }$.
Similarly, a time-dependent almost complex structure $J_t$ on $X$,
out side a compact set, is of the form $\overline{I_t}$.
Let $\overline{u}$ be a map
which satisfies the Lagrangian boundary conditions and 
${\bf (II)}$ and the decay condition ${\bf (2)}$.
($\overline{u}$ need not be pseudo-holomorphic.)
Then, the differential $Df_{\overline{u}}(0)$ is
\[
Df_{\overline{u}}(0)\eta =\nabla _{\partial /\partial \tau}\eta
+J(\overline{u}(\tau ,t))\nabla _{\partial /\partial t}\eta
+(\nabla _{\eta }J)\partial _t\overline{u}.
\]
As $\tau \to -\infty $, the right hand side is 
\begin{equation}\label{4}
\nabla _{\partial /\partial \tau}\eta
+\overline{I_t}(\gamma _-(T_-t))
\nabla _{\partial /\partial t}\eta ,
\end{equation}
and if we denote
$\eta =\eta _1\frac{\partial }{\partial \theta}
+\eta _2X_{\lambda }+\sum_{i=3}^{2n+2}\eta _if_i$, 
then an equation $(\ref{4})=0$ splits into
\[
\frac{\partial }{\partial \tau}
\left[
\begin{array}{c}
\eta _1 \\ \eta _2
\end{array}
\right]
+
\left[
\begin{array}{cc}
0 & -1 \\ 
1 & 0
\end{array}
\right]
\frac{\partial }{\partial t}
\left[
\begin{array}{c}
\eta _1 \\ \eta _2
\end{array}
\right]
=0,
\]
\[
\frac{\partial }{\partial \tau}
\left[
\begin{array}{c}
\eta _3
\\
\vdots
\\ 
\eta _{2n+2}
\end{array}
\right]
+
J_0
\frac{\partial }{\partial t}
\left[
\begin{array}{c}
\eta _3 
\\
\vdots
\\ \eta _{2n+2}
\end{array}
\right]
-J_0\frac{\partial A}{\partial t}A^{-1}
\left[
\begin{array}{c}
\eta _3 
\\
\vdots
\\ \eta _{2n+2}
\end{array}
\right]
=0.
\]
If we consider the double operator of $E_{\overline{u}}$ on 
${\bf R}\times S^1:=
({\bf R}\times [0,1])\cup _{\partial ({\bf R}\times [0,1])}
({\bf R}\times [0,1])$
(the double is exactly similar to that in the last paragraph),
then $Q_{-\infty }$ has eigenvectors 
$(\eta _1,\eta _2,\eta _3,\ldots ,\eta _{2n+2})=
(1,0,0,\ldots ,0)$ and $(0,1,0,\ldots ,0)$ of eigenvalue $0$.
(If we simply consider $E_{\overline{u}}$ 
on ${\bf R}\times [0,1]$,
then $Q_{-\infty}$ on $\mathcal{I}$ has an eigenvector 
$(1,0,\ldots , 0)$ of eigenvalue $0$
under the Lagrangian boundary conditions.)
\\

Now we introduce weighted Sobolev spaces for
the Fredholm theory of $E_{\overline{u}}$.
For $0<\delta <2\pi $ 
and $\tau _-$ as in the decay condition ${\bf (2)}$, 
we define a smooth decreasing function $\sigma $ by
\[
\sigma (\tau ):=
\left\{
\begin{array}{ll}
-\delta (\tau -\tau _-+1), & \tau \leq \tau _--2,
\\
0, & \tau _-\leq \tau ,
\end{array}
\right.
\]
and a cut function $\beta :{\bf R}\to [0,1]$ by
\[
\beta (\tau ):=
\left\{
\begin{array}{ll}
1, & \tau \leq \tau _--1,
\\
0, & \tau _-\leq \tau .
\end{array}
\right.
\]
For a section $\eta $ of $\overline{u}^*TX|
_{(-\infty ,\tau _-]\times M}$,
we denote 
$\eta =\eta _E+\eta _{\xi}$, 
where $\eta _E$ is the $E$ component and 
$\eta _{\xi }$ is the $\xi $ component.
Then we define weighted Sobolev norms by
\begin{eqnarray*}
\|\eta \|_{L^p_{1;\sigma }}
&:= &
\left\{
\int _{[\tau _-,\infty)\times [0,1]}
(1-\beta (\tau))
(|\eta |^p+|\nabla \eta |^p)d\tau dt
\right.
\\
&&
+\left.
\int _{(-\infty ,\tau _-]\times [0,1]}
e^{\sigma (\tau )}\beta (\tau )
(|\eta _E|^p+|\nabla \eta _E|^p)
+\beta (\tau )(|\eta _{\xi }|^p+|\nabla \eta _{\xi }|^p)
d\tau dt
\right\}^{1/p},
\end{eqnarray*}
and similarly 
\begin{eqnarray*}
\|\eta \|_{L^p_{0;\sigma }}:= 
\left\{
\int _{[\tau _-,\infty)\times [0,1]}
(1-\beta (\tau))|\eta |^pd\tau dt
+\int _{(-\infty ,\tau _-]\times [0,1]}
e^{\sigma (\tau )}\beta (\tau )|\eta _E|^p
+\beta (\tau )|\eta _{\xi }|^pd\tau dt
\right\}^{1/p}.
\end{eqnarray*}
We define $L^p_{1;\sigma }$ to be the set of sections $f$ 
such that $\|f\|_{L^p_{1;\sigma }}<\infty $, 
and also $L^p_{0;\sigma }$.
Let $I:L^2\to L^2_{0;\sigma }$ be the following 
isometric transformation.
\[
I(\eta )(\tau ,t)
:=\left\{
\begin{array}{ll}
e^{-\sigma (\tau )/2}\eta _E+\eta _{\xi },& \tau \leq \tau _-
\\
\eta ,& \tau _-\leq \tau
\end{array}
\right.
\]
Then we obtain
\[
E_{\overline{u};\sigma }(\eta )(\tau ,t)
:=I^{-1}E_{\overline{u}}I(\eta )(\tau ,t)
=(E_{\overline{u}}\eta )(\tau ,t)
-\frac{1}{2}\frac{d\sigma }{d\tau}\eta _E.
\]
As $\tau \to -\infty $,
$E_{\overline{u};\sigma }$ is equal to 
\[
\frac{\partial }{\partial \tau}-
Q_{-\infty }+\frac{\delta }{2}\pi _E.
\]
Since $-Q_{-\infty }+\frac{\delta }{2}\pi _E$ is invertible
(and $Q_{\infty }$ is also invertible, see the last paragraph),
$E_{\overline{u};\sigma }:L^2_1\to L^2$ is a Fredholm operator,
which implies that 
$E_{\overline{u}}:L^2_{1;\sigma }\to L^2_{0;\sigma }$ is Fredholm.
$\mbox{Ind}E_{\overline{u};\sigma }$ 
denotes the index of $E_{\overline{u};\sigma }$.
Note that,
if $\overline{u}$ is pseudo-holomorphic,
then $E_{\overline{u}}\frac{d}{da}(a*\overline{u})|_{a=0}=0$.
But $\frac{d}{da}(a*\overline{u})|_{a=0}$
is not an element of $L^p_{1;\sigma }$,
because $\frac{d}{da}(a*\overline{u})|_{a=0}$,
as $\tau \to -\infty $, is close to 
$-T_-(\frac{\partial }{\partial \alpha })$
which is not in $L^p_{1;\sigma }$. 
On the other hand
a map $E_{\overline{u}}:
\langle \frac{d}{da}(a*\overline{u})|_{a=0}\rangle
\to L^p_{0;\sigma }$ makes sense, hence we use the same notation 
$E_{\overline{u}}:
\langle \frac{d}{da}(a*\overline{u})|_{a=0}\rangle \oplus 
L^p_{1;\sigma }
(\overline{u}^*TX,\overline{u}^*TL_0,\overline{u}^*TL_1)
\to L^p_{0;\sigma }(\overline{u}^*TX)$.
(The tangent space of $\mathcal{M}_{II}(\gamma _-,p_+)$
at a smooth point $\overline{u}$ is 
$\mbox{Ker}\{E_{\overline{u}}:
\langle \frac{d}{da}(a*\overline{u})|_{a=0}\rangle \oplus 
L^p_{1;\sigma }
(\overline{u}^*TX,\overline{u}^*TL_0,\overline{u}^*TL_1)
\to L^p_{0;\sigma }(\overline{u}^*TX)$.)
$\mbox{Ind}E_{\overline{u}}$ denotes the index of 
$E_{\overline{u}}:
\langle \frac{d}{da}(a*\overline{u})|_{a=0}\rangle \oplus 
L^p_{1;\sigma }
(\overline{u}^*TX,\overline{u}^*TL_0,\overline{u}^*TL_1)
\to L^p_{0;\sigma }(\overline{u}^*TX)$.
(Hence $\mbox{Ind}E_{\overline{u}}
=\mbox{Ind}E_{\overline{u};\sigma }+1$.)
\\

Similarly we can introduce weighted Sobolev spaces
and indexes for strips of the type 
${\bf (II')}$ and ${\bf (III)}$.
\\

For strips of the type ${\bf (IV)}$
we can also define weighted Sobolev spaces and indexes 
$\mbox{Ind}E_{\overline{u};\sigma}$
for 
$E_{\overline{u};\sigma }:L^p_1\to L^p$.
Note that, if $\overline{u}:=(\alpha ,u)$ 
is pseudo-holomorphic strip of the type ${\bf (IV)}$, 
then $E_{\overline{u}}\frac{d}{da}(a*\overline{u})|_{a=0}=0$.
But $\frac{d}{da}(a*\overline{u})|_{a=0}$
is not an element of $L^p_{1;\sigma }$,
because $\frac{d}{da}(a*\overline{u})|_{a=0}$,
as $\tau \to -\infty $, is close to 
$T_-(\frac{\partial }{\partial \alpha })$
and, as $\tau \to \infty $, is close to 
$-T_+(\frac{\partial }{\partial \alpha})$
which is not in $L^p_{1;\sigma }$. 
Also $E_{\overline{u}}
\frac{d}{db}(b\sharp \overline{u})|_{b=0}=0$.
But $\frac{d}{db}(b\sharp \overline{u})|_{b=0}$
is also not an element of $L^p_{1;\sigma }$,
because $\frac{d}{db}(b\sharp \overline{u})|_{b=0}$,
as $\tau \to -\infty $, is close to 
$-(\frac{\partial }{\partial \alpha })$
and, $\tau \to \infty$, is also close to 
$-(\frac{\partial }{\partial \alpha})$
which is not in $L^p_{1;\sigma }$. 
On the other hand a map
$E_{\overline{u}}:\langle 
\frac{d}{da}(a*\overline{u})|_{a=0},
\frac{d}{db}(b\sharp \overline{u})|_{b=0}\rangle 
\to L^p_{0;\sigma }$
makes sense,
hence we use the same notation
$E_{\overline{u}}:\langle 
\frac{d}{da}(a*\overline{u})|_{a=0},
\frac{d}{db}(b\sharp \overline{u})|_{b=0}\rangle 
\oplus L^p_{1;\sigma }(\overline{u}^*TX, \overline{u}^*TL_0,
\overline{u}^*TL_1)
\to L^p_{0;\sigma }(\overline{u}^*TX)$.
(The tangent space of $\mathcal{M}_{IV}(\gamma _-,\gamma _+)$ 
at a smooth point $\overline{u}$ is
$\mbox{Ker}\{E_{\overline{u}}:
\langle \frac{d}{da}(a*\overline{u})|_{a=0},
\frac{d}{db}(b\sharp \overline{u})|_{b=0}\rangle 
\oplus L^p_{1;\sigma }(\overline{u}^*TX, \overline{u}^*TL_0,
\overline{u}^*TL_1)
\to L^p_{0;\sigma }(\overline{u}^*TX)\}$.)
$\mbox{Ind}E_{\overline{u}}$ denotes the index of 
$E_{\overline{u}}:
\langle \frac{d}{da}(a*\overline{u})|_{a=0},
\frac{d}{db}(b\sharp \overline{u})|_{b=0}\rangle 
\oplus L^p_{1;\sigma }
(\overline{u}^*TX, \overline{u}^*TL_0, \overline{u}^*TL_1)
\to L^p_{0;\sigma }(\overline{u}^*TX)$.
(Then $\mbox{Ind}E_{\overline{u}}=
\mbox{Ind}E_{\overline{u};\sigma }+2$.)

%%%%%%%%%%%%%%%%%%%%%%%%%%%%%%%%%%%%%%%%%%%%%%%%%%%%%%%%%%%%%%%%%%
\section{Floer's chain complexes for Lagrangian submanifolds %%%%%
in symplectic manifolds with concave/convex ends}%%%%%%%%%%%%%%%%%
%%%%%%%%%%%%%%%%%%%%%%%%%%%%%%%%%%%%%%%%%%%%%%%%%%%%%%%%%%%%%%%%%%

We will propose Floer's chain complexes for Lagrangian
submanifolds in symplectic manifolds with concave/convex ends.

Let $X$ be a symplectic manifold with 
finitely many concave or convex ends.
We shall assume that
the contact manifolds are compact without boundaries.
Let $L_0$ and $L_1$ be Lagrangian submanifolds in $X$ 
which satisfy 
Assumption \ref{ass1.1}, \ref{ass2.1} and \ref{ass2.2}.
\begin{ass}\label{ass3.1}%%%%%%%%%%%%%%%%%%%%%%%%%%%%%%%%%%%%%%%%%
For pseudo-holomorphic strips $\overline{u}$ 
of the type ${\bf (I)}$
the linear operators 
$E_{\overline{u}}:
L^p_1(\overline{u}^*TX,\overline{u}^*TL_0,\overline{u}^*TL_1)\to 
L^p(\overline{u}^*TX)$
are surjective, 
and for pseudo-holomorphic strips $\overline{u}$ 
of the type ${\bf (II)}$, ${\bf (II')}$ and ${\bf (III)}$ 
the linear operators
$E_{\overline{u}}:
\langle \frac{d}{da}(a*\overline{u})|_{a=0}\rangle 
\oplus L^p_{1;\sigma }
(\overline{u}^*TX,\overline{u}^*TL_0,\overline{u}^*TL_1)\to 
L^p_{0;\sigma }(\overline{u}^*TX)$
are also surjective,
and
pseudo-holomorphic strips $\overline{u}$ 
of the type ${\bf (IV)}$, ${\bf (V)}$ and ${\bf (V')}$ 
the linear operators
$E_{\overline{u}}:
\langle \frac{d}{da}(a*\overline{u})|_{a=0},
\frac{d}{db}(b\sharp \overline{u})|_{b=0}
\rangle 
\oplus L^p_{1;\sigma }
(\overline{u}^*TX,\overline{u}^*TL_0,\overline{u}^*TL_1)\to 
L^p_{0;\sigma }(\overline{u}^*TX)$
are surjective.
\end{ass}%%%%%%%%%%%%%%%%%%%%%%%%%%%%%%%%%%%%%%%%%%%%%%%%%%%%%%%%%
(The surjectivity or transversality problem will be observed 
in a forthcoming paper \cite{a2}.)
{}From the above assumption we can conclude that 
the moduli spaces $\mathcal{M}_*(*,*)$ are smooth manifolds 
whose dimension at $\overline{u}$ is equal to 
$\mbox{Ind}E_{\overline{u}}$,
and $\hat{\mathcal{M}}_*(*,*)$ are also
smooth manifolds whose dimension at $\overline{u}$
is equal to $\mbox{Ind}E_{\overline{u}}-1$ for the type 
${\bf (I)}$, ${\bf (II)}$, ${\bf (II')}$ and ${\bf(III)}$,
and $\mbox{Ind}E_{\overline{u}}-2$ 
for the type ${\bf (IV)}$, ${\bf (V)}$ and ${\bf (V')}$.
We denote a subset of $\mathcal{M}_*(*,*)$ or
$\hat{\mathcal{M}}_*(*,*)$ 
consists of the pseudo-holomorphic strips whose dimension is $d$
by $\mathcal{M}^d_*(*,*)$ or $\hat{\mathcal{M}}^d_*(*,*)$,
respectively.

For appropriate compactifications of moduli spaces
we need Assumption \ref{ass1.2}, \ref{ass1.4} and the following.
\begin{ass}\label{ass3.2}%%%%%%%%%%%%%%%%%%%%%%%%%%%%%%%%%%%%%%%%%
There are no contractible closed characteristics and
contractible Reeb chords from a Legendrian submanifold to itself
in the contact manifolds of concave ands.
\end{ass}%%%%%%%%%%%%%%%%%%%%%%%%%%%%%%%%%%%%%%%%%%%%%%%%%%%%%%%%%
A closed characteristic in $M$ is contractible iff
it represents $0$ in $\pi _1(M)$, 
and similarly a Reeb chord from $\Lambda $ to $\Lambda $
is contractible iff it represents $0$ of $\pi _1(M,\Lambda )$.

Moreover, for very technical reasons for the exponential decay
conditions, we need the following assumption.
\begin{ass}\label{ass3.3}%%%%%%%%%%%%%%%%%%%%%%%%%%%%%%%%%%%%%%%%
There is an open neighborhood $U\subset [0,1]\times {\bf R}^{2n}$
of $[0,1]\times \{0\}$ and an open neighborhood $V\subset M$
of a Reeb chord $x(Tt)$ of length $T$ 
and an immersion $\varphi :U\to V$
mapping $[0,1]\times \{0\}$ onto $x(Tt)$ such that
$\varphi ^*\lambda =f\lambda _0$
with $\lambda _0:=dt+\sum_{i=1}^nx_ndy_n $ and a positive smooth
function $f:U\to {\bf R}$ satisfying
$f(t,0,0)=T$ and $df(t,0,0)=0$ for all $t\in [0,1]$.
\end{ass}%%%%%%%%%%%%%%%%%%%%%%%%%%%%%%%%%%%%%%%%%%%%%%%%%%%%%%%%%
{}From these assumptions we can conclude that.
\begin{theo}\label{theo3.4}%%%%%%%%%%%%%%%%%%%%%%%%%%%%%%%%%%%%%%%
$\hat{\mathcal{M}}^0_*(*,*)$ is compact,
and $\hat{\mathcal{M}}^1_*(*,*)$
can be compactified whose boundaries are:
\begin{eqnarray*}%%%%%%%%%%%%%%%%%%%%%%%%%%%%%%%%%%%%%%%%%%%%
{\bf (a)}\ 
\partial \hat{\mathcal{M}}^1_I(p_-,p_+)
&=&
\bigcup_{q\in L_0\cap L_1} 
\hat{\mathcal{M}}^0_I(p_-,q)\times \hat{\mathcal{M}}^0_I(q,p_+)
\\
&&
\cup 
\bigcup_{(\gamma _-,\gamma _+)} 
\hat{\mathcal{M}}^0_{II'}(p_-,\gamma _-)
\times 
\hat{\mathcal{M}}^0_{IV}(\gamma _-,\gamma _+)
\times
\hat{\mathcal{M}}^0_{II}(\gamma _+,p_+),
\end{eqnarray*}%%%%%%%%%%%%%%%%%%%%%%%%%%%%%%%%%%%%%%%%%%%%%%
\begin{eqnarray*}%%%%%%%%%%%%%%%%%%%%%%%%%%%%%%%%%%%%%%%%%%%%
{\bf (b)}\
\partial \hat{\mathcal{M}}^1_{II}(\gamma _-,p_+)
&=&
\bigcup_{q\in L_0\cap L_1} 
\hat{\mathcal{M}}^0_{II}(\gamma _-,q)
\times \hat{\mathcal{M}}^0_I(q,p_+)
\\
&&
\cup 
\bigcup_{(\gamma '_-,\gamma '_+)} 
\hat{\mathcal{M}}^0_{III}(\gamma _-,\gamma '_-)
\times 
\hat{\mathcal{M}}^0_{IV}(\gamma '_-,\gamma '_+)
\times
\hat{\mathcal{M}}^0_{II}(\gamma '_+,p_+)
\\
&&
\cup 
\bigcup_{\gamma } 
\hat{\mathcal{M}}^0_{V}(\gamma _-,\gamma )
\times 
\hat{\mathcal{M}}^0_{II}(\gamma ,p_+),
\end{eqnarray*}%%%%%%%%%%%%%%%%%%%%%%%%%%%%%%%%%%%%%%%%%%%%%%
\begin{eqnarray*}%%%%%%%%%%%%%%%%%%%%%%%%%%%%%%%%%%%%%%%%%%%%
{\bf (b')}\
\partial \hat{\mathcal{M}}^1_{II'}(p_-,\gamma _+)
&=&
\bigcup_{q\in L_0\cap L_1} 
\hat{\mathcal{M}}^0_{I}(p_-,q)
\times \hat{\mathcal{M}}^0_{II'}(q,\gamma _+)
\\
&&
\cup 
\bigcup_{(\gamma '_-,\gamma '_+)} 
\hat{\mathcal{M}}^0_{II'}(p_-,\gamma '_-)
\times 
\hat{\mathcal{M}}^0_{IV}(\gamma '_-,\gamma '_+)
\times
\hat{\mathcal{M}}^0_{III}(\gamma '_+,\gamma _+)
\\
&&
\cup 
\bigcup_{\gamma } 
\hat{\mathcal{M}}^0_{II'}(p_-,\gamma )
\times 
\hat{\mathcal{M}}^0_{V'}(\gamma ,\gamma _+),
\end{eqnarray*}%%%%%%%%%%%%%%%%%%%%%%%%%%%%%%%%%%%%%%%%%%%%%%
\begin{eqnarray*}%%%%%%%%%%%%%%%%%%%%%%%%%%%%%%%%%%%%%%%%%%%%
{\bf (c)}\ 
\partial \hat{\mathcal{M}}^1_{III}(\gamma _-,\gamma _+)
&=&
\bigcup_{q\in L_0\cap L_1} 
\hat{\mathcal{M}}^0_{II}(\gamma _-,q)
\times \hat{\mathcal{M}}^0_{II'}(q,\gamma _+)
\\
&&
\cup 
\bigcup_{(\gamma '_-,\gamma '_+)}
\hat{\mathcal{M}}^0_{III}(\gamma _-,\gamma '_-)
\times \hat{\mathcal{M}}^0_{IV}(\gamma '_-,\gamma '_+)
\times \hat{\mathcal{M}}^0_{III}(\gamma '_+,\gamma _+)
\\
&&
\cup 
\bigcup_{\gamma } 
\hat{\mathcal{M}}^0_{V}(\gamma _-,\gamma )
\times 
\hat{\mathcal{M}}^0_{III}(\gamma ,\gamma _+)
\cup 
\bigcup_{\gamma } 
\hat{\mathcal{M}}^0_{III}(\gamma _-,\gamma )
\times 
\hat{\mathcal{M}}^0_{V'}(\gamma ,\gamma _+),
\end{eqnarray*}%%%%%%%%%%%%%%%%%%%%%%%%%%%%%%%%%%%%%%%%%%%%%%
where we used the notation $\partial \hat{\mathcal{M}}_*^1(*,*)$
to denote the boundary of the compactification of 
$\hat{\mathcal{M}}_*^1(*,*)$.
\end{theo}%%%%%%%%%%%%%%%%%%%%%%%%%%%%%%%%%%%%%%%%%%%%%%%%%%%%%%%%
Note that, from the maximum principle,
families of pseudo-holomorphic curves can not grow
toward $+\infty$ of convex ends, see Section 5.
\\

Let $C$ be the free ${\bf Z}_2$-vector space over 
$L_0\cap L_1$ and 
$\bigcup _{(\gamma _-,\gamma _+)}
\hat{\mathcal{M}}^0_{IV}(\gamma _-,\gamma _+)$.
We define a linear map $\partial :C\to C$
in terms of the canonical bases:
for $p_-\in L_0\cap L_1$
\[
\partial p_-
:=\sum _{p_+\in L_0\cap L_1}
\sharp \hat{\mathcal{M}}^0_I(p_-,p_+)p_+
+
\sum _{
\stackrel{(\gamma _-, \gamma _+),}
{u\in \hat{\mathcal{M}}^0_{IV}(\gamma _-,\gamma _+)}}
\sharp \hat{\mathcal{M}}^0_{II'}(p_-,\gamma _-)
u,
\]
where $\sharp \hat{\mathcal{M}}^0_I(p_-,p_+)$ is the modulo 2
number of the elements of $\hat{\mathcal{M}}^0_I(p_-,p_+)$
and similarly
where $\sharp \hat{\mathcal{M}}^0_{II'}(p_-,\gamma _-)$
is the modulo 2 number of the elements of 
$\hat{\mathcal{M}}^0_{II'}(p_-,\gamma _-)$ 
and the second sum ranges over 
all pairs of Reeb chords
$(\gamma _-,\gamma _+)$ with respect to the concave ends.
and for $u\in \hat{\mathcal{M}}^0_{IV}(\gamma _-,\gamma _+)$
\[
\partial u
:=\sum _{p_+\in L_0\cap L_1}
\sharp \hat{\mathcal{M}}^0_{II}(\gamma _+,p_+)p_+
+
\sum _{
\stackrel{(\gamma '_-, \gamma '_+),}
{v\in \hat{\mathcal{M}}^0_{IV}(\gamma '_-,\gamma' _+)}}
\sharp \hat{\mathcal{M}}^0_{III}(\gamma _+,\gamma '_-)v,
\]
where $\sharp \hat{\mathcal{M}^0_{II}}(\gamma _+,p_+)$ 
is the modulo 2 number of the elements of 
$\hat{\mathcal{M}}^0_{II}(\gamma _+,p_+)$
and similarly
where $\sharp \hat{\mathcal{M}}^0_{III}(\gamma _+,\gamma '_-)$
is the modulo 2 number of the elements of 
$\hat{\mathcal{M}}^0_{III}(\gamma _+,\gamma '_-)$
and the second sum ranges over 
all pairs of Reeb chords
$(\gamma '_-,\gamma '_+)$ with respect to the concave ends.
\begin{ass}\label{ass3.5}%%%%%%%%%%%%%%%%%%%%%%%%%%%%%%%%%%%%%%%%%
There are no non-trivial pseudo-holomorphic strips of the type 
${\bf (V)}$ and ${\bf (V')}$.
\end{ass}%%%%%%%%%%%%%%%%%%%%%%%%%%%%%%%%%%%%%%%%%%%%%%%%%%%%%%%%%
(It seems that the existence of 
non-trivial pseudo-holomorphic strips of 
the type ${\bf (V)}$ and ${\bf (V')}$ 
is an obstruction to $\partial ^2=0$.)
Then, we can prove 
\begin{theo}\label{theo3.6}%%%%%%%%%%%%%%%%%%%%%%%%%%%%%%%%%%%%%%%
$\partial \circ \partial =0$.
\end{theo}%%%%%%%%%%%%%%%%%%%%%%%%%%%%%%%%%%%%%%%%%%%%%%%%%%%%%%%%
{\it Proof}.
For $p\in L_0\cap L_1$
\begin{eqnarray*}
\lefteqn{
\partial \partial p
=
\partial 
\sum _{q\in L_0\cap L_1}
\sharp \hat{\mathcal{M}}^0_I(p,q)q
+
\partial 
\sum _{
\stackrel{(\gamma _-,\gamma _+),}
{u\in \hat{\mathcal{M}}^0_{IV}(\gamma _-,\gamma _+)}}
\sharp \hat{\mathcal{M}}^0_{II'}(p,\gamma _-)u
}\\
&=&
\sum _{q\in L_0\cap L_1}
\sharp \hat{\mathcal{M}}^0_I(p,q)
\left(
\sum _{r\in L_0\cap L_1}
\sharp \hat{\mathcal{M}}^0_I(q,r)r
+
\sum _{
\stackrel{(\gamma '_-,\gamma '_+),}
{v\in \hat{\mathcal{M}}^0_{IV}(\gamma'_-,\gamma '_+)}}
\sharp \hat{\mathcal{M}}^0_{II'}(q,\gamma '_-)v
\right)
\\
&&
+
\sum _{
\stackrel{(\gamma _-,\gamma _+),}
{u\in \hat{\mathcal{M}}^0_{IV}(\gamma _-,\gamma _+)}}
\sharp \hat{\mathcal{M}}^0_{II'}(p,\gamma _-)
\left(
\sum _{r\in L_0\cap L_1}
\sharp \hat{\mathcal{M}}^0_{II}(\gamma _+,r)r
+
\sum _{
\stackrel{(\gamma '_-,\gamma '_+),}
{v\in \hat{\mathcal{M}}^0_{IV}(\gamma '_-,\gamma' _+)}}
\sharp \hat{\mathcal{M}}^0_{III}(\gamma _+,\gamma '_-)v
\right)
\\
&=&
\sum_{r\in L_0\cap L_1}
\left(
\sum _{q\in L_0\cap L_1}
\sharp \hat{\mathcal{M}}^0_I(p,q)
\sharp \hat{\mathcal{M}}^0_I(q,r)
+
\sum _{
\stackrel{(\gamma _-,\gamma _+),}
{u\in \hat{\mathcal{M}}^0_{IV}(\gamma _-,\gamma _+)}}
\sharp \hat{\mathcal{M}}^0_{II'}(p,\gamma _-)
\sharp \hat{\mathcal{M}}^0_{II}(\gamma _+,r)\right)r
\end{eqnarray*}
\[
+\! \! \! \! \! \! 
\sum_{
\stackrel{(\gamma '_-,\gamma '_+),}
{v\in \hat{\mathcal{M}}^0_{IV}(\gamma '_-,\gamma '_+)}}
\left(
\sum _{q\in L_0\cap L_1}
\sharp \hat{\mathcal{M}}^0_I(p,q)
\sharp \hat{\mathcal{M}}^0_{II'}(q,\gamma '_-)
+
\sum_{
\stackrel{(\gamma _-,\gamma _+),}
{u\in \hat{\mathcal{M}}^0_{IV}(\gamma _-,\gamma _+)}}
\sharp \hat{\mathcal{M}}^0_{II'}(p,\gamma _-)
\sharp \hat{\mathcal{M}}^0_{III}(\gamma _+,\gamma '_-)\right)v.
\]
The number
\[
\sum _{q\in L_0\cap L_1}
\sharp \hat{\mathcal{M}}^0_I(p,q)
\sharp \hat{\mathcal{M}}^0_I(q,r)
+
\sum _{\stackrel{(\gamma _-,\gamma _+),}
{u\in \hat{\mathcal{M}}^0_{IV}(\gamma _-,\gamma _+)}}
\sharp \hat{\mathcal{M}}^0_{II'}(p,\gamma _-)
\sharp \hat{\mathcal{M}}^0_{II}(\gamma _+,r)
\]
is nothing but the one of the boundary components of
the compactification of 
$\hat{\mathcal{M}}^1_I(p,r)$, and similarly the number
\[
\sum _{q\in L_0\cap L_1}
\sharp \hat{\mathcal{M}}^0_I(p,q)
\sharp \hat{\mathcal{M}}^0_{II'}(q,\gamma '_-)
+
\sum_{\stackrel{(\gamma _-,\gamma _+),}
{u\in \hat{\mathcal{M}}^0_{IV}(\gamma _-,\gamma _+)}}
\sharp \hat{\mathcal{M}}^0_{II'}(p,\gamma _-)
\sharp \hat{\mathcal{M}}^0_{III}(\gamma _+,\gamma '_-)
\]
is the one of the boundary components of 
the compactification of 
$\hat{\mathcal{M}}^1_{II'}(p,\gamma '_-)$,
and hence $\partial \partial p=0$.
For $u\in \hat{\mathcal{M}}^0_{IV}(\gamma _-,\gamma _+)$
\begin{eqnarray*}
\lefteqn{
\partial \partial u
=
\partial 
\sum _{p\in L_0\cap L_1}
\sharp \hat{\mathcal{M}}^0_{II}(\gamma _+,p)p
+
\partial 
\sum _{
\stackrel{(\gamma '_-,\gamma '_+),}
{v\in \hat{\mathcal{M}}^0_{IV}(\gamma '_-,\gamma' _+)}}
\sharp \hat{\mathcal{M}}^0_{III}(\gamma _+,\gamma '_-)v
}\\
&=&
\sum _{p\in L_0\cap L_1}
\sharp \hat{\mathcal{M}}^0_{II}(\gamma _+,p)
\left(
\sum _{q\in L_0\cap L_1}
\sharp \hat{\mathcal{M}}^0_I(p,q)q
+
\sum _{
\stackrel{(\gamma ''_-,\gamma ''_+),}
{w\in \hat{\mathcal{M}}^0_{IV}(\gamma ''_-,\gamma ''_+)}}
\sharp \hat{\mathcal{M}}^0_{II'}(p,\gamma ''_-)w
\right)
\\
&&
+
\sum _{
\stackrel{(\gamma '_-,\gamma '_+),}
{v\in \hat{\mathcal{M}}^0_{IV}(\gamma '_-,\gamma' _+)}}
\sharp \hat{\mathcal{M}}^0_{III}(\gamma _+,\gamma '_-)
\left(
\sum _{q\in L_0\cap L_1}
\sharp \hat{\mathcal{M}}^0_{II}(\gamma '_+,q)q
+
\sum _{
\stackrel{(\gamma ''_-,\gamma ''_+),}
{w\in \hat{\mathcal{M}}^0_{IV}(\gamma ''_-,\gamma'' _+)}}
\sharp \hat{\mathcal{M}}^0_{III}(\gamma'_+,\gamma''_-)w
\right)
\\
&=&
\sum_{q\in L_0\cap L_1}
\left(
\sum_{p\in L_0\cap L_1}
\sharp \hat{\mathcal{M}}^0_{II}(\gamma _+,p)
\sharp \hat{\mathcal{M}}^0_I(p,q)
+
\sum _{
\stackrel{(\gamma '_-,\gamma '_+),}
{v\in \hat{\mathcal{M}}^0_{IV}(\gamma '_-,\gamma' _+)}}
\sharp \hat{\mathcal{M}}^0_{III}(\gamma _+,\gamma '_-)
\sharp \hat{\mathcal{M}}^0_{II}(\gamma '_+,q)\right)q
\end{eqnarray*}
\[
+\! \! \! \! \! \! 
\sum _{
\stackrel{(\gamma ''_-,\gamma ''_+),}
{w\in \hat{\mathcal{M}}^0_{IV}(\gamma ''_-,\gamma ''_+)}}
\left(
\sum _{p\in L_0\cap L_1}
\sharp \hat{\mathcal{M}}^0_{II}(\gamma _+,p)
\sharp \hat{\mathcal{M}}^0_{II'}(p,\gamma ''_-)
+
\sum _{
\stackrel{(\gamma '_-,\gamma '_+),}
{v\in \hat{\mathcal{M}}^0_{IV}(\gamma '_-,\gamma' _+)}}
\sharp \hat{\mathcal{M}}^0_{III}(\gamma _+,\gamma '_-)
\sharp \hat{\mathcal{M}}^0_{III}(\gamma'_+,\gamma''_-)\right)w.
\]
The number
\[
\sum_{p\in L_0\cap L_1}
\sharp \hat{\mathcal{M}}^0_{II}(\gamma _+,p)
\sharp \hat{\mathcal{M}}^0_I(p,q)
+
\sum _{\stackrel{(\gamma '_-,\gamma '_+),}
{v\in \hat{\mathcal{M}}^0_{IV}(\gamma '_-,\gamma' _+)}}
\sharp \hat{\mathcal{M}}^0_{III}(\gamma _+,\gamma '_-)
\sharp \hat{\mathcal{M}}^0_{II}(\gamma '_+,q)
\]
is nothing but the one of the boundary components of
the compactification of
$\hat{\mathcal{M}}^1_{II}(\gamma _+,q)$, and similarly the number
\[
\sum _{p\in L_0\cap L_1}
\sharp \hat{\mathcal{M}}^0_{II}(\gamma _+,p)
\sharp \hat{\mathcal{M}}^0_{II'}(p,\gamma ''_-)
+
\sum _{\stackrel{(\gamma '_-,\gamma '_+),}
{v\in \hat{\mathcal{M}}^0_{IV}(\gamma '_-,\gamma' _+)}}
\sharp \hat{\mathcal{M}}^0_{III}(\gamma _+,\gamma '_-)
\sharp \hat{\mathcal{M}}^0_{III}(\gamma'_+,\gamma''_-)
\]
is the one of the boundary components of
the compactification of 
$\hat{\mathcal{M}}^1_{III}(\gamma _+,\gamma ''_-)$,
and hence $\partial \partial u=0$.
\qed
\\

We obtain the chain complex $(C,\partial )$
for $L_0$ and $L_1$ in $X$.
If complex structures on contact structures vary, then the set 
$\bigcup _{(\gamma _-,\gamma _+)}
\hat{\mathcal{M}}^0_{IV}(\gamma _-,\gamma _+)$ may change.
This implies that the generators of the type
$\bigcup _{(\gamma _-,\gamma _+)}
\hat{\mathcal{M}}^0_{IV}(\gamma _-,\gamma _+)$
may appear and disappear, hence, at the time of this writing,
the author does not know whether  
the homology is invariant under the variation 
of complex structures on contact structures.
Concerning Assumption \ref{ass3.2} and \ref{ass3.5},
he hopes that there are some relations between 
symplectic field theory \cite{egh} and our chain complexes.

Similarly we can construct Floer's chain complexes 
for periodic orbits of Hamiltonian flows on symplectic manifolds 
with concave ends, which will appear in a forthcoming paper.

%%%%%%%%%%%%%%%%%%%%%%%%%%%%%%%%%%%%%%%%%%%%%%%%%%%%%%%%%%%%%%%%%%
\section{Gluing arguments for pseudo-holomorphic strips}%%%%%%%%%%
%%%%%%%%%%%%%%%%%%%%%%%%%%%%%%%%%%%%%%%%%%%%%%%%%%%%%%%%%%%%%%%%%%

For our purpose we need the following gluing arguments.
(We will define the notation soon later.)
\begin{theo}\label{theo4.1}%%%%%%%%%%%%%%%%%%%%%%%%%%%%%%%%%%%%%%%
For the compactification of the type ${\bf (a)}$
as in Theorem \ref{theo3.4}
we need the following ${\bf (i)}$ and ${\bf (ii)}$:
\\
\\
${\bf (i)}$
Let $\hat{K}\subset \hat{\mathcal{M}}^0_I(p_-,q)$ and 
$\hat{K}'\subset \hat{\mathcal{M}}^0_I(q,p_+)$ be compact sets.
Then there exist constants $\rho _0$ and a smooth map
\[
\hat{\sharp }:\hat{K}\times \hat{K}'\times [\rho _0,\infty )
\to \hat{\mathcal{M}}^1_I(p_-,p_+).
\]
Moreover, for $\overline{u}$ and $\overline{v}$ 
in the interior of $\hat{K}$ and $\hat{K}'$,
there exist $\varepsilon >0$ and $\rho $ so that 
$\hat{\mathcal{M}}^1_I(p_-,p_+)\cap 
\hat{U}_{(\varepsilon ,\rho )}(\overline{u},\overline{v})$
is contained in the image of $\hat{\sharp }$.
\\
\\
${\bf (ii)}$
Let $\hat{K}\subset \hat{\mathcal{M}}^0_{II'}(p_-,\gamma _-)$,
$\hat{K}'\subset \hat{\mathcal{M}}^0_{IV}(\gamma _-,\gamma _+)$
and $\hat{K}''\subset \hat{\mathcal{M}}^0_{II}(\gamma _+,p_+)$
be compact subsets.
Then there exist constants $\rho _0$ and a smooth map
\[
\hat{\sharp }:\hat{K}\times \hat{K}'\times \hat{K}''\times 
[\rho _0,\infty )
\to \hat{\mathcal{M}}^1_I(p_-,p_+).
\]
Moreover, for $\overline{u}$, $\overline{w}$ 
and $\overline{v}$ in the interior of 
$\hat{K}$, $\hat{K}'$ and $\hat{K}''$,
there exist $\varepsilon >0$ and $\rho $ so that 
$\hat{\mathcal{M}}^1_I(p_-,p_+)\cap 
\hat{U}_{(\varepsilon ,\rho)}
(\overline{u},\overline{w},\overline{v})$
is contained in the image of $\hat{\sharp }$.
\\
\\
For the type ${\bf (b)}$
we need the following ${\bf (iii)}$, ${\bf (iv)}$ and ${\bf (v)}$:
\\
\\
${\bf (iii)}$
Let $\hat{K}\subset \hat{\mathcal{M}}^0_{II}(\gamma _-,q)$ and 
$\hat{K}'\subset \hat{\mathcal{M}}^0_I(q,p_+)$ be compact sets.
Then there exist constants $\rho _0$ and a smooth map
\[
\hat{\sharp }:\hat{K}\times \hat{K}'\times [\rho _0,\infty )
\to \hat{\mathcal{M}}^1_{II}(\gamma _-,p_+).
\]
Moreover, for $\overline{u}$ and 
$\overline{v}$ in the interior of $\hat{K}$ and $\hat{K}'$,
there exist $\varepsilon >0$ and $\rho $ so that 
$\hat{\mathcal{M}}^1_{II}(\gamma _-,p_+)
\cap \hat{U}_{(\varepsilon ,\rho )}(\overline{u},\overline{v})$
is contained in the image of $\hat{\sharp }$.
\\
\\
${\bf (iv)}$
Let $\hat{K}\subset 
\hat{\mathcal{M}}^0_{III}(\gamma _-,\gamma'_-)$,
$\hat{K}'\subset \hat{\mathcal{M}}^0_{IV}(\gamma '_-,\gamma '_+)$
and $\hat{K}''\subset \hat{\mathcal{M}}^0_{II}(\gamma '_+,p_+)$
be compact subsets.
Then there exist constants $\rho _0$ and a smooth map
\[
\hat{\sharp }:\hat{K}\times \hat{K}'\times \hat{K}''\times 
[\rho _0,\infty )
\to \hat{\mathcal{M}}^1_{II}(\gamma _-,p_+).
\]
Moreover, 
for $\overline{u}$, $\overline{w}$ and $\overline{v}$ 
in the interior of 
$\hat{K}$, $\hat{K}'$ and $\hat{K}''$,
there exist $\varepsilon >0$ and $\rho $ so that 
$\hat{\mathcal{M}}^1_{II}(\gamma _-,p_+)
\cap \hat{U}_{(\varepsilon ,\rho)}
(\overline{u},\overline{w},\overline{v})$
is contained in the image of $\hat{\sharp }$.
\\
\\
${\bf (v)}$
Let $\hat{K}\subset 
\hat{\mathcal{M}}^0_{V}(\gamma _-,\gamma )$ and 
$\hat{K}'\subset \hat{\mathcal{M}}^0_{II}(\gamma ,p_+)$ 
be compact sets.
Then there exist constants $\rho _0$ and a smooth map
\[
\hat{\sharp }:\hat{K}\times \hat{K}'\times [\rho _0,\infty )
\to \hat{\mathcal{M}}^1_{II}(\gamma _-,p_+).
\]
Moreover, for $\overline{u}$ and $\overline{v}$ 
in the interior of $\hat{K}$ and $\hat{K}'$,
there exists $\varepsilon >0$ and $\rho $ so that 
$\hat{\mathcal{M}}^1_{II}(\gamma _-,p_+)
\cap \hat{U}_{(\varepsilon ,\rho )}
(\overline{u},\overline{v})$
is contained in the image of $\hat{\sharp }$.
\\
\\
For the type ${\bf (b')}$
we need the following ${\bf (iii')}$, 
${\bf (iv')}$ and ${\bf (v')}$:
\\
\\
${\bf (iii')}$
Let $\hat{K}\subset \hat{\mathcal{M}}^0_I(p_-,q)$ and 
$\hat{K}'\subset \hat{\mathcal{M}}^0_{II'}(q,\gamma _+)$ 
be compact sets.
Then there exist constants $\rho _0$ and a smooth map
\[
\hat{\sharp }:\hat{K}\times \hat{K}'\times [\rho _0,\infty )
\to \hat{\mathcal{M}}^1_{II'}(p_-,\gamma _+).
\]
Moreover, for $\overline{u}$ and $\overline{v}$ 
in the interior of $\hat{K}$ and $\hat{K}'$,
there exists $\varepsilon >0$ and $\rho $ so that 
$\hat{\mathcal{M}}^1_{II'}(p_-,\gamma _+)
\cap \hat{U}_{(\varepsilon ,\rho )}(\overline{u},\overline{v})$
is contained in the image of $\hat{\sharp }$.
\\
\\
${\bf (iv')}$
Let $\hat{K}\subset \hat{\mathcal{M}}^0_{II'}(p_-,\gamma'_-)$,
$\hat{K}'\subset \hat{\mathcal{M}}^0_{IV}(\gamma '_-,\gamma '_+)$
and 
$\hat{K}''\subset \hat{\mathcal{M}}^0_{III}(\gamma '_+,\gamma _+)$
be compact subsets.
Then there exist constants $\rho _0$ and a smooth map
\[
\hat{\sharp }:\hat{K}\times \hat{K}'\times \hat{K}''\times 
[\rho _0,\infty )
\to \hat{\mathcal{M}}^1_{II'}(p_-,\gamma _+).
\]
Moreover, for $\overline{u}$, $\overline{w}$ and 
$\overline{v}$ in the interior of 
$\hat{K}$, $\hat{K}'$ and $\hat{K}''$,
there exist $\varepsilon >0$ and $\rho $ so that 
$\hat{\mathcal{M}}^1_{II'}(p_-,\gamma _+)
\cap \hat{U}_{(\varepsilon ,\rho)}
(\overline{u},\overline{w},\overline{v})$
is contained in the image of $\hat{\sharp }$.
\\
\\
${\bf (v')}$
Let $\hat{K}\subset \hat{\mathcal{M}}^0_{II'}(p_-,\gamma )$ and 
$\hat{K}'\subset \hat{\mathcal{M}}^0_{V'}(\gamma ,\gamma _+)$ 
be compact sets.
Then there exist constants $\rho _0$ and a smooth map
\[
\hat{\sharp }:\hat{K}\times \hat{K}'\times [\rho _0,\infty )
\to \hat{\mathcal{M}}^1_{II'}(p_-,\gamma _+).
\]
Moreover, for $\overline{u}$ and $\overline{v}$ 
in the interior of $\hat{K}$ and $\hat{K}'$,
there exists $\varepsilon >0$ and $\rho $ so that 
$\hat{\mathcal{M}}^1_{II'}(p_-,\gamma _+)
\cap \hat{U}_{(\varepsilon ,\rho )}(\overline{u},\overline{v})$
is contained in the image of $\hat{\sharp }$.
\\
\\
For the type ${\bf (c)}$
we need the following 
${\bf (vi)}$, ${\bf (vii)}$ and ${\bf (vii')}$:
\\
\\
${\bf (vi)}$
Let $\hat{K}\subset \hat{\mathcal{M}}^0_{II}(\gamma _-,q)$ and 
$\hat{K}'\subset \hat{\mathcal{M}}^0_{II'}(q ,\gamma _+)$ 
be compact sets.
Then there exist constants $\rho _0$ and a smooth map
\[
\sharp :\hat{K}\times \hat{K}'\times [\rho _0,\infty )
\to \hat{\mathcal{M}}^1_{III}(\gamma _-,\gamma _+).
\]
Moreover, for $\overline{u}$ and $\overline{v}$ 
in the interior of $\hat{K}$ and $\hat{K}'$,
there exists $\varepsilon >0$ and $\rho $ so that 
$\hat{\mathcal{M}}^1_{III}(\gamma _-,\gamma _+)
\cap \hat{U}_{(\varepsilon ,\rho )}
(\overline{u},\overline{v})$
is contained in the image of $\hat{\sharp }$.
\\
\\
${\bf (vii)}$
Let $\hat{K}\subset 
\hat{\mathcal{M}}^0_{III}(\gamma _-,\gamma'_-)$,
$\hat{K}'\subset \hat{\mathcal{M}}^0_{IV}(\gamma '_-,\gamma '_+)$
and $\hat{K}''\subset 
\hat{\mathcal{M}}^0_{III}(\gamma '_+,\gamma _+)$
be compact subsets.
Then there exist constants $\rho _0$ and a smooth map
\[
\hat{\sharp }:\hat{K}\times \hat{K}'\times \hat{K}''
\times [\rho _0,\infty )
\to \hat{\mathcal{M}}^1_{III}(\gamma _-,\gamma _+).
\]
Moreover, for $\overline{u}$, $\overline{w}$ and 
$\overline{v}$ in the interior of 
$\hat{K}$, $\hat{K}'$ and $\hat{K}''$,
there exist $\varepsilon >0$ and $\rho $ so that 
$\hat{\mathcal{M}}^1_{III}(\gamma _-,\gamma _+)
\cap \hat{U}_{(\varepsilon ,\rho)}
(\overline{u},\overline{w},\overline{v})$
is contained in the image of $\hat{\sharp }$.
\\
\\
${\bf (viii)}$
Let $\hat{K}\subset 
\hat{\mathcal{M}}^0_{V}(\gamma _-,\gamma )$ and 
$\hat{K}'\subset \hat{\mathcal{M}}^0_{III}(\gamma ,\gamma _+)$ 
be compact sets.
Then there exist constants $\rho _0$ and a smooth map
\[
\sharp :\hat{K}\times \hat{K}'\times [\rho _0,\infty )
\to \hat{\mathcal{M}}^1_{III}(\gamma _-,\gamma _+).
\]
Moreover, for $\overline{u}$ and $\overline{v}$ 
in the interior of $\hat{K}$ and $\hat{K}'$,
there exists $\varepsilon >0$ and $\rho $ so that 
$\hat{\mathcal{M}}^1_{III}(\gamma _-,\gamma _+)
\cap \hat{U}_{(\varepsilon ,\rho )}(\overline{u},\overline{v})$
is contained in the image of $\hat{\sharp }$.
\\
\\
${\bf (viii')}$
Let $\hat{K}\subset 
\hat{\mathcal{M}}^0_{III}(\gamma _-,\gamma )$ and 
$\hat{K}'\subset \hat{\mathcal{M}}^0_{V'}(\gamma ,\gamma _+)$ 
be compact sets.
Then there exist constants $\rho _0$ and a smooth map
\[
\hat{\sharp }:\hat{K}\times \hat{K}'\times [\rho _0,\infty )
\to \hat{\mathcal{M}}^1_{III}(\gamma _-,\gamma _+).
\]
Moreover, for $\overline{u}$ and $\overline{v}$ 
in the interior of $\hat{K}$ and $\hat{K}'$,
there exists $\varepsilon >0$ and $\rho $ so that 
$\hat{\mathcal{M}}^1_{III}(\gamma _-,\gamma _+)
\cap \hat{U}_{(\varepsilon ,\rho )}
(\overline{u},\overline{v})$
is contained in the image of $\hat{\sharp }$.
\end{theo}%%%%%%%%%%%%%%%%%%%%%%%%%%%%%%%%%%%%%%%%%%%%%%%%%%%%%%%%

In the following we denote by
$\beta :{\bf R}\to [0,1]$ a smooth function
\[
\beta (\tau ):=\left\{
\begin{array}{ll}
0,\ \tau \leq 0,
\\
1,\ 1\leq \tau,
\end{array}
\right.
\]
and by $\|f\|_{[r_1,r_2];L^p_1}$
the restriction of the $L^p_1$ norm on 
$[r_1,r_2]\times [0,1]$, i.e.,
$\|f\|_{[r_1,r_2];L^p_1}:=
\{\int_{[r_1,r_2]\times [0,1]}|f|^p+|Df|^pd\tau dt\}^{1/p}$,
and also $\|f\|_{[r_1,r_2];L^p}$.
\\

First we recall the proof of the gluing argument ${\bf (i)}$.
(This part is so standard, 
the reader may skip to the next content.)
For compact sets $K\subset \mathcal{M}^1_I(p_-,q)$ and 
$K'\subset \mathcal{M}^1_I(q,p_+)$,
we define 
$w_{\chi}(\tau ,t)\in L^p_1({\bf R}\times [0,1];X,L_0,L_1)$
for $\chi :=(\overline{u},\overline{v},\rho )
\in K\times K'\times [\rho _0,\infty )$ by 
\[
w_{\chi }(\tau ,t):=
\left\{
\begin{array}{lc}
\overline{u}(\tau +\rho ,t),& \tau \leq -1,
\\
\exp _q
\left(
\beta (-\tau )\eta  (\tau +\rho ,t)+
\beta (\tau )\zeta (\tau -\rho,t)
\right),
& -1\leq \tau \leq 1,
\\
\overline{v}(\tau -\rho ,t),
& 1\leq \tau,
\end{array}
\right.
\]
where $\overline{u}=\exp_q\eta$ for $\rho -1\leq \tau $
and
$\overline{v}=\exp_q\zeta$ for $\tau \leq -\rho +1$.
From $\overline{\partial }_J\overline{u}
=\overline{\partial }_J\overline{v}=0$ we conclude that 
\begin{eqnarray}
\|\overline{\partial }_{J_t}w_{\chi }\|_{L^p}\leq 
C\left(
e^{-d_+\rho }\|e^{d_+\tau }\eta  \| _{[\rho -1,\rho];L^p_1}+
e^{-d_-\rho}\|e^{-d_-\tau }\zeta \| _{[-\rho ,1-\rho];L^p_1}
\right),
\label{5}
\end{eqnarray}
where $C$ is a constant depending only on $K$, $K'$ and $\rho _0$.
We denote the Taylor expansion of 
$f_{w_{\chi }}$%(\eta ):=\Phi _{w_{\chi }}(\eta )^{-1}
by
\[
f_{w_{\chi }}(\xi ):=f_{w_{\chi }}(0)
+Df_{w_{\chi }}(0)\xi +N_{w_{\chi }}(\xi ).
\]
\begin{lem}\label{lem4.2}%%%%%%%%%%%%%%%%%%%%%%%%%%%%%%%%%%%%%%%%%
For $\|\eta \|_{L^p_1}\leq c$ and $\|\zeta \|_{L^p_1}\leq c$,
the nonlinear part $N_{w_{\chi }}$ satisfies the estimate
\begin{eqnarray}
\|N_{w_{\chi }}(\eta )-N_{w_{\chi }}(\zeta )\|_{L^p}
\leq C
(\|\eta \|_{L^p_1}+\|\zeta \|_{L^p_1})\|\eta -\zeta \|_{L^p_1},
\label{6}
\end{eqnarray}
where $C$ is a constant 
depending only on $\|\nabla w_{\chi }\|_{L^p}$ and $c$.
\end{lem}%%%%%%%%%%%%%%%%%%%%%%%%%%%%%%%%%%%%%%%%%%%%%%%%%%%%%%%%%
{\it Proof}. Basically it is done by the Taylor expansion.
(In the following we shall use $f$ to denote $f_{w_{\chi }}$.)
\begin{eqnarray*}
& &N_{w_{\chi }}(\eta )-N_{w_{\chi }}(\zeta)
\\
&=&
\int_0^1(1-t)
\{(d^2f)_{t\eta }(\eta,\eta)
-
(d^2f)_{t\zeta }(\zeta,\zeta)\}dt
\\
&=&
\int_0^1(1-t)
\{(d^2f)_{t\eta }(\eta,\eta-\zeta)
+(d^2f)_{t\eta }(\eta,\zeta)-(d^2f)_{t\zeta }(\eta,\zeta)
+
(d^2f)_{t\zeta}(\eta-\zeta,\zeta)\}dt
\\
&=&
\int_0^1(1-t)
\left\{(d^2f)_{t\eta }(\eta,\eta-\zeta)
+
\int_0^1(d^3f)_{(1-s)t\eta +st\zeta }(t\eta -t\zeta ,\eta,\zeta)ds
+
(d^2f)_{t\zeta}(\eta-\zeta,\zeta)\right\}dt.
\end{eqnarray*}
Then we can conclude
\[
\|N_{w_{\chi }}(\eta )-N_{w_{\chi }}(\zeta )\|_{L^p}
\leq 
C(\|\eta \|_{L^p_1}+\|\eta \|_{L^p_1}\|\zeta \|_{L^p_1}
+\|\zeta \|_{L^p_1})\|\eta -\zeta \|_{L^p_1},
\]
where the constant $C$ depends only on 
$\|\nabla w_{\chi }\|_{L^p}$ and $c$. 
The reason of the independence of $C$ from 
$\|w_{\chi }\|_{L^{\infty }}$ is 
the compactness of the contact manifolds of the cylinders
and boundedness of metrics and connections and so on.
\qed  
\\

Let $\chi _i:=(\overline{u}_i,\overline{v}_i,\rho _i)$, 
$i=1,2,\ldots $, be a sequence of 
$K\times K'\times [\rho _0,\infty )$ such that 
$\rho _i\to \infty $.
We assume that $w_{\chi _i}((-\rho _i,\rho _i),[0,1])$
is contained in the Gaussian coordinate of $q$.
For $\xi _i\in 
L^p_1(w_{\chi _i}^*TX,w_{\chi _i}^*TL_0,w_{\chi _i}^*TL_1)$,
define $\xi _{0i}\in 
L^p_1([-\rho _i+\sqrt{\rho _i}-1,\rho _i-\sqrt{\rho _i}+1]
\times [0,1];T_qX,T_qL_0,T_qL_1)$ 
such that 
\[
D\exp _q(\xi _{0i}(\tau ,t))=\xi _i(\tau ,t).
\]
Note that, if we put $\beta _i(\tau ):=
\beta (-\tau +\rho _i-\sqrt{\rho _i}+1)
\beta (\tau +\rho _i-\sqrt{\rho _i}+1)$,
then $\beta _i\xi _{0i}$ is an element of 
$L^p_1(T_qX,T_qL_0,T_qL_1)$.
Similarly define differential operators
$E_{0i}$ on 
$[-\rho _i+\sqrt{\rho _i}-1,\rho _i-\sqrt{\rho _i}+1]\times [0,1]$
by
\[
D\exp _q(E_{0i}\xi _{0i})=E_i\xi _i,
\]
where $E_i$ denotes the differential $Df_{w_{\chi _i}}(0)$. 
Note that the sequence of $\{E_{0i}\}_{i=1,2,\ldots }$ 
converges to the standard Cauchy-Riemann operator 
$\overline{\partial }_0$ on $[R_1,R_2]\times [0,1]$. 
(The convergence means that, if we denote 
$E_{0i}=a_i\frac{\partial }{\partial \tau }
+b_i\frac{\partial }{\partial t}+c_i$,
then $a_i\to 1$, $b_i\to J_0$ and $c_i\to 0$
in the $C^{\infty }$ topology.)
\begin{pro}\label{pro4.3}%%%%%%%%%%%%%%%%%%%%%%%%%%%%%%%%%%%%%%%%%
If
$\|\xi _i\|_{L^p_1}=1$ and $\|E_i\xi _i\|_{L^p}\to 0$, 
then there exists a subsequence $\{\xi _{i_l}\}$ such that 
\[
\|\xi _{i_l}\|_{[R_1,R_2];L^p_1}\to 0.
\]
\end{pro}%%%%%%%%%%%%%%%%%%%%%%%%%%%%%%%%%%%%%%%%%%%%%%%%%%%%%%%%%
{\it Proof}.
{}From the assumption $\|\xi _i\|_{L^p_1}=1$,
there is a constant $C$ such that 
$\|\beta _i\xi _{0i}\|_{L^p_1}\leq C$.
By the Rellich's theorem, there exists 
$\xi _{(R_1,R_2)}\in L^p([R_1,R_2]\times [0,1];
T_qX,T_qL_0,T_qL_1)$ and a subsequence 
$\{\beta _{i_l}\xi _{0i_l}\}$
such that 
\[
\|\xi _{(R_1,R_2)}-\beta _{i_l}\xi _{0i_l}\|_{[R_1,R_2];L^p}
\to 0.
\]
For simplicity, we use $\xi _{0i}$ to denote $\beta _i\xi _{0i}$
and assume $\{\beta _i\xi _{0i}\}$ satisfies 
$\|\beta _i\xi _{0i}- \xi _{(R_1,R_2)}\|_{[R_1,R_2];L^p}\to 0$ 
and 
$\xi _{(R'_1,R'_2)}|_{[R_1,R_2]\times [0,1]}=\xi _{(R_1,R_2)}$
for $[R_1,R_2]\subset [R'_1,R'_2]$.
By the G\"{a}rding's inequality 
\[
\|\xi _{0i}-\xi _{0j}\|_{[R_1,R_2];L^p_1}
\leq C
(\|E_{0i}(\xi _{0i}-\xi _{0j})\|_{[R_1+\delta ,R_2+\delta ];L^p}
+\|\xi _{0i}-\xi _{0j}\|_{[R_1+\delta ,R_2+\delta ];L^p}),
\]
where $C$ is a constant depending on $R_1,R_2$ and $\delta $.
We already know  
$\|\xi _{0i}-\xi _{0j}\|_{[R_1+\delta ,R_2+\delta ];L^p}\to 0$,
and from
\begin{eqnarray*}
&&
\| E_{0i}(\xi _{0i}-\xi _{0j})\|_{[R_1+\delta ,R_2+\delta ];L^p}
\\
& \leq &
\| E_{0i}\xi _{0i}\|_{[R_1+\delta ,R_2+\delta ];L^p}
+\| (E_{0i}-E_{0j})\xi _{0j}\|_{[R_1+\delta ,R_2+\delta ];L^p}
+\| E_{0j}\xi _{0j}\|_{[R_1+\delta ,R_2+\delta ];L^p}
\\
& \leq &
\| E_{0i}\xi _{0i}\|_{[R_1+\delta ,R_2+\delta ];L^p}
+\|E_{0i}-E_{0j}\|_{[R_1+\delta ,R_2+\delta ];L^p}
\|\xi _{0j}\|_{[R_1+\delta ,R_2+\delta ];L^p_1}
+\| E_{0j}\xi _{0j}\|_{[R_1+\delta ,R_2+\delta ];L^p}
\\ 
&\to & 0,
\end{eqnarray*}
we can conclude 
$\|E_{0i}(\xi _{0i}-\xi _{0j})\|_{[R_1+\delta ,R_2+\delta ];L^p}
\to 0$.
Then $\{\xi _{0i}\}$ has a subsequence which converges to 
$\xi _{(R_1,R_2)}$ in the norm $\| \cdot \|_{[R_1,R_2];L^p_1}$.
Moreover, from $\|\beta _i\xi _{0i}\|_{L^p_1}\leq C$,
we can conclude $\|\xi _{\infty }\|_{L^p_1}\leq C$,
where $\xi _{\infty }|_{[R_1,R_2]\times [0,1]}
:=\xi _{(R_1,R_2)}$. 
On the other hand, from 
\begin{eqnarray*}
&&\|\overline{\partial }_0\xi _{\infty }\|_{[R_1,R_2];L^p}
\\
& \leq &
\|E_{0i}\xi _{0i}\|_{[R_1,R_2];L^p}+
\|(E_{0i}-\overline{\partial }_0)\xi _{0i}\|_{[R_1,R_2];L^p}
+\|\overline{\partial }_0(\xi _{0i}-\xi _{\infty })
\|_{[R_1,R_2];L^p}
\\
&\leq &
\|E_{0i}\xi _{0i}\|_{[R_1,R_2];L^p}+
\|E_{0i}-\overline{\partial }_0\|_{[R_1,R_2];L^p}
\|\xi _{0i}\|_{[R_1,R_2];L^p_1}
+C\|\xi _{0i}-\xi _{\infty }\|_{[R_1,R_2];L^p_1}
\\
&\to & 0,
\end{eqnarray*}
$\overline{\partial }\xi _{\infty }=0$
and hence $\xi _{\infty }=0$.
Finally, from $\|\xi _i\|_{[R_1,R_2];L^p_1}\leq 
C\|\xi _{0i}\|_{[R_1,R_2];L^p_1}$, 
there exists a subsequence $\{\xi _{i_l}\}$
such that 
$\|\xi _{i_l}\|_{[R_1,R_2];L^p_1}\to 0$.
\qed
\\

For $\eta \in \mbox{Ker}E_{\overline{u}}$ and  
$\zeta \in \mbox{Ker}E_{\overline{v}}$,
we define 
$\eta \sharp _{\rho }\zeta \in 
L^p_1(w_{\chi }^*TX,w_{\chi }^*TL_0,w_{\chi }^*TL_1)$ by
\[
(\eta \sharp _{\rho }\zeta )(\tau ,t)
:=
\left\{
\begin{array}{ll}
\beta (-\tau -2)\eta (\tau +\rho _i), &\tau \leq 2,
\\
0, &-2\leq \tau \leq 2,
\\
\beta (\tau -2)\zeta (\tau -\rho _i), &2\leq \tau.
\end{array}
\right.
\]
Let 
$W_{w_{\chi }}^{\perp}$ 
be the $L^2$-orthogonal compliment of 
$W_{w_{\chi }}:=\{\eta \sharp _{\rho }\zeta |
\eta \in \mbox{Ker}E_{\overline{u}},
\zeta \in \mbox{Ker}E_{\overline{v}}\}$ in 
$L^p_1(w_{\chi }^*TX,w_{\chi }^*TL_0,w_{\chi }^*TL_1)$.
(Note that the dimension of $W_{w_{\chi }}$ is equal to 
$\dim \mbox{Ker}E_{\overline{u}} +
\dim \mbox{Ker}E_{\overline{v}}$.)
\begin{pro}\label{pro4.4}%%%%%%%%%%%%%%%%%%%%%%%%%%%%%%%%%%%%%%%%%
There exist constants $\rho _0$ and $C$ such that 
for $\chi \in K\times K'\times [\rho _0,\infty )$
and $\xi \in W_{w_{\chi }}^{\perp}$
\[
\|\xi \|_{L^p_1}\leq C\|E_{w_{\chi }}\xi \|_{L^p}.
\]
\end{pro}%%%%%%%%%%%%%%%%%%%%%%%%%%%%%%%%%%%%%%%%%%%%%%%%%%%%%%%%%
{\it Proof}.
Let $\chi _i:=(\overline{u}_i,\overline{v}_i,\rho _i),
i=1,2,\ldots $, 
be a sequence such that $\rho _i\to \infty$ 
and
there exist $\xi _i\in W_{w_{\chi _i}}^{\perp }$ 
satisfying
$\|\xi _i\|_{L^p_1}=1$ and $\|E_{w_{\chi _i}}\xi _i\|_{L^p}\to 0$.
Then we can conclude from Proposition \ref{pro4.3}
there exists a subsequence $\{\xi _{i_l}\}$ such that 
$\|\xi _{i_l}\|_{[-3,3];L^p_1}\to 0$.
For simplicity we denote this subsequence by $\{\xi _i\}$
and define
\begin{eqnarray*}
\eta _i(\tau ,t)&
:=&\beta (-\tau +\rho _i-1)\xi _i(\tau -\rho _i,t)
\in L^p_1(\overline{u}_i^*TX,
\overline{u}_i^*TL_0,\overline{u}_i^*TL_1),
\\
\zeta _i(\tau ,t)&
:=&\beta (\tau +\rho _i-1)\xi _i(\tau +\rho _i,t)
\in L^p_1(\overline{v}_i^*TX,
\overline{v}_i^*TL_0,\overline{v}_i^*TL_1).
\end{eqnarray*}
Split $\eta _i$ into $k_i+n_i$, where 
$k_i\in \mbox{Ker}E_{\overline{u}_i}$ and 
$n_i\in (\mbox{Ker}E_{\overline{u}_i})^{\perp}$, then
\begin{eqnarray*}
\|\eta _i\|_{L^p_1}
&\leq &
\|k_i\|_{L^p_1}+\|n_i\|_{L^p_1}
\\
&\leq &
\|k_i\|_{L^p_1}+C_{\overline{u}_i}\|E_{\overline{u}_i}n_i\|_{L^p}
\\
&\leq &
\|k_i\|_{L^p_1}+C_{\overline{u}_i}
\|E_{\overline{u}_i}\eta _i\|_{L^p}
\\
&=&
\|k_i\|_{L^p_1}+C_{\overline{u}_i}\|E_{w_{\chi _i}}
\left(\beta (-\tau +\rho _i-1)
\xi _i(\tau -\rho _i,t)\right)\|_{L^p}
\\
&\leq &
\|k_i\|_{L^p_1}+
C_{\overline{u}_i}'\|\xi _i\|_{[-2,-1];L^p_1}
+C_{\overline{u}_i}''\|E_{w_{\chi _i}}\xi _i\|_{L^p}
\end{eqnarray*}
(Note that, from the compactness of $K$,
$C_{\overline{u}_i}'$ and $C_{\overline{u}_i}''$ are bounded.)
We already know
$\|\xi _i\|_{[-2,-1];L^p_1}\to 0$ and $\|E_{w_{\chi _i}}
\xi _i\|_{L^p}\to 0$.
Let $\{e^i_1,\ldots ,e^i_r\}$ be orthogonal bases of 
$\mbox{Ker}E_{\overline{u}_i}$. 
\begin{eqnarray*}
|\langle k_i,e^i_j\rangle |
&=&
|\langle \eta _i,e^i_j\rangle |
\\
&=&
|\langle \eta _i,(1-\beta (-\tau -\rho _i+2))e^i_j\rangle 
|_{[\rho _i-2,\rho _i-1]\times [0,1]}|
\\
&\leq &
\|\eta _i\|_{[\rho _i-2,\rho _i-1];L^p}
\|e ^i_j\|_{[\rho _i-2,\rho _-1];L^{p/(p-1)}}
\\
&\leq &
C\|\xi _i\|_{[-2,-1];L^p}
\|e ^i_j\|_{[\rho _i-2,\rho _i-1];L^{p/(p-1)}}.
\end{eqnarray*}
Hence, from $\|k_i\|_{L^p_1}\leq \sum _j
|\langle k_i,e^i_j\rangle |
\|e^i_j\|_{L^p_1}$, 
we conclude $\|k_i\|_{L^p_1}\to 0$.
(Note that the compactness of $K$ induces 
the boundedness of norms of $e^i_j$.)
Then we obtain $\|\eta _i\|_{L^p_1}\to 0$.
Similarly we can prove also $\|\zeta _i\|_{L^p_1}\to 0$.
Put together them with
\[
\|\xi _i\|_{L^p_1}
\leq 
\|\eta _i\|_{L^p_1}+\|\xi _i\|_{[-3,3];L^p_1}
+\|\zeta _i\|_{L^p_1},
\]
we obtain $\|\xi _i\|_{L^p_1}\to 0$ which is a contradiction to 
the assumption $\|\xi _i\|_{L^p_1}=1$.
We finish proving the proposition.
\qed
\begin{pro}\label{pro4.5}%%%%%%%%%%%%%%%%%%%%%%%%%%%%%%%%%%%%%%%%%
There exist constants $\rho _0$ and $C$ 
such that 
for $\chi \in K\times K'\times [\rho _0,\infty )$
there exists a map
$G_{w_{\chi }}:L^p(w_{\chi }^*TX)\to W_{w_{\chi }}^{\perp}
\subset L^p_1(w_{\chi }^*TX, w_{\chi }^*TL_0,w_{\chi }^*TL_1)$
such that
\begin{eqnarray}
E_{w_{\chi }}G_{w_{\chi }}&=&\mbox{{\rm id}},
\nonumber
\\
\|G_{w_{\chi }}\xi \|_{L^p_1}&\leq & C\|\xi \|_{L^p}.
\label{7}
\end{eqnarray}
\end{pro}%%%%%%%%%%%%%%%%%%%%%%%%%%%%%%%%%%%%%%%%%%%%%%%%%%%%%%%%%
{\it Proof}.
From Proposition \ref{pro4.4}, if 
$\xi \in \mbox{Ker}E_{w_{\chi }}\cap W_{w_{\chi }}^{\perp}$, 
then $\xi =0$ and 
\[
\dim \mbox{Ker}E_{w_{\chi }}
\leq 
\dim \mbox{Ker}E_{\overline{u}}
+\dim \mbox{Ker}E_{\overline{v}}.
\]
On the other hand, from Assumption \ref{ass3.1},
\begin{eqnarray*}
\dim \mbox{Ker}E_{w_{\chi }}
& \geq & \mbox{Ind}E_{w_{\chi }}
\\
&=&
\mbox{Ind}E_{\overline{u}}+\mbox{Ind}E_{\overline{v}}
\\
&=&
\dim \mbox{Ker}E_{\overline{u}}
+\dim \mbox{Ker}E_{\overline{v}}.
\end{eqnarray*}
Hence we obtain 
\[
\mbox{Ker}E_{w_{\chi }}\oplus W_{w_{\chi }}^{\perp}=
L^p_1(w_{\chi }^*TX,w_{\chi }^*TL_0,w_{\chi }^*TL_1),
\]
and $E_{w_{\chi }}$ is surjective.
From the estimate of Proposition \ref{pro4.4} 
we can obtain $G_{w_{\chi }}$ as in the proposition.
\qed
\\

So far we obtain the following:
Let $K\subset \mathcal{M}^1_I(p_-,q)$ and
$K'\subset \mathcal{M}^1_I(q,p_+)$ be compact sets.
Then there are constants $\rho _0$ such that,
for $\chi :=(\overline{u},\overline{v},\rho)\in
K\times K'\times [\rho _0,\infty )$, a map 
$
f_{w_{\chi }}:
L^p_1(w_{\chi }^*TX,w_{\chi }^*TL_0,w_{\chi }^*TL_1)\to 
L^p(w_{\chi }^*TX)$
satisfies (\ref{5}) and (\ref{6}) and 
$Df_{w_{\chi }}(0)$ possesses a right inverse 
$G_{w_{\chi }}:L^p(w_{\chi }^*TX)\to 
L^p_1(w_{\chi }^*TX,w_{\chi }^*TL_0,w_{\chi }^*TL_1)$
satisfying (\ref{7}).
Then, from the Newton's method in Appendix A,
we can conclude that there are constants $\rho _0$ and 
$C$ and  a smooth map 
\[
\sharp :K\times K'\times [\rho _0,\infty )\to 
\mathcal{M}_I^2(p_-,p_+),\
\chi \mapsto \exp _{w_{\chi }}(\xi _{\chi })
\]
with 
$\|\xi _{\chi }\|_{L^p_1}\leq 
C\|\overline{\partial }_{J_t}w_{\chi}\|_{L^p}$.
Divide them by the ${\bf R}$ actions, and we obtain a gluing map
$\hat{\sharp }:\hat{K}\times \hat{K}'\times [\rho _0,\infty )\to 
\hat{\mathcal{M}}_I^1(p_-,p_+)$.
\\

The next step is to show the surjectivity of
$\sharp :K\times K'\times [\rho _0,\infty )\to 
\mathcal{M}_I^2(p_-,p_+)\cap U_{(\varepsilon ,\rho )}
(\overline{u},\overline{v})$.
Let $w$ be a map which satisfies 
the Lagrangian boundary conditions and 
${\bf (I)}$ and the decay conditions ${\bf (1)}$ and ${\bf (1')}$,
and $w(\tau ,t)=\exp _q\eta (\tau ,t)$ for $-1\leq \tau \leq 1$,
where $q\in L_0\cap L_1$.
Then we define 
\[
x_{\rho }(\tau ,t):=\left\{
\begin{array}{ll}
w(\tau -\rho ,t), & \tau \leq \rho -1,
\\
\exp _q\beta (-\tau +\rho )\eta (\tau -\rho ,t), 
& \rho -1\leq \tau ,
\end{array}
\right.
\]
\[
y_{\rho }(\tau ,t):=\left\{
\begin{array}{ll}
\exp _q\beta (\tau +\rho )\eta (\tau +\rho ,t), & 
\tau \leq -\rho +1,
\\
w(\tau +\rho ,t), & -\rho +1\leq \tau .
\end{array}
\right.
\]
Moreover we define 
$
U_{(\varepsilon ,\rho _0)}(\overline{u},\overline{v})
$,
for $(\overline{u},\overline{v})\in K\times K'$,
to be the set of $w$ such that,
for $\rho >\rho _0$,
$\|\overline{u}-x_{\rho }\|_{L^p_1}<\varepsilon $
and $\|\overline{v}-y_{\rho }\|_{L^p_1}<\varepsilon $.
(For simplicity, we shall use a letter $x$ to denote $x_{\rho }$, 
and also $y:=y_{\rho }$.)
If $w\in \mathcal{M}_I^2(p_-,p_+)\cap 
U_{(\varepsilon ,\rho _0)}(\overline{u},\overline{v})$, then 
for a smooth map
$f_x$
there are constants $C$ and $C'$ such that
\begin{eqnarray}
\|f_x(0)\|_{L^p}
&\leq &C\|\eta \|_{[-1,0];L^p_1},
\label{8}
\\
\|N_x(\xi )-N_x(\xi ')\|_{L^p}
&\leq &C'(\|\xi \|_{L^p_1}+\|\xi '\|_{L^p_1})
\|\xi -\xi '\|_{L^p_1},
\label{9}
\end{eqnarray}
where $\|\xi \|_{L^p_1}\leq c$ and $\|\xi ' \|_{L^p_1}\leq c$ and
$C'$ depends on $c$,
and also for $f_y$.
(The proofs of the above estimates are similar to 
those of (\ref{5}) and (\ref{6}).)
For $\xi \in \mbox{Ker}E_w$,
we define 
$\hat{\xi }:=(\xi _x,\xi _y)\in 
L^p_1(x^*TX,x^*TL_0,x^*TL_1)\oplus L^p_1(y^*TX,y^*TL_0,y^*TL_1)$
by 
\begin{eqnarray*}
\xi _x&:=&
\beta (-\tau +\rho +1)\xi (\tau -\rho ,t),
\\
\xi _y&:=&
\beta (\tau +\rho +1)\xi (\tau +\rho ,t).
\end{eqnarray*}
Let $W_{(w,\rho )}^{\perp }$ 
the $L^2$-orthogonal compliment of $W_{(w,\rho )}
:=\{\hat{\xi }|
\xi \in \mbox{Ker}E_w
\}$ in
$L^p_1(x^*TX,x^*TL_0,x^*TL_1)\oplus L^p_1(y^*TX,y^*TL_0,y^*TL_1)$.
\begin{pro}\label{pro4.6}%%%%%%%%%%%%%%%%%%%%%%%%%%%%%%%%%%%%%%%%%
There are constants $\varepsilon _0>0$ and $\rho _0$ and $C$ 
such that for 
$w\in \mathcal{M}_I^2(p_-,p_+)\cap 
U_{(\varepsilon _0,\rho _0)}(\overline{u},\overline{v})$ and 
$\xi :=(\xi _x,\xi _y) \in W_{(w,\rho )}^{\perp}$
\[
\|\xi \|_{L^p_1}\leq C\|(E_x\xi _x,E_y\xi _y)\|_{L^p}.
\]
\end{pro}%%%%%%%%%%%%%%%%%%%%%%%%%%%%%%%%%%%%%%%%%%%%%%%%%%%%%%%%%
{\it Proof}.
Let $\{\varepsilon _i \}_{i=1,2,\ldots }$ 
and $\{\rho _i\}_{i=1,2,\ldots }$ and 
$w_i\in \mathcal{M}_I^2(p_-,p_+)\cap 
U_{(\varepsilon _i,\rho _i)}(\overline{u},\overline{v})$
be sequences such that 
$\varepsilon _i\to 0$ and $\rho _i\to \infty $
and there exist $\xi _i=(\xi _{x_i},\xi _{y_i})
\in W_{(w_i,\rho _i)}^{\perp }$ satisfying
$\|\xi _i\|_{L^p_1}=1$ and 
$\|(E_{x_i}\xi _{x_i}, E_{y_i}\xi _{y_i})\|_{L^p}\to 0$. 
Then, in a similar way to the proof of Proposition \ref{pro4.4}, 
we can prove that there exists a subsequence 
$\{\xi _{i_l}\}$ such that 
$\|\xi _{i_l}\|_{L^p_1}\to 0$, which contradicts the 
assumption $\|\xi _i\|_{L^p_1}=1$.
\qed
\begin{pro}\label{pro4.7}%%%%%%%%%%%%%%%%%%%%%%%%%%%%%%%%%%%%%%%%%
There exist constants 
$\varepsilon _0>0$ and $\rho _0$ and $C$ 
such that 
for 
$w\in \mathcal{M}_I^2(p_-,p_+)\cap 
U_{(\varepsilon _0,\rho _0)}(\overline{u},\overline{v})$
there exists a map 
$G_{(w,\rho )}
:L^p(x^*TX)\oplus L^p(y^*TX)\to W_{(w,\rho )}^{\perp }$ such that
\begin{eqnarray}
(E_x\oplus E_y)G_{(w,\rho )} &=& \mbox{{\rm id}},
\nonumber
\\
\|G_{(w,\rho )}\xi \|_{L^p_1} &\leq & C\|\xi \|_{L^p}.
\label{10}
\end{eqnarray}
\end{pro}%%%%%%%%%%%%%%%%%%%%%%%%%%%%%%%%%%%%%%%%%%%%%%%%%%%%%%%%%
{\it Proof}.
From Proposition \ref{pro4.6}, if 
$\xi \in \mbox{Ker}(E_x\oplus E_y)\cap W_{(w,\rho )}^{\perp }$,
then $\xi =0$ and
\[
\dim \mbox{Ker}(E_x\oplus E_y)\leq \dim W_{(w,\rho )}.
\]
On the other hand, from Assumption \ref{ass3.1},
\begin{eqnarray*}
\dim \mbox{Ker}(E_x\oplus E_y)&\geq &\mbox{Ind}E_x +\mbox{Ind}E_y
\\
&=& \dim W_{(w,\rho )}.
\end{eqnarray*}
Hence we obtain
\[
\mbox{Ker}(E_x\oplus E_y)\oplus W_{(w,\rho )}^{\perp}
=L^p_1(x^*TX,x^*TL_0,x^*TL_1)\oplus L^p_1(y^*TX,y^*TL_0,y^*TL_1),
\]
and $E_x\oplus E_y$ is surjective.
From the estimate of Proposition \ref{pro4.6} we can obtain
$G_{(w,\rho )}$ as in the proposition.
\qed
\\

So far we obtain the following:
There are constants $\varepsilon _0>0$ and 
$\rho _0$ such that,
for $w\in \mathcal{M}_I^2(p_-,p_+)\cap 
U_{(\varepsilon _0,\rho _0)}(\overline{u},\overline{v})$,
a map 
$
f_{(x,y)}:=(f_x,f_y):
L^p_1(x^*TX,x^*TL_0,x^*TL_1)\oplus L^p_1(y^*TX,y^*TL_0,y^*TL_1)
\to 
L^p(x^*TX)\oplus L^p(y^*TX)$ satisfies
(\ref{8}) and (\ref{9}),
and $Df_{(x,y)}(0)$ possesses a right inverse 
$G_{(w,\rho )}:L^p(x^*TX)\oplus L^p(y^*TX)
\to L^p_1(x^*TX,x^*TL_0,x^*TL_1)\oplus 
L^p_1(y^*TX,y^*TL_0,y^*TL_1)$ satisfying (\ref{10}).
Then, from the Newton's method in Appendix A,
we can conclude that there are constants 
$\varepsilon >0$ and $\rho $ and $C$ and  a smooth map 
\[
\sharp ':\mathcal{M}_I^2(p_-,p_+)\cap 
U_{(\varepsilon ,\rho )}(\overline{u},\overline{v})\to 
K\times K',\
w\mapsto (\exp _x(\xi  _{w;x}),
\exp _y(\xi _{w;y}))
\]
with 
$\|\xi _{w;x}\|_{L^p_1}\leq 
C\|\overline{\partial }_{J_t}x\|_{L^p}$ and
$\|\xi _{w;y}\|_{L^p_1}\leq 
C\|\overline{\partial}_{J_t}y\|_{L^p}$.
Divide them by the ${\bf R}$ actions, and we obtain 
a map $\hat{\sharp }':\hat{\mathcal{M}}_I^1(p_-,p_+)\cap 
\hat{U}_{(\varepsilon ,\rho )}(\overline{u},\overline{v})\to 
\hat{K}\times \hat{K}'$.

{}From the construction of $\sharp $ and $\sharp '$,
if $\rho _0$ is large and $\varepsilon $ is small enough,
then $\sharp \circ \sharp '$ and $\sharp '\circ \sharp $
are diffeomorphisms.
We finish proving the gluing argument ${\bf (i)}$.
\\

Next we will prove the gluing argument ${\bf (ii)}$.
Take a lift of 
$\hat{\mathcal{M}}^0_{II'}(p_-,\gamma _-)\times 
\hat{\mathcal{M}}^0_{IV}(\gamma _-,\gamma _+)\times
\hat{\mathcal{M}}^0_{II}(\gamma _+,p_+)$ 
in 
$\mathcal{M}^1_{II'}(p_-,\gamma _-)\times 
\mathcal{M}^2_{IV}(\gamma _-,\gamma _+)\times
\mathcal{M}^1_{II}(\gamma _+,p_+)
$ 
and consider the orbit of the lift 
by the following ${\bf R}^2$-action:
for $(a,0)\in {\bf R}^2$
\[
(a,0)\cdot (\overline{u},\overline{w},\overline{v})
:=\left(
\overline{u},a*\overline{w},\overline{u}
\right),
\]
and for $(0,b)\in {\bf R}^2$
\[
(0,b)\cdot (\overline{u},\overline{w},\overline{v})
:=\left(
\overline{u},b\sharp \overline{w},
\overline{u}
\right).
\]
Note that the orbit is diffeomorphic to
$\hat{\mathcal{M}}^0_{II'}(p_-,\gamma _-)\times 
\hat{\mathcal{M}}^0_{IV}(\gamma _-,\gamma _+)\times
\hat{\mathcal{M}}^0_{II}(\gamma _+,p_+)
\times {\bf R}^2$. 
We choose a compact set $S$ in the orbit,
and we will construct a gluing map
$\sharp :S\times [\rho _0,\infty )\to \mathcal{M}^2_I(p_-,p_+)$.
We consider a concave end which is isomorphic to 
$(-\infty ,R]\times M$.
Fix $(\overline{u},\overline{w},\overline{v})\in S$.
For simplicity, 
for $\overline{u}$ %\in \hat{\mathcal{M}}^0_{II'}(p_-,\gamma _-)$
we assume $c_+=0$ and $\tau _+=\tau _0$, 
where $c_+$ and $\tau _0 $ are constants as in the decay 
condition ${\bf (2')}$,
and we denote $(\alpha , \theta _+,z_+)$ as in ${\bf (2')}$ 
by $(\alpha _u,\theta _+^u,z_+^u )$.
Also for 
$\overline{v}$ %\in \hat{\mathcal{M}}^0_{II}(\gamma _+,p_+)$
we assume $c_-=0$ and $\tau _-=\tau _0$, as in ${\bf (2)}$,
and we denote 
$(\alpha , \theta _-,z_-)$ as in ${\bf (2)}$ 
by $(\alpha _v,\theta _-^v,z_-^v )$.
Moreover, for 
$\overline{w}$ %\in \hat{\mathcal{M}}^0_{IV}(\gamma _-,\gamma _+)$
we assume $c_+=c_-=0$ and
$\tau _+=-\tau _-=\tau _0$,
and we denote $(\alpha , \theta _-,z_-)$ %as in ${\bf (2')}$ 
by $(\alpha _w,\theta _-^w,z_-^w )$
and
$(\alpha , \theta _+,z_+)$ %as in ${\bf (2)}$ 
by $(\alpha _w,\theta _+^w,z_+^w )$.
We define $w_{\chi }(\tau ,t)\in 
L^p_1({\bf R}\times [0,1];X,L_0,L_1)$ for 
$\chi :=(\overline{u},\overline{w}, \overline{v},\rho )\in
S\times [\rho _0,\infty )$ by $w_{\chi }(\tau ,t):=$
%%%%%%%%%%%%%%%%%%%%%%%%%%%%%%%%%%%%%%%%%%%%%%%%%%%%%%%%%%%%%%%%%%
\begin{figure}[h]
\setlength{\unitlength}{1mm}
\begin{picture}(0,150)(0,0)
%%%% waku
%\put(0,0){\framebox(130,140){}}
%\put(65,0){\line(0,1){140}}
%\put(0,60){\line(1,0){130}}
%\put(0,70){\line(1,0){130}}
%\put(0,80){\line(1,0){130}}
%%%% hidari no tsutsu
\put(82.5,0){\line(0,1){135}}
\bezier{50}(82.45,134.98)(82.45,140)(81,140)
\put(112.5,0){\line(0,1){135}}
\bezier{50}(112.4,134.98)(112.4,140)(114,140)
%%%% hidari no ue wakka
\bezier{170}(82.5,80)(82.5,75)(97.5,75)
\bezier{170}(97.5,75)(112.5,75)(112.5,80)
\bezier{20}(82.5,80)(82.5,85)(97.5,85)
\bezier{20}(97.5,85)(112.5,85)(112.5,80)
%%%% hidari no shita wakka
\bezier{170}(82.5,60)(82.5,55)(97.5,55)
\bezier{170}(97.5,55)(112.5,55)(112.5,60)
\bezier{20}(82.5,60)(82.5,65)(97.5,65)
\bezier{20}(97.5,65)(112.5,65)(112.5,60)
%%%% migiue no tsutsu
\put(17.5,80){\line(0,1){55}}
\bezier{50}(17.45,134.98)(17.45,140)(16,140)
\put(47.5,80){\line(0,1){55}}
\bezier{50}(47.4,134.98)(47.4,140)(49,140)
%%%% migiue no  wakka
\bezier{170}(82.5,125)(82.5,120)(97.5,120)
\bezier{170}(97.5,120)(112.5,120)(112.5,125)
\bezier{20}(82.5,125)(82.5,130)(97.5,130)
\bezier{20}(97.5,130)(112.5,130)(112.5,125)
%%%% hidariue ue no \Lambda tachi
\put(86,124){\line(2,1){9}}
\put(109,126){\line(-2,-1){9}}
%%%% hidariue no ue wakka
\bezier{170}(17.5,125)(17.5,120)(32.5,120)
\bezier{170}(32.5,120)(47.5,120)(47.5,125)
\bezier{20}(17.5,125)(17.5,130)(32.5,130)
\bezier{20}(32.5,130)(47.5,130)(47.5,125)
%%%% hidariue no shita wakka
\bezier{170}(17.5,95)(17.5,90)(32.5,90)
\bezier{170}(32.5,90)(47.5,90)(47.5,95)
\bezier{20}(17.5,95)(17.5,100)(32.5,100)
\bezier{20}(32.5,100)(47.5,100)(47.5,95)
%%%% migi no wakka
\bezier{170}(82.5,70)(82.5,65)(97.5,65)
\bezier{170}(97.5,65)(112.5,65)(112.5,70)
\bezier{20}(82.5,70)(82.5,75)(97.5,75)
\bezier{20}(97.5,75)(112.5,75)(112.5,70)
%%%% hidari no wakka
\bezier{170}(17.5,70)(17.5,65)(32.5,65)
\bezier{170}(32.5,65)(47.5,65)(47.5,70)
\bezier{170}(17.5,70)(17.5,75)(32.5,75)
\bezier{170}(32.5,75)(47.5,75)(47.5,70)
%%%% hidari no \Lambda tachi with Reeb chords
\put(21,69){\line(2,1){9}}
\put(44,71){\line(-2,-1){9}}
\put(36,67){\vector(-4,1){12.8}}
\put(29,73){\vector(4,-1){12.8}}
%%%% hidari no Reeb chords no namae
\put(23,65){\makebox(5,5){$\gamma _-$}}
\put(37,70){\makebox(5,5){$\gamma _+$}}
%%%% migishita no tsutsu
\put(17.5,0){\line(0,1){60}}
\put(47.5,0){\line(0,1){60}}
%%%% hidarishita no ue wakka
\bezier{170}(17.5,45)(17.5,40)(32.5,40)
\bezier{170}(32.5,40)(47.5,40)(47.5,45)
\bezier{20}(17.5,45)(17.5,50)(32.5,50)
\bezier{20}(32.5,50)(47.5,50)(47.5,45)
%%%% hidarishita no shita wakka
\bezier{170}(17.5,15)(17.5,10)(32.5,10)
\bezier{170}(32.5,10)(47.5,10)(47.5,15)
\bezier{20}(17.5,15)(17.5,20)(32.5,20)
\bezier{20}(32.5,20)(47.5,20)(47.5,15)
%%%% yajirushi
\put(60,65){\makebox(10,10){$\Longrightarrow$}}
%%%% migi ue \Lambda tachi with Reeb chords
\put(86,79){\line(2,1){9}}
\put(109,81){\line(-2,-1){9}}
\put(88,80){\line(5,-1){13.53}}
\put(107,80){\line(-5,1){13.53}}
%%%% migi shita \Lambda tachi with Reeb chods
\put(86,59){\line(2,1){9}}
\put(109,61){\line(-2,-1){9}}
\put(88,60){\line(5,-1){13.53}}
\put(107,60){\line(-5,1){13.53}}
%%%% migi no \Lambda tachi with Reeb chords
\put(86,69){\line(2,1){9}}
\put(109,71){\line(-2,-1){9}}
\put(88,70){\line(5,-1){13.53}}
\put(107,70){\line(-5,1){13.53}}
%%%% hidariue ue no \Lambda tachi
\put(21,124){\line(2,1){9}}
\put(44,126){\line(-2,-1){9}}
%%%% hidariue shita no \Lambda tachi
\put(21,94){\line(2,1){9}}
\put(44,96){\line(-2,-1){9}}
%%%% hidarishita ue no \Lambda tachi
\put(21,44){\line(2,1){9}}
\put(44,46){\line(-2,-1){9}}
%%%% hidarishita shita no \Lambda tachi
\put(21,14){\line(2,1){9}}
\put(44,16){\line(-2,-1){9}}
%%%%
\put(10.5,68){\makebox(10,10){$\Lambda ^-_0$}}
\put(44.5,62){\makebox(10,10){$\Lambda ^-_1$}}
%%%% migishita no strip
\bezier{400}(88,60)(88,5)(91.25,5)
\bezier{50}(93.5,62)(93.5,5)(91.25,5)
\bezier{400}(101,57.4)(101,2.4)(104.25,2.4)
\bezier{400}(106.9,59.4)(106.5,2.4)(104.25,2.4)
\put(91.25,5){\line(5,-1){13}}
%%%% hidarishita no strip
\bezier{400}(23,55)(23,5)(26.25,5)
\bezier{50}(28.5,57)(28.5,5)(26.25,5)
\bezier{400}(36,52.4)(36,2.4)(39.25,2.4)
\bezier{400}(41.5,54.4)(41.5,2.4)(39.25,2.4)
\put(26.25,5){\line(5,-1){13}}
%%%% migiue no strips
\put(88,80){\line(0,1){50}}
\put(101,77.4){\line(0,1){50}}
\put(107,80){\line(0,1){50}}
\multiput(94,82.6)(0,2){25}{\line(0,1){0.5}}
%%%% migiue no strips
\put(23,80){\line(0,1){50}}
\put(36,77.4){\line(0,1){50}}
\put(42,80){\line(0,1){50}}
\multiput(29,82.6)(0,2){25}{\line(0,1){0.5}}
%%%%% migi no tensen
\multiput(94,62)(0,2){15}{\line(0,1){0.5}}
\put(88,60){\line(0,1){20}}
\put(101,57.4){\line(0,1){20}}
\put(107,59.4){\line(0,1){20}}
%%%%% hidari no tensen
\multiput(23,55)(0,2){15}{\line(0,1){0.5}}
\multiput(36,52.4)(0,2){15}{\line(0,1){0.5}}
\multiput(29,57)(0,2){15}{\line(0,1){0.5}}
\multiput(42,54.4)(0,2){15}{\line(0,1){0.5}}
%%%% \infty toka
\put(114.5,123){\makebox(5,5){$R$}}
\put(114.5,-2){\makebox(5,5){$-\infty $}}
\put(114.5,57){\makebox(25,5){$-T_+T_-\rho -1$}}
\put(114.5,67){\makebox(15,5){$-T_+T_-\rho $}}
\put(114.5,78){\makebox(25,5){$-T_+T_-\rho +1$}}
\put(49.5,-2){\makebox(5,5){$-\infty $}}
\put(49.5,57){\makebox(5,5){$\infty $}}
\put(49.5,78){\makebox(5,5){$-\infty $}}
\put(49.5,123){\makebox(5,5){$R$}}
%%%% strips no namae
\put(26,134){\makebox(5,5){$\overline{u}$}}
\put(34,136){\makebox(5,5){$\overline{v}$}}
\put(29,-2.5){\makebox(5,5){$\overline{w}$}}
\put(94,-2.5){\makebox(5,5){$w_{\chi }$}}
\put(114.5,27){\makebox(49,5)
{$(\alpha _w(\tau ,t)-2T_-T_+\rho ,w(\tau ,t))$}}
\put(114.5,105){\makebox(25,5)
{$\overline{v}(\tau -2T_+\rho ,t)$}}
\put(55.5,101){\makebox(25,5)
{$\overline{u}(\tau +2T_+\rho ,t)$}}
%%%%
\end{picture}
\begin{center}
Figure 3
\end{center}
\end{figure}
%%%%%%%%%%%%%%%%%%%%%%%%%%%%%%%%%%%%%%%%%%%%%%%%%%%%%%%%%%%%%%%%%%
\begin{itemize}

\item 
$\overline{u}(\tau +2T_+\rho ,t)$, 
for $\tau \leq -T_+\rho -1/T_-$,

\item 
$(\beta (-T_-\tau -T_-T_+\rho )\alpha _u(\tau +2T_+\rho ,t)
+\{1-\beta (-T_-\tau -T_-T_+\rho )\}(-T_-\tau -2T_-T_+\rho ),
\beta (-T_-\tau -T_-T_+\rho )\theta _+^u(\tau +2T_+\rho ,1-t)
+\{1-\beta (-T_-\tau -T_-T_+\rho )\}(1-t),
\beta (-T_-\tau -T_-T_+\rho )z_+^u(\tau +2T_+\rho ,t))$, 
for $-T_+\rho -1/T_-\leq \tau \leq -T_+\rho $,

\item 
$(\beta (T_-\tau +T_-T_+\rho )(\alpha _w(\tau ,t)-2T_-T_+\rho )+
\{1-\beta (T_-\tau +T_-T_+\rho )\}(-T_-\tau -2T_-T_+\rho ),
\beta (T_-\tau +T_-T_+\rho )\theta _-^w(\tau ,1-t)
+\{1-\beta (T_-\tau +T_-T_+\rho )\}(1-t),
\beta (T_-\tau +T_-T_+\rho )z_-^w(\tau ,t))$,
for $-T_+\rho \leq \tau \leq -T_+\rho +1/T_-$,

\item 
$(\alpha _w(\tau ,t)-2T_-T_+\rho ,w(\tau ,t))$,
for $-T_+\rho +1/T_-\leq \tau \leq T_-\rho -1/T_+$,

\item 
$(\beta (-T_+\tau +T_-T_+\rho )(\alpha _w(\tau ,t)-2T_-T_+\rho )
+\{1-\beta (-T_+\tau +T_-T_+\rho )\}(T_+\tau -2T_-T_+\rho ),
\beta (-T_+\tau +T_-T_+\rho )\theta _+^w(\tau ,t)
+\{1-\beta (-T_+\tau +T_-T_+\rho )\}t,
\beta (-T_+\tau +T_-T_+\rho )z_+^w(\tau ,t))$,
for $T_-\rho -1/T_+\leq \tau \leq T_-\rho $,

\item 
$(\beta (T_+\tau -T_-T_+\rho )\alpha _v(\tau -2T_-\rho ,t)
+\{1-\beta (T_+\tau -T_-T_+\rho )\}(T_+\tau -2T_-T_+\rho ),
\beta (T_+\tau -T_-T_+\rho )\theta _-^v(\tau -2T_-\rho ,t)
+\{1-\beta (T_+\tau -T_-T_+\rho )\}t,
\beta (T_+\tau -T_-T_+\rho )z_-^v(\tau -2T_-\rho ,t))$, 
for $T_-\rho \leq \tau \leq T_-\rho +1/T_+$,

\item 
$\overline{v}(\tau -2T_-\rho ,t)$
for $T_-\rho +1/T_+\leq \tau $,

\end{itemize}
see Figure 3.
Then $\overline{\partial }_{J_t}w_{\chi }=0$
for $\tau \leq -T_+\rho -1/T_-$,
$-T_+\rho +1/T_-\leq \tau \leq T_-\rho -1/T_+$ 
and $T_-\rho +1/T_+\leq \tau $, and there are constants 
$C_1>0$ and $d>0$ such that 
\begin{eqnarray}
\|\overline{\partial }_{J_t}w_{\chi }\|_{L^p}
\leq C_1e^{-d\rho },
\label{11}
\end{eqnarray}
where $C_1$ depends only on $S$ and $\rho _0$,
and $N_{w_{\chi }}$ satisfies 
\begin{eqnarray}
\|N_{w_{\chi }}(\xi )-N_{w_{\chi }}(\xi ' )\|_{L^p}
\leq C_2(\|\xi \|_{L^p_1}+\|\xi '\|_{L^p_1})
\|\xi -\xi '\|_{L^p_1},
\label{12}
\end{eqnarray}
where $\|\xi \|_{L^p_1}\leq c$ and $\|\xi '\|_{L^p_1}\leq c$ 
and $C_2$ is a constant depending only on 
$\|\nabla w_{\chi }\|_{L^p}$ and $c$.
(The proof is similar to that of (\ref{6}).)

{}From a spectral flow we can conclude
\[
\mbox{Ind}E_{w_{\chi }}
=
\mbox{Ind}E_{\overline{u};\sigma }
+\dim \mbox{Ker}Q_{\infty }
+\mbox{Ind}E_{\overline{w};\sigma }
+\dim \mbox{Ker}Q'_{\infty }
+\mbox{Ind}E_{\overline{v};\sigma},
\]
where $E_{\overline{u}}=\frac{\partial }{\partial \tau}-Q_{\tau}$
and $E_{\overline{v}}=\frac{\partial }{\partial \tau}-Q'_{\tau}$.
In Section 2 we know 
$\dim \mbox{Ker}Q_{\infty }=\dim \mbox{Ker}Q'_{\infty }=1$,
then 
$\mbox{Ind}E_{w_{\chi }}
=
\mbox{Ind}E_{\overline{u};\sigma }
+\mbox{Ind}E_{\overline{w};\sigma }
+\mbox{Ind}E_{\overline{v};\sigma}+2$.
For simplicity, in the following, we assume that 
$\mbox{Ind}E_{\overline{u};\sigma }
=\mbox{Ind}E_{\overline{w};\sigma }
=\mbox{Ind}E_{\overline{v};\sigma}=0$ (and hence 
$\mbox{Ind}E_{w_{\chi }}=2$).
We will introduce the following 
two sections $e ^0_{\rho }$ and $e^1_{\rho}$
of $w_{\chi }^*TX$:
\begin{itemize}
\item
Note that 
$\frac{d}{da}(a*\overline{u})|_{a=0}$ is a section of
$\overline{u}^*TX$ which, as $\tau \to \infty $,
is close to $T_-(\frac{\partial }{\partial \alpha})$,
and 
$\frac{d}{da}(a*\overline{w})|_{a=0}
-T_+\frac{d}{db}(b\sharp \overline{w})|_{b=0}$ 
is a section of
$\overline{w}^*TX$ which, as $\tau \to -\infty $,
is close to $(T_-+T_+)(\frac{\partial }{\partial \alpha})$
and, as $\tau \to \infty $, 
is close to $0$.
In a similar way to construct $w_{\chi }$, we glue 
$\beta (\tau -\tau _0+1)%(T_-+T_+)
\frac{d}{da}(a*\overline{u})|_{a=0}$
and 
$\frac{1}{T_-+T_+}\{\frac{d}{da}(a*\overline{w})|_{a=0}-T_+
\frac{d}{db}(b\sharp \overline{w})|_{b=0}\}$
and the zero-section of $\overline{v}^*TX$ 
to construct $e^0_{\rho}$.
Then $E_{w_{\chi }}e^0 _{\rho }=0$ 
for $\tau \leq -2T_+\rho +\tau _0-1$
$-2T_+\rho +\tau _0\leq \tau \leq -T_+\rho -1/T_-$,
$-T_+\rho +1/T_-\leq \tau \leq T_-\rho -1/T_+$ 
and $T_-\rho +1/T_+\leq \tau $.

\item
Note that 
$\frac{d}{da}(a*\overline{w})|_{a=0}
+T_-\frac{d}{db}(b\sharp \overline{w})|_{b=0}$ 
is a section of
$\overline{w}^*TX$ which, as $\tau \to -\infty $,
is close to $0$
and, as $\tau \to \infty $,
is closed 
$-(T_-+T_+)(\frac{\partial }{\partial \alpha})$.
In a similar way to construct $w_{\chi }$, we glue 
the zero-section of $\overline{u}^*TX$ and 
$-\frac{1}{T_-+T_+}\{\frac{d}{da}(a*\overline{w})|_{a=0}+
T_-\frac{d}{db}(b\sharp \overline{w})|_{b=0}\}$
and 
$\beta (-\tau -\tau _0+1)
\frac{d}{da}(a*\overline{v})|_{a=0}$
to construct $e^1_{\rho}$.
Then $E_{w_{\chi }}e^1_{\rho }=0$ 
for $\tau \leq -T_+\rho -1/T_-$,
$-T_+\rho +1/T_-\leq \tau \leq T_-\rho -1/T_+$, 
$T_-\rho +1/T_+\leq \tau \leq 2T_-\rho -\tau _0$
and $2T_-\rho -\tau _0+1\leq \tau $.

\end{itemize} From the construction, we can conclude that
$e^0_{\rho }$ is close to $(\frac{\partial }{\partial \alpha})$
on $H_-:=[\tau _0-2T_+\rho ,-\tau _0]$ 
and $0$ on $H_+:=[\tau _0,-\tau _0+2T_-\rho ]$,
and $e^1 _{\rho }$ is close to $0$ on $H_-$ 
and $(\frac{\partial }{\partial \alpha})$ on $H_-$.
Hence
\begin{eqnarray*}
&&\lim_{\rho \to \infty }
\frac{\|e^0_{\rho }\|_{H;L^2}}{\|e^0_{\rho }\|_{L^2}}
= 1,
\ 
\lim_{\rho \to \infty }
\frac{\|e^1 _{\rho }\|_{H;L^2}}{\|e^1 _{\rho }\|_{L^2}}
= 1,
\\
&&\lim_{\rho \to \infty }\int_{{\bf R}\times [0,1]}
\left(
e^0 _{\rho }/\|e^0_{\rho }\|_{L^2},
e^1 _{\rho }/\|e^1_{\rho }\|_{L^2}
\right)
d\tau dt
= 0.
\end{eqnarray*}
\begin{pro}\label{pro4.8}%%%%%%%%%%%%%%%%%%%%%%%%%%%%%%%%%%%%%%%%
Let $\chi _i:=
(\overline{u}_i,\overline{w}_i,\overline{v}_i,\rho _i)$,
$i=1,2,\ldots $, 
be a sequence of $S\times [\rho _0,\infty )$
such that $\rho _i\to \infty $.
If $\|\xi _i\|_{L^p_1}=1$ and 
$\|E_{w_{\chi _i}}\xi _i\|_{L^p}\to 0$, then 
there exists a subsequence $\{\xi _{i_l}\}$
such that
\[
\|\xi _{i_l} \|_{[-T_+\rho +R_1,-T_+\rho +R_2];L^p_1}\to 0,
\]
where $R_1\leq 0\leq R_2$.
\end{pro}%%%%%%%%%%%%%%%%%%%%%%%%%%%%%%%%%%%%%%%%%%%%%%%%%%%%%%%%%
{\it Proof}.
The proof is similar to that of Proposition \ref{pro4.3}.
Define $\xi '_i(\tau ,t):=\xi _i(\tau -T_+\rho ,t)$.
By Rellich's theorem, there exists 
$\xi' _{(R_1,R_2)}\in L^p([R_1,R_2]\times [0,1];
E\oplus\gamma _-^*\xi ,
{\bf R}(\frac{\partial }{\partial \alpha })
\oplus T_{\gamma _-(0)}\Lambda ^-_0,
{\bf R}(\frac{\partial }{\partial \alpha })
\oplus T_{\gamma _-(1)}\Lambda ^-_1)
$
and a subsequence $\{\xi '_{i_l}\}$ which converges
$\xi '_{(R_1,R_2)}$ in the norm 
$\|\cdot \|_{[R_1,R_2];L^p}$.
For simplicity, 
we assume $\|\xi '_{(R_1,R_2)}-\xi '_i\|_{L^p}\to 0$
and $\xi '_{(R'_1,R'_2)}|_{[R_1,R_2]\times [0,1]}
=\xi '_{(R_1,R_2)}$ for $[R_1,R_2]\subset [R'_1,R'_2]$.
By the G\"{a}rding's inequality,
we can conclude that $\{\xi '_i\}$ has a subsequence 
which converges to $\xi '_{(R_1,R_2)}$ in the norm 
$\|\cdot \|_{[R_1,R_2];L^p_1}$.
Moreover, from $\|\xi '_i\|_{L^p_1}\leq C$, 
we can conclude $\|\xi '_{\infty }\|_{L^p_1}\leq C$,
where $\xi '_{\infty }|_{[R_1,R_2]\times [0,1]}
:=\xi '_{(R_1,R_2)}$.
On the other hand, 
$\xi '_{\infty }$ splits into 
$E$ component $(\xi '_{\infty })_E$
and $\xi $ component $(\xi '_{\infty })_{\xi }$,
and
\begin{eqnarray*}
\frac{\partial }{\partial \tau }(\xi '_{\infty })_E
+\sqrt{-1}\frac{\partial }{\partial t}(\xi '_{\infty })_E
&=&0,
\\
\frac{\partial }{\partial \tau }(\xi '_{\infty })_{\xi }
+J_0\frac{\partial }{\partial t}(\xi '_{\infty })_{\xi }
-J_0\frac{\partial A}{\partial t}A^{-1}(\xi '_{\infty })_{\xi }
&=&0,
\end{eqnarray*}
where $-J_0\frac{\partial A}{\partial t}A^{-1}$ is the 
symmetric.
Hence $\xi '_{\infty }=0$ (see \cite{s}), and then
$\|\xi '_{i_l}\|_{[R_1,R_2];L^p_1}\to 0$.
\qed
\\

Let $W_{w_{\chi }}^{\perp }$ be the $L^2$-orthogonal compliment
of 
$W_{w_{\chi }}:=\langle e^0 _{\rho },e^1 _{\rho }\rangle$
in 
$L^p_1(w_{\chi }^*TX,w_{\chi }^*TL_0,w_{\chi }^*TL_1)$.
\begin{pro}\label{pro4.9}%%%%%%%%%%%%%%%%%%%%%%%%%%%%%%%%%%%%%%%%%
There exist constants $\rho _0$ and $C_3$ such that 
for $\chi \in S\times [\rho _0,\infty )$
and $\xi \in W_{w_{\chi }}^{\perp }$
\[
\|\xi \|_{L^p_1}\leq C_3\rho ^{\frac{3}{2}-\frac{1}{p}}
\|E_{w_{\chi }}\xi \|_{L^p}.
\]
\end{pro}%%%%%%%%%%%%%%%%%%%%%%%%%%%%%%%%%%%%%%%%%%%%%%%%%%%%%%%%%
{\it Proof}.
Let $\chi _i:=
(\overline{u}_i,\overline{w}_i,\overline{v}_i,\rho _i), 
i=1,2,\ldots $, be a sequence such that $\rho _i\to \infty $
and there exist $\xi _i\in W_{w_{\chi }}^{\perp }$
satisfying
$\|\xi _i\|_{L^p_1}=1$ and
$\rho _i^{\frac{3}{2}-\frac{1}{p}}\|E_{w_{\chi _i}}\xi _i\|_{L^p}
\to 0$.  
We define the following smooth functions:
\begin{itemize}
\item
$\beta _i(\tau):=\beta (-T_-\tau +T_+T_-\rho _i)$,
\item
$\beta '_i(\tau):=\beta (T_-\tau +T_+T_-\rho )$,
\item
$\beta ''_i(\tau):=\beta (-T_+\tau -T_+T_-\rho )$,
\item
$\beta '''_i(\tau):=\beta (T_+\tau -T_+T_-\rho )$.
\end{itemize}
Then 
\begin{eqnarray*}
\|\xi _i\|_{L^p_1}
&\leq & 
\|\beta _i\xi _i\|_{L^p_1}
+
\|\xi _i\|_{[-T+\rho -1/T_-,-T_+\rho +1/T_-];L^p_1}
+
\|\beta '_i\beta ''_i\xi _i\|_{L^p_1}
\\
&& 
+
\|\xi _i\|_{[T_-\rho -1/T_+,T_-\rho +1/T_+];L^p_1}
+
\|\beta '''\xi _i\|_{L^p_1}.
\end{eqnarray*} From Proposition \ref{pro4.8}, we can conclude
$\|\xi _i\|_{[-T+\rho -1/T_-,-T_+\rho +1/T_-];L^p_1}\to 0$
and
$\|\xi _i\|_{[-T+\rho -1/T_-,-T_+\rho +1/T_-];L^p_1}\to 0$.
Since we assume $\mbox{Ker}E_{\overline{u};\sigma }=0$,
there exist constants $C$ such that
$\|\beta _i\xi _i\|_{L^p_1}\leq 
C\|E_{\overline{u};\sigma }(\beta _i\xi _i)\|_{L^p}$,
and then
\begin{eqnarray*}
\|\beta _i\xi _i\|_{L^p_1}
&\leq &
C\|E_{\overline{u};\sigma }(\beta _i\xi _i)\|_{L^p}
\\
&\leq &
C\|E_{\overline{u}}(\beta _i\xi _i)-\frac{d\sigma }{d\tau }
\pi _E(\beta _i\xi _i)\|_{L^p}
\\
&= &
C\|E_{w_{\chi _i}}(\beta _i\xi _i)-\frac{d\sigma }{d\tau }
\pi _E(\beta _i\xi _i)\|_{L^p}
\\
&\leq &
C\|\beta _iE_{w_{\chi _i}}\xi _i\|_{L^p}
+C\|\dot{\beta }_i\xi _i\|_{L^p}
+C'\|\pi _E\xi _i\|_{[\tau _0-2T+\rho ,-\tau _0];L^p}
\end{eqnarray*}
From the assumption 
$\rho ^{\frac{3}{2}-\frac{1}{p}}\|E_{w_{\chi _i}}
\xi _i\|_{L^p}\to 0$
we know $\|\beta _iE_{w_{\chi _i}}\xi _i\|_{L^p}\to 0$,
and from Proposition \ref{pro4.8}
we conclude $\|\dot{\beta }_i\xi _i\|_{L^p}\to 0$.
Regard $\pi _E\xi _i$ on 
$[\tau _0-2T_+\rho ,-\tau _0]\times [0,1]$
as a function on $T_i:={\bf R}/2(T_+\rho -\tau _0)\times [0,1]$,
and split $\pi _E\xi _i$ into $n_i\in \mbox{Ker}E_{w_{\chi _i}}$ 
and $k_i\in (\mbox{Ker}E_{w_{\chi _i}})^{\perp }$.
Since $\xi _i\in W_{w_{\chi _i}}^{\perp }$,
we obtain 
$\lim_{i\to 0}\langle \pi _E\xi _i,
e^0 _{\rho _i}/\|e^0 _{\rho _i}\|_{L^2}
\rangle =0$,
and then $\lim_{i\to \infty }\|n_i\|_{L^2(T_i)}= 0$.
Moreover, since $E_{w_{\chi _i}}|_E
=\frac{\partial }{\partial \tau }
+\sqrt{-1}\frac{\partial }{\partial t}$,
then
$\|k_i\|_{L^2(T_i)}\leq \frac{T_+\rho _i-\tau _0}{\pi}
\|E_{w_{\chi _i}}k_i\|_{L^2(T_i)}$.
By the H\"{o}lder's inequality 
$\|y\|_{L^2(T_i)}\leq 
(2T_+\rho _i-2\tau _0)^{\frac{1}{2}-\frac{1}{p}}\|y\|_{L^p(T_i)}$
(we need $p>2$)
and the assumption 
$\rho _i^{\frac{3}{2}-\frac{1}{p}}
\|E_{w_{\chi _i}}\xi _i\|_{L^p}\to 0$,
we can conclude
$\lim _{i\to \infty }\|k_i\|_{L^2(T_i)}= 0$,
and then $\lim_{i\to \infty }\|\pi _E\xi _i\|_{L^2(T_i)}=0$.
From the assumption $\|\xi _i\|_{L^p_1}=1$ and the Sobolev's
embedding theorem (we need $p>2$), we can conclude $|\xi _i|<c$.
Then 
$\|\pi _E\xi _i\|_{L^p(T_i)}\leq c^{1-2/p}
\|\pi _E\xi _i\|_{L^2(T_i)}^{2/p}$ (also we need $p>2$),
and $\lim_{i\to \infty }\|\pi _E\xi _i\|_{L^p(T_i)}= 0$.
Finally, we finish proving 
$\lim _{i\to \infty }\|\beta _i\xi _i\|_{L^p_1}=0$.
Similarly, we can prove 
$\lim _{i\to \infty }\|\beta '_i\beta ''_i\xi _i\|_{L^p_1}=0$
and
$\lim _{i\to \infty }\|\beta '''_i\xi _i\|_{L^p_1}=0$,
and then $\|\xi _i\|_{L^p_1}\to 0$
which is a contradiction to the assumption $\|\xi _i\|_{L^p_1}=1$.
\qed 
\begin{pro}\label{pro4.10}%%%%%%%%%%%%%%%%%%%%%%%%%%%%%%%%%%%%%%%%
There exist constants $\rho _0$ and $C_3$ such that 
for $\chi \in S\times [\rho _0,\infty )$
there exists a map 
$G_{w_{\chi }}:L^p(w_{\chi }TX)\to W_{w_{\chi }}^{\perp}$
such that 
\begin{eqnarray}
E_{w_{\chi }}G_{w_{\chi }}&=&\mbox{{\rm id}},
\nonumber
\\
\|G_{w_{\chi }}\xi \|_{L^p_1}&\leq &
C_3\rho ^{\frac{3}{2}-\frac{1}{p}}\|\xi \|_{L^p}.
\label{13}
\end{eqnarray}
\end{pro}%%%%%%%%%%%%%%%%%%%%%%%%%%%%%%%%%%%%%%%%%%%%%%%%%%%%%%%%%
{\it Proof}. From Proposition \ref{pro4.9}, if 
$\xi \in \mbox{Ker}E_{w_{\chi }}\cap W_{w_{\chi }}^{\perp}$, 
then $\xi =0$ and 
\[
\dim \mbox{Ker}E_{w_{\chi }}
\leq 
2.
\]
On the other hand, from Assumption \ref{ass3.1},
\begin{eqnarray*}
\dim \mbox{Ker}E_{w_{\chi }}
& \geq & \mbox{Ind}E_{w_{\chi }}
\\
&=&
2.
\end{eqnarray*}
Hence we obtain 
\[
\mbox{Ker}E_{w_{\chi }}\oplus W_{w_{\chi }}^{\perp}=
L^p_1(w_{\chi }^*TX,w_{\chi }^*TL_0,w_{\chi }^*TL_1),
\]
and $E_{w_{\chi }}$ is surjective. From 
the estimate of Proposition \ref{pro4.9} 
we can obtain $G_{w_{\chi }}$ as in the proposition.
\qed
\\

So far we obtain the following:
Let $S$ be a compact set in the orbit.
Then there are constants $\rho _0$ such that,
for $\chi :=(\overline{u},\overline{w}, \overline{v},\rho)\in
S\times [\rho _0,\infty )$, a map 
$
f_{w_{\chi }}:
L^p_1(w_{\chi }^*TX,w_{\chi }^*TL_0,w_{\chi }^*TL_1)\to 
L^p(w_{\chi }^*TX)$
satisfies
(\ref{11}) and (\ref{12}) and 
$Df_{w_{\chi }}(0)$
possesses a right inverse 
$G_{w_{\chi }}:L^p(w_{\chi }^*TX)\to 
L^p_1(w_{\chi }^*TX,w_{\chi }^*TL_0,w_{\chi }^*TL_1)$
satisfying (\ref{13}). From $(\ref{11})$ and $(\ref{13})$
\[
\|G_{w_{\chi }}f_{w_{\chi }}(0)\|_{L^p_1}\leq C_1C_2
\rho ^{\frac{3}{2}-\frac{1}{p}}e^{-d\rho },
\]
and from $(\ref{12})$ and $(\ref{13})$
\[
\|G_{w_{\chi }}N_{w_{\chi }}(\xi )-
G_{w_{\chi }}N_{w_{\chi }}(\xi ')\|_{L^p_1}
\leq C_2C_3\rho ^{\frac{3}{2}-\frac{1}{p}}
(\|\xi \|_{L^p_1}+\|\xi '\|_{L^p_1})
\|\xi -\xi '\|_{L^p_1}.
\]
Denote $C:=C_2C_3\rho ^{\frac{3}{2}-\frac{1}{p}}$
and choose $\rho _0$ large enough, then
\[
\|G_{w_{\chi }}f_{w_{\chi }}(0)\|_{L^p_1}\leq \frac{1}{8C}.
\] 
Then, from the Newton's method in Appendix A,
we can conclude that there are constants $\rho _0$ and 
$C'$ and  a smooth map 
\[
\sharp :S\times [\rho _0,\infty )\to 
\mathcal{M}_I^2(p_-,p_+),\
\chi \mapsto \exp _{w_{\chi }}(\xi _{\chi })
\]
with 
$\|\xi _{\chi }\|_{L^p_1}\leq 
C'\|\overline{\partial }_{J_t}w_{\chi}\|_{L^p}$.
Divide them by the ${\bf R}$ actions, then we obtain
a gluing map 
$
\hat{\sharp }:\hat{S}\times [\rho _0,\infty )\to 
\hat{\mathcal{M}}_I^1(p_-,p_+)$.
\\

The next step is to show the surjectivity of
$\sharp :S\times [\rho _0,\infty )\to 
\mathcal{M}_I^2(p_-,p_+)\cap U_{(\varepsilon ,\rho )}
(\overline{u},\overline{w}, \overline{v})$.
Let $\overline{h}$
be a map which satisfies the Lagrangian boundary conditions
and ${\bf (I)}$ and the decay conditions 
${\bf (1)}$ and ${\bf (1')}$,
and
$\overline{h}(\tau ,t)=(\alpha _h,\theta _-^h,z_-^h)
\in \mbox{Im}\iota _-$ for 
$-\tau _0 \leq \tau \leq \tau _0-2T_+\rho $
and $\overline{h}(\tau ,t)=(\alpha _h,\theta _+^h,z_+^h)
\in \mbox{Im}\iota _+$ for
$\tau _0 \leq \tau \leq -\tau _0+2T_-\rho $.
Then we define
$x_{\rho }(\tau +2T_+\rho ,t)$ by
\begin{itemize}
\item 
$\overline{h}(\tau ,t)$ for $\tau \leq -T_+\rho -1/T_-$,
\item
$
(\beta (-T_-\tau -T_-T_+\rho )\alpha _h(\tau ,t) 
+\{1-\beta (-T_-\tau -T_-T_+\rho )\}(-T_-\tau -2T_-T_+\rho ),
\beta (-T_-\tau -T_-T_+\rho )\theta ^h_-(\tau ,t)  
+\{1-\beta (-T_-\tau -T_-T_+\rho )\}(1-t),
\beta (-T_-\tau -T_-T_+\rho )z^h_-(\tau ,t))
$ for $-T_+\rho -1/T_-\leq \tau $,
\end{itemize}
$w_{\rho }(\tau ,t)$ by
\begin{itemize}
\item 
$
(\beta (T_-\tau +T_-T_+\rho )(\alpha _h(\tau ,t)+2T_-T_+\rho ),
+\{1-\beta (T_-\tau +T_-T_+\rho )\}(-T_-\tau ),
\beta (T_-\tau +T_-T_+\rho )\theta ^h_-(\tau ,1-t)
+\{1-\beta (T_-\tau +T_-T_+\rho )\}(1-t),
\beta (T_-\tau +T_-T_+\rho )\}z_-^h(\tau ,t)
)
$,
for $\tau \leq -T_+\rho +1/T_-$,
\item 
$(\alpha _h(\tau ,t)+2T_-T_+\rho ,h(\tau ,t))$,
for $-T_+\rho +1/T_-\leq \tau \leq T_-\rho -1/T_+$,
\item
$
(\beta (-T_+\tau +T_-T_+\rho )\alpha _h(\tau ,t)+
\{1-\beta (-T_+\tau +T_-T_+\rho )\}T_+\tau ,
\beta (-T_+\tau +T_-T_+\rho )\theta _+^h+
\{1-\beta (-T_+\tau +T_-T_+\rho )\}t,
\beta (-T_+\tau +T_-T_+\rho )z_+^h(\tau ,t)
)$,
for $T_-\rho -1/T_+\leq \tau $, 
\end{itemize}
and $y_{\rho }(\tau -2T_+\rho ,t)$ by
\begin{itemize}
\item
$
(
\beta (T_+\tau -T_-T_+\rho )\alpha _h(\tau ,t)+
\{1-\beta (T_+\tau -T_-T_+\rho )\}(T_+\tau -2T_-T_+\rho ),
\beta (T_+\tau -T_-T_+\rho )\theta _+^h(\tau ,t)+
\{1-\beta (T_+\tau -T_-T_+\rho )\}t,
\beta (T_+\tau -T_-T_+\rho )z_+^h(\tau ,t)
)
$,
for $\tau \leq T_-\rho +1/T_+$,
\item
$\overline{h}(\tau ,t)$ for $T_-\rho +1/T_+\leq \tau $.
\end{itemize}
Moreover we define 
$U_{(\varepsilon ,\rho _0)}
(\overline{u},\overline{w},\overline{v})$,
for $(\overline{u},\overline{w},\overline{v})
\in S$, to be the set of 
$\overline{h}$ such that,
for $\rho >\rho _0$,
$\overline{u}-x_{\rho }$
satisfies the 
$e^{-\frac{1}{\varepsilon }|\tau |}$-exponential decay condition
(see (\ref{14})),
and also 
$\overline{w}-z_{\rho }$ and
$\overline{u}-y_{\rho }$.
(For simplicity,  we shall use a letter $x$ to denote $x_{\rho }$,
and also $w:=w_{\rho }$ and $y:=y_{\rho }$.) 
If $\overline{h}\in \mathcal{M}_I^2(p_-,p_+)\cap 
U_{(\varepsilon ,\rho _0)}
(\overline{u},\overline{w},\overline{v})$, then 
for a smooth map
$f_x$
there are constants $C$ and $C'$
such that
\begin{eqnarray}
\|f_x(0)\|_{L^p_{0;\sigma }}
&\leq &C
e^{-\frac{1}{\varepsilon }\rho }
\label{14}
\\
\|N_x(\xi )-N_x(\xi ')\|_{L^p_{0;\sigma }}
&\leq &C'
(\|\xi \|_{L^p_{1;\sigma }}
+\|\xi '\|_{L^p_{1;\sigma }})
\|\xi -\xi '\|_{L^p_{1;\sigma }},
\label{15}
\end{eqnarray}
where $\|\xi \|_{L^p_{1;\sigma }}\leq c$ 
and $\|\xi \|_{L^p_{1;\sigma }}\leq c$ and
$C'$ depends on $c$,
and also for $f_w$ and $f_y$.
(The proofs of the above estimates are similar to 
those of (\ref{5}) and (\ref{6}).)
We define two sections 
$e^0 _{\rho }:=(e^0 _x,e^0 _z,e^0 _y)$ and 
$e^1 _{\rho }:=(e^1 _x,e^1 _z,e^1 _y)\in 
(\langle \frac{d}{da}(a*x)|_{a=0}\rangle \oplus 
L^p_{1;\sigma }(x^*TX,x^*TL_0,x^*TL_1))
\oplus (\langle \frac{d}{da}(a*x)|_{a=0},
\frac{d}{db}(b\sharp w)|_{b=0}\rangle \oplus 
L^p_{1;\sigma }(w^*TX,w^*TL_0,w^*TL_1))
\oplus 
(\langle \frac{d}{da}(a*x)|_{a=0}\rangle \oplus 
L^p_{1;\sigma }(y^*TX,y^*TL_0,y^*TL_1))$
by 
\begin{eqnarray*}
e^0_{\rho }&:=&
\left(
\beta (\tau -\tau _0+1)\frac{d}{da}(a*x)\bigg|_{a=0},
\frac{1}{T_-+T_+}
\left\{\frac{d}{da}(a*w)\bigg|_{a=0}-T_+\frac{d}{db}(b\sharp w
)\bigg|_{b=0}\right\},
0
\right),
\\
e^1_{\rho }&:=&
\left(0,
-\frac{1}{T_-+T_+}\left\{
\frac{d}{da}(a*w)\bigg|_{a=0}
+T_-\frac{d}{db}(b\sharp w)\bigg|_{b=0}
\right\},
\beta (-\tau -\tau _0+1)
\frac{d}{da}(a*y)\bigg|_{a=0}
\right).
\end{eqnarray*}
Define $W_{(\overline{h},\rho )}:=
\langle 
e^0_{\rho },e^1_{\rho } 
\rangle $
and 
$
H_{(\overline{h},\rho )}:=W_{(\overline{h},\rho )}\oplus 
L^p_{1;\sigma }
(x^*TX,x^*TL_0,x^*TL_1)
\oplus L^p_{1;\sigma }(w^*TX,w^*TL_0,w^*TL_1)
\oplus L^p_{1;\sigma }(y^*TX,y^*TL_0,y^*TL_1)
$,
and 
the $L^2$-inner product on $H_{(\overline{h},\rho )}$ by
\begin{eqnarray*}
&&\langle \xi ,\xi '\rangle _{H_{(\overline{h},\rho )}}:=
\langle \xi ,\xi '\rangle _{L^2},
\\
&&\langle e^0_{\rho } ,\xi \rangle _{H_{(\overline{h},\rho )}}=
\langle e^1_{\rho } ,\xi \rangle _{H_{(\overline{h},\rho )}}
:=0,
\\
&&\langle e^0_{\rho } ,e^1_{\rho }\rangle 
_{H_{(\overline{h},\rho )}}:=0,
\end{eqnarray*}
where $\xi ,\xi ' \in 
L^p_{1;\sigma }
(x^*TX,x^*TL_0,x^*TL_1)
\oplus L^p_{1;\sigma }(w^*TX,w^*TL_0,w^*TL_1)
\oplus L^p_{1;\sigma }(y^*TX,y^*TL_0,y^*TL_1).
$
Let $W_{(\overline{h},\rho )}^{\perp }$ be 
the $L^2$-orthogonal compliment of $W_{(\overline{h},\rho )}$ in
$H_{(\overline{h},\rho )}$.
\begin{pro}\label{pro4.11}%%%%%%%%%%%%%%%%%%%%%%%%%%%%%%%%%%%%%%%%
There are constants $\varepsilon _0>0$ and $\rho _0$ and $C$ 
such that for 
$\overline{h}\in \mathcal{M}_I^2(p_-,p_+)\cap 
U_{(\varepsilon _0,\rho _0)}
(\overline{u},\overline{w},\overline{v})$ and 
$\xi :=(\xi _x,\xi _w, \xi _y) 
\in W_{(\overline{h},\rho )}^{\perp}$
\[
\|\xi \|_{L^p_{1;\sigma }}
\leq C%_3\rho ^{\frac{3}{2}-\frac{1}{p}}
\|(E_x\xi _x,E_w\xi _w,E_y\xi _y)\|_{L^p_{0;\sigma }}.
\]
\end{pro}%%%%%%%%%%%%%%%%%%%%%%%%%%%%%%%%%%%%%%%%%%%%%%%%%%%%%%%%%
{\it Proof}.
Let $\{\varepsilon _i \}_{i=1,2,\ldots }$ 
and $\{\rho _i\}_{i=1,2,\ldots }$ and 
$\overline{h}_i\in \mathcal{M}_I^2(p_-,p_+)\cap 
U_{(\varepsilon _i,\rho _i)}(\overline{u},\overline{w},
\overline{v})$
be sequences such that
$\varepsilon _i\to 0$ and $\rho _i\to \infty $
and
there exist $\xi _i=(\xi _{x_i},\xi _{w_i},\xi _{y_i})
\in W_{(\overline{h}_i,\rho _i)}^{\perp }$ satisfying
$\|\xi _i\|_{L^p_{1;\sigma }}=1$ and 
$
\|(E_{x_i}\xi _{x_i},E_{w_i}\xi _{w_i}, 
E_{y_i}\xi _{y_i})\|_{L^p_{0;\sigma }}\to 0$. 
Then, 
in a similar way to the proof of Proposition \ref{pro4.6}, 
we can prove that
there exists a subsequence 
$\{\xi _{i_l}\}$ such that 
$\|\xi _{i_l}\|_{L^p_{1;\sigma }}\to 0$, which contradicts the 
assumption $\|\xi _i\|_{L^p_{1;\sigma }}=1$.
\qed
\begin{pro}\label{pro4.12}%%%%%%%%%%%%%%%%%%%%%%%%%%%%%%%%%%%%%%%%
There exist constants 
$\varepsilon _0>0$ and $\rho _0$ and $C$ 
such that for 
$\overline{h}\in \mathcal{M}_I^2(p_-,p_+)\cap 
U_{(\varepsilon _0,\rho _0)}
(\overline{u},\overline{w},\overline{v})$
there exists a map 
$G_{(\overline{h},\rho )}
:L^p_{0;\sigma }(x^*TX)\oplus 
L^p_{0;\sigma }(w^*TX)\oplus L^p_{0;\sigma }(y^*TX)
\to W_{(\overline{h},\rho )}^{\perp }$ such that
\begin{eqnarray}
(E_x\oplus E_w\oplus E_y)G_{(\overline{h},\rho )} 
&=& \mbox{{\rm id}},
\nonumber
\\
\|G_{(\overline{h},\rho )}\xi \|_{L^p_{1;\sigma }} &\leq & 
C
\|\xi \|_{L^p_{0;\sigma }}.
\label{16}
\end{eqnarray}
\end{pro}%%%%%%%%%%%%%%%%%%%%%%%%%%%%%%%%%%%%%%%%%%%%%%%%%%%%%%%%%
{\it Proof}. From Proposition \ref{pro4.11}, if 
$\xi \in \mbox{Ker}(E_x\oplus E_w\oplus E_y)
\cap W_{(\overline{h},\rho )}^{\perp }$,
then $\xi =0$ and
\[
\dim \mbox{Ker}(E_x\oplus E_w\oplus E_y)\leq 
\dim W_{(\overline{h},\rho )}.
\]
On the other hand, from Assumption \ref{ass3.1},
\begin{eqnarray*}
\dim \mbox{Ker}(E_x\oplus E_w\oplus E_y)
&\geq &\mbox{Ind}E_{x;\sigma } +\mbox{Ind}E_{w;\sigma }
+\mbox{Ind}E_{y;\sigma }+2
\\
&=& \dim W_{(\overline{h},\rho )}.
\end{eqnarray*}
Hence we obtain
\[
\mbox{Ker}(E_x\oplus E_w\oplus E_y)\oplus 
W_{(\overline{h},\rho )}^{\perp}
=H_{(\overline{h},\rho )},
\]
and $E_x\oplus E_w\oplus E_y$ is surjective.
{}From the estimate of Proposition \ref{pro4.11} we can obtain
$G_{(\overline{h},\rho )}$ as in the proposition.
\qed
\\

So far we obtain the following:
There are constants $\varepsilon _0>0$ and 
$\rho _0$ such that,
for $\overline{h}\in \mathcal{M}_I^2(p_-,p_+)\cap 
U_{(\varepsilon _0,\rho _0)}
(\overline{u},\overline{w}, \overline{v})$,
a map 
$
f_{(x,w,y)}:=(f_x,f_w,f_y):
H_{(\rho ,\varepsilon )}
\to 
L^p_{0;\sigma }(x^*TX)\oplus L^p_{0;\sigma }(w^*TX)
\oplus L^p_{0;\sigma }(y^*TX)$
satisfies
(\ref{14}) and (\ref{15}), 
and $Df_{(x,w,y)}(0)$ possesses a right inverse 
$G_{(\overline{h},\rho )}:
L^p_{0;\sigma }(x^*TX)
\oplus L^p_{0;\sigma }(w^*TX)\oplus L^p_{\sigma }(y^*TX)
\to H_{(\overline{h},\rho )}$
satisfying (\ref{16}).
If we choose $\rho _0$ large enough, then 
\[
\|G_{(\overline{h},\rho )}f_{(x,w,y)}(0)\|_{L^p_{1;\sigma }}
\leq \frac{1}{8C}
\] 
{}From the Newton's method in Appendix A,
we can conclude that there are constants 
$\varepsilon >0$ and $\rho $ and
$C'$ and  a smooth map 
\begin{eqnarray*}
\sharp ':\mathcal{M}_I^2(p_-,p_+)\cap 
U_{(\varepsilon ,\rho )}
(\overline{u},\overline{w}, \overline{v})
\to  
S,\ 
\overline{h} 
\mapsto 
(\exp _x(\xi  _{\overline{h};x}),
\exp _w(\xi _{\overline{h};w}),
\exp _y(\xi _{\overline{h};y}))
\end{eqnarray*}
with 
$\|\xi _{\overline{h};x}\|_{L^p_{1;\sigma }}\leq 
C'
\|\overline{\partial }_{J_t}x\|_{L^p_{0;\sigma }}$, 
and also 
$w$ and $y$.
Divide them by the ${\bf R}$ actions, then we obtain a map
$\hat{\sharp }':\hat{\mathcal{M}}_I^1(p_-,p_+)\cap 
\hat{U}_{(\varepsilon ,\rho )}
(\overline{u},\overline{w}, \overline{v})
\to  
\hat{S}$.

{}From the construction of $\sharp $ and $\sharp '$,
if $\rho _0$ is large and $\varepsilon $ is small enough,
then $\sharp \circ \sharp '$ and $\sharp '\circ \sharp $
are diffeomorphisms.
We finish proving the gluing argument ${\bf (ii)}$.
\\

Next we will prove the gluing argument $({\bf iii})$.
(Most of the proof is similar to that of $({\bf i})$,
we will show a sketch.)
In a similar way of the proof for
${\bf (i)}$,
for compact sets $K\subset \mathcal{M}_{II}^1(\gamma _-,q)$
and $K'\subset \mathcal{M}_I^1(q,p_+)$
we construct a strip $w_{\chi }(\tau ,t)$
for $\chi :=(\overline{u},\overline{v},\rho )\in K\times K'\times
[\rho _0,\infty )$
which satisfies the Lagrangian boundary conditions 
and ${\bf (II)}$ 
and the decay conditions ${\bf (2)}$ and ${\bf (1')}$
and 
\begin{eqnarray*}
\|\overline{\partial }_{J_t}w_{\chi }\|_{L^p_{0;\sigma }}
\leq Ce^{-d\rho },
\end{eqnarray*}
where $C$ and $d>0$ are constants depending only on $K$, $K'$
and $\rho _0$.
Moreover,
for $\|\xi \|_{L^p_{1;\sigma}}\leq c$ 
and $\|\xi '\|_{L^p_{1;\sigma }}\leq c$, 
\begin{eqnarray*}
\|N_{w_{\chi }}(\xi )-N_{w_{\chi }}(\xi ')\|_{L^p_{0;\sigma }}
\leq C(\|\xi \|_{L^p_{1;\sigma }}
+\|\xi '\|_{L^p_{1;\sigma }})
\|\xi -\xi '\|_{L^p_{1;\sigma }},
\end{eqnarray*}
where $C$ is a constant depending only on 
$\|\nabla w_{\chi }\|_{L^p_{0;\sigma }}$ and $c$.
We use maps
$E_{\overline{u}}:
\langle \frac{d}{da}(a*\overline{u})|_{a=0}\rangle 
\oplus 
L^p_{1;\sigma }
(\overline{u}^*TX,\overline{u}^*TL_0,\overline{u}^*L_1)
\to 
L^p_{0;\sigma }(\overline{u}^*TX)$
and 
$E_{\overline{v}}:
L^p_1
(\overline{v}^*TX,\overline{v}^*TL_0,\overline{v}^*L_1)
\to 
L^p(\overline{v}^*TX)$.
For $\eta \in \mbox{Ker}E_{\overline{u}}$
and 
$\zeta \in \mbox{Ker}E_{\overline{v}}$,
we define similar 
$\eta \sharp _{\rho }\zeta $ as in the proof of ${\bf (i)}$.
Let $W_{w_{\chi }}^{\perp }$ be the $L^2$-
orthogonal compliment of 
$W_{w_{\chi }}
:=\{\eta \sharp _{\rho }\zeta |
\eta \in \mbox{Ker}E_{\overline{u}},
\zeta \in \mbox{Ker}E_{\overline{v}}\}
$ in
$\langle \frac{d}{da}(a*w_{\chi })|_{a=0}\rangle \oplus
L^p_{1;\sigma }(w_{\chi }^*TX,w_{\chi }^*TL_0,w_{\chi }^*TL_1)$.
\begin{pro}\label{pro4.13}%%%%%%%%%%%%%%%%%%%%%%%%%%%%%%%%%%%%%%%%
There exist constants $\rho _0$ and $C$ such that
for $\chi \in K\times K'\times [\rho _0,\infty )$
and
$\xi \in W_{w_{\chi }}^{\perp }$
\[
\|\xi \|_{L^p_{1;\sigma }}
\leq C\|E_{w_{\chi }}\xi \|_{L^p_{0;\sigma }}.
\]
\end{pro}%%%%%%%%%%%%%%%%%%%%%%%%%%%%%%%%%%%%%%%%%%%%%%%%%%%%%%%%%
\begin{pro}\label{pro4.14}%%%%%%%%%%%%%%%%%%%%%%%%%%%%%%%%%%%%%%%%
There exist constants $\rho _0$ and $C$ such that for
$\chi \in K\times K'\times [\rho _0,\infty )$
there exists a map
$G_{w_{\chi }}:L^p_{0;\sigma }\to W_{w_{\chi }}^{\perp }$
such that
\begin{eqnarray*}
E_{w_{\chi }}G_{w_{\chi }}
&=&
\mbox{{\rm id}},
\\
\|G_{w_{\chi }}\xi \|_{L^p_{1;\sigma }}
&\leq &C\|\xi \|_{L^p_{0;\sigma }}. 
\end{eqnarray*}
\end{pro}%%%%%%%%%%%%%%%%%%%%%%%%%%%%%%%%%%%%%%%%%%%%%%%%%%%%%%%%%
{}From these propositions and the Newton's method,
we can conclude that there are constants $\rho _0$ and $C$ 
and a smooth map 
\[
\sharp :K\times K'\times [\rho _0,\infty )
\to 
\mathcal{M}_{II}^2(\gamma _-,p_+),\ 
\chi \mapsto \exp _{w_{\chi}}(\xi _{\chi })
\]
with $\|\xi _{\chi }\|_{L^p_{1;\sigma }}
\leq C\|\overline{\partial }_{J_t}w_{\chi }\|_{L^p_{0;\sigma }}$.
Divide them by the ${\bf R}$ actions, we obtain 
a gluing map
$\hat{\sharp }:\hat{K}\times \hat{K}'\times [\rho _0,\infty )
\to \hat{\mathcal{M}}_{II}^1(\gamma _-,p_+)$. 
\\

The next step is to show the surjectivity of
$\sharp :K\times K'\times [\rho _0,\infty )
\to 
\mathcal{M}_{II}^2(\gamma _-,p_+)\cap U_{(\varepsilon ,\rho )}
(\overline{u},\overline{v})$. In a similar way of the 
proof for ${\bf (i)}$, for a map $w$ 
which satisfies the Lagrangian boundary conditions and 
${\bf (II)}$ and the decay conditions ${\bf (2)}$ 
and ${\bf (1')}$,
we define $x:=x_{\rho }$ and $y:=y_{\rho }$ and
$U_{(\varepsilon ,\rho _0)}(\overline{u},\overline{v})$.
If $w\in \mathcal{M}_{II}^2(\gamma _-,p_+)\cap 
U_{(\varepsilon ,\rho )}(\overline{u},\overline{v})$,
then for a smooth map $f_x$ there are constants 
$C$ and $C'$ and $d>0$ such that 
\begin{eqnarray*}
\|f_x(0)\|_{L^p_{0;\sigma }}
&\leq &Ce^{-d\rho },
\\
\|N_x(\xi )-N_x(\xi ')\|_{L^p_{0;\sigma }}
&\leq &
C'(\|\xi \|_{L^p_{1;\sigma }}+\|\xi '\|_{L^p_{1;\sigma }})
\|\xi -\xi '\|_{L^p_{1;\sigma }},
\end{eqnarray*}
where $\|\xi \|_{L^p_{1;\sigma }}\leq c$ and
$\|\xi '\|_{L^p_{1;\sigma }}\leq c$ and $C'$ depends on $c$,
and also for $f_y$ 
with norms $\|\cdot \|_{L^p_1}$ and $\|\cdot \|_{L^p}$.
For $\xi \in \mbox{Ker}\{E_w:
\langle \frac{d}{da}(a*w)|_{a=0}\rangle \oplus 
L^p_{1;\sigma }\to L^p_{0;\sigma }\}$,
we define
similar 
$\hat{\xi }:=(\xi _x,\xi _y)\in 
(\langle \frac{d}{da}(a*x)|_{a=0}\rangle \oplus L^p_{1;\sigma }
(x^*TX,x^*TL_0,x^*TL_1))\oplus 
L^p_1(y^*TX,y^*TL_0,y^*TL_1)$
as in the proof of ${\bf (i)}$.
Let $W_{(w,\rho )}^{\perp }$ be the $L^2$-orthogonal
compliment of 
$W_{(w,\rho )}:=\{\hat{\xi }|\xi \in \mbox{Ker}E_w\}$
in
$(\langle \frac{d}{da}(a*x)|_{a=0}\rangle \oplus L^p_{1;\sigma }
(x^*TX,x^*TL_0,x^*TL_1))\oplus 
L^p_1(y^*TX,y^*TL_0,y^*TL_1)$.
\begin{pro}\label{pro4.15}%%%%%%%%%%%%%%%%%%%%%%%%%%%%%%%%%%%%%%%%
There are constants $\varepsilon _0>0$ and $\rho _0$ and $C$ 
such that for $w\in \mathcal{M}_{II}^2(\gamma _-,p_+)\cap
U_{(\varepsilon _0,\rho _0)}(\overline{u},\overline{v})$
and $\xi :=(\xi _x,\xi _y)\in W_{(w,\rho )}^{\perp }$
\[
\|\xi \|_{L^p_{1;\sigma }\oplus L^p_1}
\leq C\|(E_x\xi _x,E_y\xi _y)\|_{L^p_{0;\sigma}\oplus L^p}.
\]
\end{pro}%%%%%%%%%%%%%%%%%%%%%%%%%%%%%%%%%%%%%%%%%%%%%%%%%%%%%%%%%
\begin{pro}\label{pro4.16}%%%%%%%%%%%%%%%%%%%%%%%%%%%%%%%%%%%%%%%%
There exist constants $\varepsilon _0>0$ and $\rho _0$ and $C$ 
such that for $w\in \mathcal{M}_{II}^2(\gamma _-,p_+)
\cap U_{(\varepsilon _0,\rho _0)}(\overline{u},\overline{v})$
there exists a map $G_{(w,\rho )}:L^p_{0;\sigma }(x^*TX)\oplus 
L^p(y^*TX)\to W_{(w,\rho )}^{\perp }$ such that 
\begin{eqnarray*}
(E_x\oplus E_y)G_{(w,\rho )}
&=&\mbox{{\rm id}},
\\
\|G_{(w,\rho )}\xi \|_{L^p_{1;\sigma }\oplus L^p_1}
&\leq &C\|\xi \|_{L^p_{0;\sigma }\oplus L^p}.
\end{eqnarray*}
\end{pro}%%%%%%%%%%%%%%%%%%%%%%%%%%%%%%%%%%%%%%%%%%%%%%%%%%%%%%%%%
{}From these propositions and the Newton's method,
we can conclude that there are constants $\varepsilon >0$ and 
$\rho $ and $C$ and a smooth map
\[
\sharp ':\mathcal{M}_{II}^2(\gamma _-,p_+)\cap 
U_{(\varepsilon ,\rho )}(\overline{u},\overline{v})
\to K\times K',\
w\mapsto (\exp _x(\xi _{w;x}),\exp _y(\xi _{w;y}))
\]
with $\|\xi _{w;x}\|_{L^p_{1;\sigma }}\leq C
\|\overline{\partial }_{J_t}x\|_{L^p_{0;\sigma }}$
and
$\|\xi _{w;y}\|_{L^p_1}\leq C
\|\overline{\partial }_{J_t}y\|_{L^p}$.
Divide them by the ${\bf R}$ actions,
and we obtain a map 
$\hat{\sharp }':\hat{\mathcal{M}}_{II}^1(\gamma _-,p_+)\cap 
\hat{U}_{(\varepsilon ,\rho )}(\overline{u},\overline{v})
\to \hat{K}\times \hat{K}'$.

{}From the construction of $\sharp $ and $\sharp '$,
if $\rho _0$ is large and $\varepsilon $ is small enough,
then $\sharp \circ \sharp '$ and $\sharp '\circ \sharp $
are diffeomorphisms. We finish proving the gluing
argument ${\bf (iii)}$.
\\

Next we will prove the gluing argument ${\bf (iv)}$.
(Most of the proof is similar to that of ${\bf (ii)}$,
we will show a sketch.)
Take a lift of  
$\hat{\mathcal{M}}_{III}^0(\gamma _-,\gamma _-')
\times \hat{\mathcal{M}}_{IV}^0(\gamma _-',\gamma _+')
\times \hat{\mathcal{M}}_{II}^0(\gamma _+',p_+)$
in 
$\mathcal{M}_{III}^1(\gamma _-,\gamma _-')
\times \mathcal{M}_{IV}^2(\gamma _-',\gamma _+')
\times \mathcal{M}_{II}^1(\gamma _+',p_+)$
and consider the orbit of the lift by  
${\bf R}^2$-action:
for $(a,0)\in {\bf R}^2$,
$(a,0)\cdot (\overline{u},\overline{w},\overline{v}):=
(\overline{u},a*\overline{w},\overline{v})$,
and for $(0,b)\in {\bf R}^2$,
$(0,b)\cdot (\overline{u},\overline{w},\overline{v}):=
(\overline{u},b\sharp \overline{w},\overline{v})$.
Note that the orbit is diffeomorphic to 
$\hat{\mathcal{M}}_{III}^0(\gamma _-,\gamma _-')
\times \hat{\mathcal{M}}_{IV}^0(\gamma _-',\gamma _+')
\times \hat{\mathcal{M}}_{II}^0(\gamma _+',p_+)
\times {\bf R}^2$.
We choose a compact set $S$ in the orbit, and we will construct
a gluing map $\sharp :S\times [\rho _0,\infty )\to 
\mathcal{M}_{II}^2(\gamma _-,p_+)$.
In a similar way of the proof for ${\bf (ii)}$
we construct a strip $w_{\chi }(\tau ,t)$ 
for $\chi :=(\overline{u},\overline{w},\overline{v},\rho )
\in S\times [\rho _0,\infty )$
which satisfies the Lagrangian boundary conditions 
and ${\bf (II)}$ and decay conditions 
${\bf (2)}$ and ${\bf (1')}$ and
\begin{eqnarray*}
\|\overline{\partial }_{J_t}w_{\chi }\|_{L^p_{0;\sigma }}
&\leq &Ce^{-d\rho},
\end{eqnarray*}
where $C$ and $d>0$ are constants depending only on $S$
and $\rho _0$.
Moreover, for $\|\xi \|_{L^p_{1;\sigma }}\leq c$ 
and $\|\xi '\|_{L^p_{1;\sigma }}\leq c$,
\begin{eqnarray*}
\|N_{w_{\chi }}(\xi )-N_{w_{\chi }}(\xi ')\|_{L^p_{0;\sigma }}
&\leq &
C(\|\xi \|_{L^p_{1;\sigma }}+\|\xi '\|_{L^p_{1;\sigma }})
\|\xi -\xi '\|_{L^p_{1;\sigma }},
\end{eqnarray*}
where $C$ is a constant depending only on 
$\|\nabla w_{\chi }\|_{L^p_{0;\sigma }}$ and $c$.
We define similar $e^0_{\rho }$ and $e^1_{\rho }$
as in the proof of ${\bf (ii)}$.
Let $W_{w_{\chi }}^{\perp }$ be 
the $L^2$- orthogonal compliment of
$
W_{w_{\chi }}:=
\langle -\beta (-\tau -2\tau _0-2T_-\rho )
\frac{d}{da}(a*\overline{u})|_{a=0}(\tau +2T_-\rho ,t)+
e^0 _{\rho },e^1 _{\rho }
\rangle 
$
in 
$\langle \frac{d}{da}(a*w_{\chi })|_{a=0}\rangle 
\oplus L^p_{1;\sigma }
(w_{\chi }^*TX,w_{\chi }^*TL_0,w_{\chi }^*TL_1)$.
\begin{pro}\label{pro4.17}%%%%%%%%%%%%%%%%%%%%%%%%%%%%%%%%%%%%%%%%
There exist constants $\rho _0$ and $C$ such that 
for $\chi \in S\times [\rho _0,\infty )$
and $W_{w_{\chi }}^{\perp }$
\[
\|\xi \|_{L^p_{1;\sigma }}
\leq C\rho ^{\frac{3}{2}-\frac{1}{p}}
\|E_{w_{\chi }}\xi \|_{L^p_{0;\sigma }}.
\]
\end{pro}%%%%%%%%%%%%%%%%%%%%%%%%%%%%%%%%%%%%%%%%%%%%%%%%%%%%%%%%%
\begin{pro}\label{pro4.18}%%%%%%%%%%%%%%%%%%%%%%%%%%%%%%%%%%%%%%%%
There exist constants $\rho _0$ and $C$ such that 
for $\chi \in S\times [\rho _0,\infty )$
there exists a map 
$G_{w_{\chi }}:L^p_{0;\sigma }\to W_{w_{\chi}}^{\perp } $
such that 
\begin{eqnarray*}
E_{w_{\chi }}G_{w_{\chi }}
&=&
\mbox{{\rm id}},
\\
\|G_{w_{\chi }}\xi \|_{L^p_{1;\sigma }}
&\leq &
C\rho ^{\frac{3}{2}-\frac{1}{p}}
\|\xi \|_{L^p_{0;\sigma }}.
\end{eqnarray*}
\end{pro}%%%%%%%%%%%%%%%%%%%%%%%%%%%%%%%%%%%%%%%%%%%%%%%%%%%%%%%%%
{}From these propositions and the Newton's method,
we can conclude that 
there are constants $\rho _0$ and $C$ and a smooth map 
\[
\sharp :S\times [\rho _0,\infty )\to 
\mathcal{M}_{II}^2(\gamma _-,p_+),\
\xi \mapsto \exp _{w_{\chi }}(\xi _{\chi })
\]
with $\|\xi _{\chi }\|_{L^p_{1;\sigma }}\leq 
C\|\overline{\partial }_{J_t}w_{\chi }\|_{L^p_{0;\sigma }}$.
Divide them by the ${\bf R}$ actions, we obtain a gluing map
$\hat{\sharp }:\hat{S}\times [\rho _0,\infty )\to 
\hat{\mathcal{M}}_{II}^1(\gamma _-,p_+)$.
\\

The next step is to show the surjectivity of
$\sharp :S\times [\rho _0,\infty )\to 
\mathcal{M}_{II}^2(\gamma _-,p_+)$.
In a similar way of the proof for ${\bf (ii)}$,
for a map $\overline{h}$ which satisfies the Lagrangian 
boundary conditions and ${\bf (II)}$ 
and the decay conditions ${\bf (2)}$ and ${\bf (1')}$,
we define
$x:=x_{\rho }$, $w:=w_{\rho }$ and $y:=y_{\rho }$
and $U_{(\varepsilon ,\rho _0)}
(\overline{u},\overline{w},\overline{v})$.
If $\overline{h}\in \mathcal{M}_{II}^2(\gamma _-,p_+)\cap 
U_{(\varepsilon ,\rho _0)}
(\overline{u},\overline{w},\overline{v})$,
then for a smooth map $f_x$ there are constants $C$ and $C'$
\begin{eqnarray*}
\|f_x(0)\|_{L^p_{0;\sigma }}
&\leq &Ce^{-\frac{1}{\varepsilon }\rho },
\\
\|N_x(\xi )-N_x(\xi ')\|_{L^p_{0;\sigma }}
&\leq &
C'(\|\xi \|_{L^p_{1;\sigma }}+\|\xi '\|_{L^p_{;\sigma }})
\|\xi -\xi '\|_{L^p_{1;\sigma }},
\end{eqnarray*}
where $\|\xi \|_{L^p_{1;\sigma }}\leq c$
and $\|\xi '\|_{L^p_{1;\sigma }}\leq c$
and $C'$ depend on $c$, and also for $f_w$ and $f_y$.
We define similar $e^0_{\rho }$ and $e^1_{\rho }$ 
as in the proof of ${\bf (ii)}.$
Define $W_{(\overline{h},\rho )}
:=\langle 
-\beta (-\tau -\tau _0)\frac{d}{da}(a*x)|_{a=0}
+e^0_{\rho },e^1_{\rho }
\rangle $ and 
$H_{(\overline{h},\rho )}:=W_{(\overline{h},\rho )}\oplus 
L^p_{1;\sigma }(x^*TX,x^*TL_0,x^*TL_1)
\oplus L^p_{1;\sigma }(w^*TX,w^*TL_0,w^*TL_1)
\oplus L^p_{1;\sigma }(y^*TX,y^*TL_0,y^*TL_1)$
and 
the $L^2$-inner product on $H_{(\overline{h},\rho )}$ by
\begin{eqnarray*}
&&\langle \xi ,\xi '\rangle _{H_{(\overline{h},\rho )}}
:=\langle \xi ,\xi '\rangle _{L^2},
\\
&&\left\langle 
-\beta (-\tau -\tau _0)\frac{d}{da}(a*x)\bigg|_{a=0}
+e^0_{\rho },\xi \right\rangle _{H_{(\overline{h},\rho )}}
=
\langle e^1_{\rho },
\xi \rangle _{H_{(\overline{h},\rho )}}
:=0,
\\
&&
\left\langle -\beta (\tau +\tau _0)\frac{d}{da}(a*x)\bigg|_{a=0}
+e^0 _{\rho },
e^1 _{\rho }\right\rangle _{H_{(\overline{h},\rho )}}
:=0,
\end{eqnarray*}
where $\xi ,\xi '\in 
L^p_{1;\sigma }(x^*TX,x^*TL_0,x^*TL_1)
\oplus L^p_{1;\sigma }(w^*TX,w^*TL_0,w^*TL_1)
\oplus L^p_{1;\sigma }(y^*TX,y^*TL_0,y^*TL_1)$.
Let $W_{(\overline{h},\rho )}^{\perp }$ be 
the $L^2$-orthogonal compliment of $W_{(\overline{h},\rho )}$ in
$H_{(\overline{h},\rho )}$.
\begin{pro}\label{pro4.19}%%%%%%%%%%%%%%%%%%%%%%%%%%%%%%%%%%%%%%%%
There are constants $\varepsilon _0>0$ and $\rho _0$ and $C$ 
such that 
for 
$\overline{h}\in \mathcal{M}_{II}^2(\gamma _-,p_+)\cap 
U_{(\varepsilon _0,\rho _0)}
(\overline{u},\overline{w},\overline{v})$ and 
$\xi :=(\xi _x,\xi _w, \xi _y) 
\in W_{(\overline{h},\rho )}^{\perp}$
\[
\|\xi \|_{L^p_{1;\sigma }}
\leq C
\|(E_x\xi _x,E_w\xi _w,E_y\xi _y)\|_{L^p_{0;\sigma }}.
\]
\end{pro}%%%%%%%%%%%%%%%%%%%%%%%%%%%%%%%%%%%%%%%%%%%%%%%%%%%%%%%%%
\begin{pro}\label{pro4.20}%%%%%%%%%%%%%%%%%%%%%%%%%%%%%%%%%%%%%%%%
There exist constants 
$\varepsilon _0>0$ and $\rho _0$ and $C$ 
such that 
for 
$\overline{h}\in \mathcal{M}_{II}^2(\gamma _-,p_+)\cap 
U_{(\varepsilon _0,\rho _0)}
(\overline{u},\overline{w},\overline{v})$
there exists a map 
$G_{(\overline{h},\rho )}
:L^p_{0;\sigma }(x^*TX)\oplus 
L^p_{0;\sigma }(w^*TX)\oplus L^p_{0;\sigma }(y^*TX)
\to W_{(\overline{h},\rho )}^{\perp }$ such that
\begin{eqnarray*}
(E_x\oplus E_w\oplus E_y)G_{(\overline{h},\rho )} 
&=& \mbox{{\rm id}},
\\
\|G_{(\overline{h},\rho )}\xi \|_{L^p_{1;\sigma }} &\leq & 
C
\|\xi \|_{L^p_{0;\sigma }}.
\end{eqnarray*}
\end{pro}%%%%%%%%%%%%%%%%%%%%%%%%%%%%%%%%%%%%%%%%%%%%%%%%%%%%%%%%%
{}From these propositions and the Newton's method,
we can conclude that there are constants $\varepsilon >0$
and $\rho $ and $C$ and a smooth map 
\[
\sharp ':\mathcal{M}_{II}^2(\gamma _-,p_+)
\cap U_{(\varepsilon ,\rho )}
(\overline{u},\overline{w},\overline{v})\to S,\
\overline{h}\mapsto 
(\exp _x(\xi _{\overline{h};x}),
\exp _w(\xi _{\overline{h};w}),
\exp _y(\xi _{\overline{h};y}))
\]
with $\|\xi _{\overline{h};x}\|_{L^p_{1;\sigma }}
\leq C\|\overline{\partial }_{J_t}x\|_{L^p_{0;\sigma }}$,
and also $w$ and $y$.
Divide them by the ${\bf R}$ actions,
then we obtain a map 
$\hat{\sharp }':\hat{\mathcal{M}}_{II}^1(\gamma _-,p_+)
\cap 
\hat{U}_{(\varepsilon ,\rho )}
(\overline{u},\overline{w},\overline{v})\to \hat{S}$.

{}From the construction of $\sharp $ and $\sharp '$,
if $\rho _0$ is large and $\varepsilon $ is small enough,
then $\sharp \circ \sharp '$ and $\sharp '\circ \sharp $
are diffeomorphisms. We finish proving the gluing
argument ${\bf (iv)}$.
\\

Next we will prove the gluing argument ${\bf (v)}$.
(Most of the proof is similar to that of ${\bf (ii)}$,
we will show a sketch.)
Take a lift of  
$\hat{\mathcal{M}}_V^0(\gamma _-,\gamma )
\times \hat{\mathcal{M}}_{II}^0(\gamma ,p_+)$
in 
$\mathcal{M}_V^2(\gamma _-,\gamma )
\times \mathcal{M}_{II}^1(\gamma ,p_+)$
and consider the orbit of the lift by ${\bf R}^2$-action:
for $(a,0)\in {\bf R}^2$, 
$(a,0)\cdot (\overline{u},\overline{v}):=
(a*\overline{u},\overline{v})$,
and for $(0,b)\in {\bf R}^2$,
$(0,b)\cdot (b\sharp \overline{u},\overline{v})$.
Note that the orbit is diffeomorphic to 
$\hat{\mathcal{M}}_V^0(\gamma _-,\gamma )
\times \hat{\mathcal{M}}_{II}^0(\gamma ,p_+)\times {\bf R}^2$.
We choose a compact set $S$ in the orbit, and we will construct
a gluing map $\sharp :S\times [\rho _0,\infty )\to 
\mathcal{M}_{II}^2(\gamma _-,p_+)$.
In a similar way of the proof for ${\bf (ii)}$
we construct a strip $w_{\chi }(\tau ,t)$ 
for $\chi :=(\overline{u},\overline{v},\rho )
\in S\times [\rho _0,\infty )$
which satisfies the Lagrangian boundary conditions 
and ${\bf (II)}$
and decay conditions ${\bf (2)}$ and ${\bf (1')}$ and
\begin{eqnarray*}
\|\overline{\partial }_{J_t}w_{\chi }\|_{L^p_{0;\sigma }}
&\leq &Ce^{-d\rho},
\end{eqnarray*}
where $C$ and $d>0$ are constants depending only on $S$
and $\rho _0$.
Moreover, for $\|\xi \|_{L^p_{1;\sigma }}\leq c$ 
and $\|\xi '\|_{L^p_{1;\sigma }}\leq c$,
\begin{eqnarray*}
\|N_{w_{\chi }}(\xi )-N_{w_{\chi }}(\xi ')\|_{L^p_{0;\sigma }}
&\leq &
C(\|\xi \|_{L^p_{1;\sigma }}+\|\xi '\|_{L^p_{1;\sigma }})
\|\xi -\xi '\|_{L^p_{1;\sigma }},
\end{eqnarray*}
where $C$ is a constant depending only on 
$\|\nabla w_{\chi }\|_{L^p_{0;\sigma }}$ and $c$.
We define similar $e^0_{\rho }$ and $e^1_{\rho }$ 
as in the proof of ${\bf (ii)}$.
Let $W_{w_{\chi }}^{\perp }$ be 
the $L^2$- orthogonal compliment of 
$W_{w_{\chi }}:=
\langle e^0_{\rho },e^1_{\rho }\rangle $
in 
$\langle \frac{d}{da}(a*w_{\chi })|_{a=0}\rangle 
\oplus L^p_{1;\sigma }
(w_{\chi }^*TX,w_{\chi }^*TL_0,w_{\chi }^*TL_1)$.
\begin{pro}\label{pro4.21}%%%%%%%%%%%%%%%%%%%%%%%%%%%%%%%%%%%%%%%%
There exist constants $\rho _0$ and $C$ such that 
for $\chi \in S\times [\rho _0,\infty )$
and $W_{w_{\chi }}^{\perp }$
\[
\|\xi \|_{L^p_{1;\sigma }}
\leq C\rho ^{\frac{3}{2}-\frac{1}{p}}
\|E_{w_{\chi }}\xi \|_{L^p_{0;\sigma }}.
\]
\end{pro}%%%%%%%%%%%%%%%%%%%%%%%%%%%%%%%%%%%%%%%%%%%%%%%%%%%%%%%%%
\begin{pro}\label{pro4.22}%%%%%%%%%%%%%%%%%%%%%%%%%%%%%%%%%%%%%%%%
There exist constants $\rho _0$ and $C$ such that 
for $\chi \in S\times [\rho _0,\infty )$
there exists a map 
$G_{w_{\chi }}:L^p_{0;\sigma }\to W_{w_{\chi}}^{\perp } $
such that 
\begin{eqnarray*}
E_{w_{\chi }}G_{w_{\chi }}
&=&
\mbox{{\rm id}},
\\
\|G_{w_{\chi }}\xi \|_{L^p_{1;\sigma }}
&\leq &
C\rho ^{\frac{3}{2}-\frac{1}{p}}
\|\xi \|_{L^p_{0;\sigma }}.
\end{eqnarray*}
\end{pro}%%%%%%%%%%%%%%%%%%%%%%%%%%%%%%%%%%%%%%%%%%%%%%%%%%%%%%%%%
{}From these propositions and the Newton's method,
we can conclude that 
there are constants $\rho _0$ and $C$ and a smooth map 
\[
\sharp :S\times [\rho _0,\infty )\to 
\mathcal{M}_{II}^2(\gamma _-,p_+),\
\xi \mapsto \exp _{w_{\chi }}(\xi _{\chi })
\]
with $\|\xi _{\chi }\|_{L^p_{1;\sigma }}\leq 
C\|\overline{\partial }_{J_t}w_{\chi }\|_{L^p_{0;\sigma }}$.
Divide them by the ${\bf R}$ actions, we obtain a gluing map
$\hat{\sharp }:\hat{S}\times [\rho _0,\infty )\to 
\hat{\mathcal{M}}_{II}^1(\gamma _-,p_+)$.
\\

The next step is to show the surjectivity of
$\sharp :S\times [\rho _0,\infty )\to 
\mathcal{M}_{II}^2(\gamma _-,p_+)$.
In a similar way of the proof for ${\bf (ii)}$,
for a map $\overline{h}$ which satisfies the Lagrangian 
boundary conditions and 
${\bf (II)}$ and the decay conditions 
${\bf (2)}$ and ${\bf (1')}$,
we define
$x:=x_{\rho }$ and $y:=y_{\rho }$
and $U_{(\varepsilon ,\rho _0)}
(\overline{u},\overline{v})$.
If $\overline{h}\in \mathcal{M}_{II}^2(\gamma _-,p_+)\cap 
U_{(\varepsilon ,\rho _0)}
(\overline{u},\overline{v})$,
then for a smooth map $f_x$ there are constants $C$ and $C'$
\begin{eqnarray*}
\|f_x(0)\|_{L^p_{0;\sigma }}
&\leq &Ce^{-\frac{1}{\varepsilon }\rho },
\\
\|N_x(\xi )-N_x(\xi ')\|_{L^p_{0;\sigma }}
&\leq &
C'(\|\xi \|_{L^p_{1;\sigma }}+\|\xi '\|_{L^p_{;\sigma }})
\|\xi -\xi '\|_{L^p_{1;\sigma }},
\end{eqnarray*}
where $\|\xi \|_{L^p_{1;\sigma }}\leq c$
and $\|\xi '\|_{L^p_{1;\sigma }}\leq c$
and $C'$ depend on $c$, and also for $f_y$.
Define $W_{(\overline{h},\rho )}:=
\langle 
e^0_{\rho },e^1_{\rho }
\rangle $
and 
$
H_{(\overline{h},\rho )}:=W_{(\overline{h},\rho )}\oplus 
L^p_{1;\sigma }(x^*TX,x^*TL_0,x^*TL_1)
\oplus L^p_{1;\sigma }(y^*TX,y^*TL_0,y^*TL_1)$,
and 
the $L^2$-inner product on $H_{(\overline{h},\rho )}$
in a similar way to that of ${\bf (ii)}$.
Let $W_{(\overline{h},\rho )}^{\perp }$ be 
the $L^2$-orthogonal compliment of $W_{(\overline{h},\rho )}$ in
$H_{(\overline{h},\rho )}$.
\begin{pro}\label{pro4.23}%%%%%%%%%%%%%%%%%%%%%%%%%%%%%%%%%%%%%%%%
There are constants $\varepsilon _0>0$ and $\rho _0$ and $C$ 
such that 
for 
$\overline{h}\in \mathcal{M}_{II}^2(\gamma _-,p_+)\cap 
U_{(\varepsilon _0,\rho _0)}
(\overline{u},\overline{v})$ and 
$\xi :=(\xi _x, \xi _y) 
\in W_{(\overline{h},\rho )}^{\perp}$
\[
\|\xi \|_{L^p_{1;\sigma }}
\leq C%_3\rho ^{\frac{3}{2}-\frac{1}{p}}
\|(E_x\xi _x,E_y\xi _y)\|_{L^p_{0;\sigma }}.
\]
\end{pro}%%%%%%%%%%%%%%%%%%%%%%%%%%%%%%%%%%%%%%%%%%%%%%%%%%%%%%%%%
\begin{pro}\label{pro4.24}%%%%%%%%%%%%%%%%%%%%%%%%%%%%%%%%%%%%%%%%
There exist constants 
$\varepsilon _0>0$ and $\rho _0$ and $C$ 
such that 
for 
$\overline{h}\in \mathcal{M}_{II}^2(\gamma _-,p_+)\cap 
U_{(\varepsilon _0,\rho _0)}
(\overline{u},\overline{v})$
there exists a map 
$G_{(\overline{h},\rho )}
:L^p_{0;\sigma }(x^*TX)\oplus L^p_{0;\sigma }(y^*TX)
\to W_{(\overline{h},\rho )}^{\perp }$ such that
\begin{eqnarray*}
(E_x\oplus E_y)G_{(\overline{h},\rho )} 
&=& \mbox{{\rm id}},
\\
\|G_{(\overline{h},\rho )}\xi \|_{L^p_{1;\sigma }} &\leq & 
C
\|\xi \|_{L^p_{0;\sigma }}.
\end{eqnarray*}
\end{pro}%%%%%%%%%%%%%%%%%%%%%%%%%%%%%%%%%%%%%%%%%%%%%%%%%%%%%%%%%
{}From these propositions and the Newton's method,
we can conclude that there are constants $\varepsilon >0$
and $\rho $ and $C$ and a smooth map 
\[
\sharp ':\mathcal{M}_{II}^2(\gamma _-,p_+)
\cap U_{(\varepsilon ,\rho )}
(\overline{u},\overline{v})\to S,\
\overline{h}\mapsto 
(\exp _x(\xi _{\overline{h};x}),
\exp _y(\xi _{\overline{h};y}))
\]
with $\|\xi _{\overline{h};x}\|_{L^p_{1;\sigma }}
\leq C\|\overline{\partial }_{J_t}x\|_{L^p_{0;\sigma }}$,
and also $y$.
Divide them by the ${\bf R}$ actions,
then we obtain a map 
$\hat{\sharp }':\hat{\mathcal{M}}_{II}^1(\gamma _-,p_+)
\cap 
\hat{U}_{(\varepsilon ,\rho )}
(\overline{u},\overline{v})\to \hat{S}$.

{}From the construction of $\sharp $ and $\sharp '$,
if $\rho _0$ is large and $\varepsilon $ is small enough,
then $\sharp \circ \sharp '$ and $\sharp '\circ \sharp $
are diffeomorphisms. We finish proving the gluing
argument ${\bf (v)}$.
\\

Next we will prove the gluing argument $({\bf vi})$.
(Most of the proof is similar to that of $({\bf i})$,
we will show a sketch.)
In a similar way to the proof for
${\bf (i)}$,
for compact sets $K\subset \mathcal{M}_{II}^1(\gamma _-,q)$
and $K'\subset \mathcal{M}_{II'}^1(q,\gamma _+)$
we construct a strip $w_{\chi }(\tau ,t)$
for $\chi :=(\overline{u},\overline{v},\rho )\in K\times K'\times
[\rho _0,\infty )$
which satisfies the Lagrangian boundary conditions 
and ${\bf (III)}$ 
and the decay conditions ${\bf (2)}$ and ${\bf (2')}$
and 
\begin{eqnarray*}
\|\overline{\partial }_{J_t}w_{\chi }\|_{L^p_{0;\sigma }}
\leq Ce^{-d\rho },
\end{eqnarray*}
where $C$ and $d>0$ are constants depending only on $K$, $K'$
and $\rho _0$.
Moreover,
for $\|\xi \|_{L^p_{1;\sigma}}\leq c$ 
and $\|\xi '\|_{L^p_{1;\sigma }}\leq c$, 
\begin{eqnarray*}
\|N_{w_{\chi }}(\xi )-N_{w_{\chi }}(\xi ')\|_{L^p_{0;\sigma }}
\leq C(\|\xi \|_{L^p_{1;\sigma }}
+\|\xi '\|_{L^p_{1;\sigma }})
\|\xi -\xi '\|_{L^p_{1;\sigma }},
\end{eqnarray*}
where $C$ is a constant depending only on 
$\|\nabla w_{\chi }\|_{L^p_{0;\sigma }}$ and $c$.
We use maps
$E_{\overline{u}}:
\langle \frac{d}{da}(a*\overline{u})|_{a=0}\rangle 
\oplus 
L^p_{1;\sigma }
(\overline{u}^*TX,\overline{u}^*TL_0,\overline{u}^*L_1)
\to 
L^p_{0;\sigma }(\overline{u}^*TX)$
and 
$E_{\overline{v}}:
\langle \frac{d}{da}(a*\overline{v})|_{a=0}\rangle 
\oplus
L^p_1
(\overline{v}^*TX,\overline{v}^*TL_0,\overline{v}^*L_1)
\to 
L^p(\overline{v}^*TX)$.
For $\eta \in \mbox{Ker}E_{\overline{u}}$ and 
$\zeta \in \mbox{Ker}E_{\overline{v}}\}$, we define 
similar $\eta \sharp _{\rho }\zeta $ 
as in the proof of ${\bf (i)}$.
Let $W_{w_{\chi }}^{\perp }$ be the $L^2$-
orthogonal compliment of 
$W_{w_{\chi }}:=\{\eta \sharp _{\rho }\zeta |
\eta \in \mbox{Ker}E_{\overline{u}},
\zeta \in \mbox{Ker}E_{\overline{v}}\}$
in
$\langle \frac{d}{da}(a*w_{\chi })|_{a=0}\rangle \oplus
L^p_{1;\sigma }(w_{\chi }^*TX,w_{\chi }^*TL_0,w_{\chi }^*TL_1)$.
\begin{pro}\label{pro4.25}%%%%%%%%%%%%%%%%%%%%%%%%%%%%%%%%%%%%%%%%
There exist constants $\rho _0$ and $C$ such that
for $\chi \in K\times K'\times [\rho _0,\infty )$
and
$\xi \in W_{w_{\chi }}^{\perp }$
\[
\|\xi \|_{L^p_{1;\sigma }}
\leq C\|E_{w_{\chi }}\xi \|_{L^p_{0;\sigma }}.
\]
\end{pro}%%%%%%%%%%%%%%%%%%%%%%%%%%%%%%%%%%%%%%%%%%%%%%%%%%%%%%%%%
\begin{pro}\label{pro4.26}%%%%%%%%%%%%%%%%%%%%%%%%%%%%%%%%%%%%%%%%
There exist constants $\rho _0$ and $C$ such that for
$\chi \in K\times K'\times [\rho _0,\infty )$
there exists a map
$G_{w_{\chi }}:L^p_{0;\sigma }\to W_{w_{\chi }}^{\perp }$
such that
\begin{eqnarray*}
E_{w_{\chi }}G_{w_{\chi }}
&=&
\mbox{{\rm id}},
\\
\|G_{w_{\chi }}\xi \|_{L^p_{1;\sigma }}
&\leq &C\|\xi \|_{L^p_{0;\sigma }}. 
\end{eqnarray*}
\end{pro}%%%%%%%%%%%%%%%%%%%%%%%%%%%%%%%%%%%%%%%%%%%%%%%%%%%%%%%%%
{}From these propositions and the Newton's method,
we can conclude that there are constants $\rho _0$ and $C$ 
and a smooth map 
\[
\sharp :K\times K'\times [\rho _0,\infty )
\to 
\mathcal{M}_{III}^2(\gamma _-,\gamma _+),\ 
\chi \mapsto \exp _{w_{\chi}}(\xi _{\chi })
\]
with $\|\xi _{\chi }\|_{L^p_{1;\sigma }}
\leq C\|\overline{\partial }_{J_t}w_{\chi }\|_{L^p_{0;\sigma }}$.
Divide them by the ${\bf R}$ actions, we obtain 
a gluing map
$\hat{\sharp }:\hat{K}\times \hat{K}'\times [\rho _0,\infty )
\to \hat{\mathcal{M}}_{III}^1(\gamma _-,\gamma _+)$. 
\\

The next step is to show the surjectivity of
$\sharp :K\times K'\times [\rho _0,\infty )
\to 
\mathcal{M}_{III}^2(\gamma _-,\gamma _+)
\cap U_{(\varepsilon ,\rho )}
(\overline{u},\overline{v})$. In a similar way to the 
proof for ${\bf (i)}$, for a map $w$ 
which satisfies the Lagrangian boundary conditions and 
${\bf (III)}$ and the decay conditions ${\bf (2)}$
and ${\bf (2')}$,
we define $x:=x_{\rho }$ and $y:=y_{\rho }$ and
$U_{(\varepsilon ,\rho _0)}(\overline{u},\overline{v})$.
If $w\in \mathcal{M}_{III}^2(\gamma _-,\gamma _+)\cap 
U_{(\varepsilon ,\rho )}(\overline{u},\overline{v})$,
then for a smooth map $f_x$ there are constants 
$C$ and $C'$ and $d>0$ such that 
\begin{eqnarray*}
\|f_x(0)\|_{L^p_{0;\sigma }}
&\leq &Ce^{-d\rho },
\\
\|N_x(\xi )-N_x(\xi ')\|_{L^p_{0;\sigma }}
&\leq &
C'(\|\xi \|_{L^p_{1;\sigma }}+\|\xi '\|_{L^p_{1;\sigma }})
\|\xi -\xi '\|_{L^p_{1;\sigma }},
\end{eqnarray*}
where $\|\xi \|_{L^p_{1;\sigma }}\leq c$ and
$\|\xi '\|_{L^p_{1;\sigma }}\leq c$ and $C'$ depends on $c$,
and also for $f_y$.
For $\xi \in \mbox{Ker}\{E_w:
\langle \frac{d}{da}(a*w)|_{a=0}\rangle \oplus 
L^p_{1;\sigma }\to L^p_{0;\sigma }\}$,
we define
similar 
$\hat{\xi }:=(\xi _x,\xi _y)\in 
(\langle \frac{d}{da}(a*x)|_{a=0}\rangle \oplus L^p_{1;\sigma }
(x^*TX,x^*TL_0,x^*TL_1))\oplus 
(\langle \frac{d}{da}(a*y)|_{a=0}\rangle \oplus
L^p_{1;\sigma }(y^*TX,y^*TL_0,y^*TL_1))$
as in the proof of ${\bf (i)}$.
Let $W_{(w,\rho )}^{\perp }$ be the $L^2$-orthogonal
compliment of 
$W_{(w,\rho )}:=\{\hat{\xi }|\xi \in \mbox{Ker}E_w\}$
in 
$(\langle \frac{d}{da}(a*x)|_{a=0}\rangle \oplus L^p_{1;\sigma }
(x^*TX,x^*TL_0,x^*TL_1))\oplus 
(\langle \frac{d}{da}(a*y)|_{a=0}\rangle \oplus
L^p_{1;\sigma }(y^*TX,y^*TL_0,y^*TL_1))$.
\begin{pro}\label{pro4.27}%%%%%%%%%%%%%%%%%%%%%%%%%%%%%%%%%%%%%%%%
There are constants $\varepsilon _0>0$ and $\rho _0$ and $C$ 
such that for $w\in \mathcal{M}_{III}^2(\gamma _-,\gamma _+)\cap
U_{(\varepsilon _0,\rho _0)}(\overline{u},\overline{v})$
and $\xi :=(\xi _x,\xi _y)\in W_{(w,\rho )}^{\perp }$
\[
\|\xi \|_{L^p_{1;\sigma }}
\leq C\|(E_x\xi _x,E_y\xi _y)\|_{L^p_{0;\sigma}}.
\]
\end{pro}%%%%%%%%%%%%%%%%%%%%%%%%%%%%%%%%%%%%%%%%%%%%%%%%%%%%%%%%%
\begin{pro}\label{pro4.28}%%%%%%%%%%%%%%%%%%%%%%%%%%%%%%%%%%%%%%%%
There exist constants $\varepsilon _0>0$ and $\rho _0$ and $C$ 
such that for $w\in \mathcal{M}_{III}^2(\gamma _-,\gamma _+)
\cap U_{(\varepsilon _0,\rho _0)}(\overline{u},\overline{v})$
there exists a map $G_{(w,\rho )}:L^p_{0;\sigma }(x^*TX)\oplus 
L^p_{0;\sigma }(y^*TX)\to W_{(w,\rho )}^{\perp }$ such that 
\begin{eqnarray*}
(E_x\oplus E_y)G_{(w,\rho )}
&=&\mbox{{\rm id}},
\\
\|G_{(w,\rho )}\xi \|_{L^p_{1;\sigma }}
&\leq &C\|\xi \|_{L^p_{0;\sigma }}.
\end{eqnarray*}
\end{pro}%%%%%%%%%%%%%%%%%%%%%%%%%%%%%%%%%%%%%%%%%%%%%%%%%%%%%%%%%
{}From these propositions and the Newton's method,
we can conclude that there are constants $\varepsilon >0$ and 
$\rho $ and $C$ and a smooth map
\[
\sharp ':\mathcal{M}_{III}^2(\gamma _-,\gamma _+)\cap 
U_{(\varepsilon ,\rho )}(\overline{u},\overline{v})
\to K\times K',\
w\mapsto (\exp _x(\xi _{w;x}),\exp _y(\xi _{w;y}))
\]
with $\|\xi _{w;x}\|_{L^p_{1;\sigma }}\leq C
\|\overline{\partial }_{J_t}x\|_{L^p_{0;\sigma }}$
and
$\|\xi _{w;y}\|_{L^p_{1;\sigma }}\leq C
\|\overline{\partial }_{J_t}y\|_{L^p_{0;\sigma }}$.
Divide by the ${\bf R}$ actions,
and we obtain a map 
$\hat{\sharp }':\hat{\mathcal{M}}_{III}^1(\gamma _-,\gamma _+)
\cap 
\hat{U}_{(\varepsilon ,\rho )}(\overline{u},\overline{v})
\to \hat{K}\times \hat{K}'$.

{}From the construction of $\sharp $ and $\sharp '$,
if $\rho _0$ is large and $\varepsilon $ is small enough,
then $\sharp \circ \sharp '$ and $\sharp '\circ \sharp $
are diffeomorphisms. We finish proving the gluing
argument ${\bf (vi)}$.
\\

Next we will prove the gluing argument ${\bf (vii)}$.
(Most of the proof is similar to that of ${\bf (ii)}$,
we will show a sketch.)
Take a lift of  
$\hat{\mathcal{M}}_{III}^0(\gamma _-,\gamma _-')
\times \hat{\mathcal{M}}_{IV}^0(\gamma _-',\gamma _+')
\times \hat{\mathcal{M}}_{III}^0(\gamma _+',\gamma _+)$
in 
$\mathcal{M}_{III}^1(\gamma _-,\gamma _-')
\times \mathcal{M}_{IV}^2(\gamma _-',\gamma _+')
\times \mathcal{M}_{III}^1(\gamma _+',\gamma _+)$
and consider the orbit of the lift by the similar 
${\bf R}^2$-action 
to that of ${\bf (ii)}$.
Note that the orbit is diffeomorphic to 
$\hat{\mathcal{M}}_{III}^0(\gamma _-,\gamma _-')
\times \hat{\mathcal{M}}_{IV}^0(\gamma _-',\gamma _+')
\times \hat{\mathcal{M}}_{III}^0(\gamma _+',\gamma _+)
\times {\bf R}^2$.
We choose a compact set $S$ in the orbit, and we will construct 
a gluing map $\sharp :S\times [\rho _0,\infty )\to 
\mathcal{M}_{III}^2(\gamma _-,\gamma _+)$.
In a similar way of the proof for ${\bf (ii)}$
we construct a strip $w_{\chi }(\tau ,t)$ 
for $\chi :=(\overline{u},\overline{w},\overline{v},\rho )
\in S\times [\rho _0,\infty )$
which satisfies the Lagrangian boundary conditions 
and ${\bf (III)}$
and decay conditions ${\bf (2)}$ and ${\bf (2')}$ and
\begin{eqnarray*}
\|\overline{\partial }_{J_t}w_{\chi }\|_{L^p_{0;\sigma }}
&\leq &Ce^{-d\rho},
\end{eqnarray*}
where $C$ and $d>0$ are constants depending only on $S$
and $\rho _0$.
Moreover, for $\|\xi \|_{L^p_{1;\sigma }}\leq c$ 
and $\|\xi '\|_{L^p_{1;\sigma }}\leq c$,
\begin{eqnarray*}
\|N_{w_{\chi }}(\xi )-N_{w_{\chi }}(\xi ')\|_{L^p_{0;\sigma }}
&\leq &
C(\|\xi \|_{L^p_{1;\sigma }}+\|\xi '\|_{L^p_{1;\sigma }})
\|\xi -\xi '\|_{L^p_{1;\sigma }},
\end{eqnarray*}
where $C$ is a constant depending only on 
$\|\nabla w_{\chi }\|_{L^p_{0;\sigma }}$ and $c$.
We define similar $e^0_{\rho }$ and $e^1_{\rho }$ as in the proof 
of ${\bf (ii)}$.
Let 
$W_{w_{\chi }}^{\perp }$ be the $L^2$-orthogonal compliment of
$W_{w_{\chi }}:=
\langle -\beta (-\tau -2\tau _0-2T_-\rho )
\frac{d}{da}(a*\overline{u}|_{a=0}(\tau +2T_-\rho ,t)+e^0_{\rho},
\beta (\tau -2\tau _0-2T_+\rho )
\frac{d}{da}(a*\overline{w})|_{a=0}(\tau -2T_+\rho ,t)
+e^1_{\rho}\rangle $
in 
$\langle \frac{d}{da}(a*w_{\chi })|_{a=0}\rangle 
\oplus L^p_{1;\sigma }
(w_{\chi }^*TX,w_{\chi }^*TL_0,w_{\chi }^*TL_1)$.
\begin{pro}\label{pro4.29}%%%%%%%%%%%%%%%%%%%%%%%%%%%%%%%%%%%%%%%%
There exist constants $\rho _0$ and $C$ such that 
for $\chi \in S\times [\rho _0,\infty )$
and $W_{w_{\chi }}^{\perp }$
\[
\|\xi \|_{L^p_{1;\sigma }}
\leq C\rho ^{\frac{3}{2}-\frac{1}{p}}
\|E_{w_{\chi }}\xi \|_{L^p_{0;\sigma }}.
\]
\end{pro}%%%%%%%%%%%%%%%%%%%%%%%%%%%%%%%%%%%%%%%%%%%%%%%%%%%%%%%%%
\begin{pro}\label{pro4.30}%%%%%%%%%%%%%%%%%%%%%%%%%%%%%%%%%%%%%%%%
There exist constants $\rho _0$ and $C$ such that 
for $\chi \in S\times [\rho _0,\infty )$
there exists a map 
$G_{w_{\chi }}:L^p_{0;\sigma }\to W_{w_{\chi}}^{\perp } $
such that 
\begin{eqnarray*}
E_{w_{\chi }}G_{w_{\chi }}
&=&
\mbox{{\rm id}},
\\
\|G_{w_{\chi }}\xi \|_{L^p_{1;\sigma }}
&\leq &
C\rho ^{\frac{3}{2}-\frac{1}{p}}
\|\xi \|_{L^p_{0;\sigma }}.
\end{eqnarray*}
\end{pro}%%%%%%%%%%%%%%%%%%%%%%%%%%%%%%%%%%%%%%%%%%%%%%%%%%%%%%%%%
{}From these propositions and the Newton's method,
we can conclude that 
there are constants $\rho _0$ and $C$ and a smooth map 
\[
\sharp :S\times [\rho _0,\infty )\to 
\mathcal{M}_{III}^2(\gamma _-,\gamma _+),\
\xi \mapsto \exp _{w_{\chi }}(\xi _{\chi })
\]
with $\|\xi _{\chi }\|_{L^p_{1;\sigma }}\leq 
C\|\overline{\partial }_{J_t}w_{\chi }\|_{L^p_{0;\sigma }}$.
Divide them by the ${\bf R}$ actions, we obtain a gluing map
$\hat{\sharp }:\hat{S}\times [\rho _0,\infty )\to 
\hat{\mathcal{M}}_{III}^1(\gamma _-,\gamma _+)$.
\\

The next step is to show the surjectivity of
$\sharp :S\times [\rho _0,\infty )\to 
\mathcal{M}_{III}^2(\gamma _-,\gamma _+)$.
In a similar way of the proof for ${\bf (ii)}$,
for a map $\overline{h}$ which satisfies the Lagrangian 
boundary conditions and ${\bf (III)}$ and 
the decay conditions ${\bf (2)}$ and ${\bf (2')}$,
we define
$x:=x_{\rho }$, $w:=w_{\rho }$ and $y:=y_{\rho }$
and $U_{(\varepsilon ,\rho _0)}
(\overline{u},\overline{w},\overline{v})$.
If $\overline{h}\in \mathcal{M}_{III}^2(\gamma _-,\gamma _+)\cap 
U_{(\varepsilon ,\rho _0)}
(\overline{u},\overline{w},\overline{v})$,
then for a smooth map $f_x$ there are constants $C$ and $C'$
\begin{eqnarray*}
\|f_x(0)\|_{L^p_{0;\sigma }}
&\leq &Ce^{-\frac{1}{\varepsilon }\rho },
\\
\|N_x(\xi )-N_x(\xi ')\|_{L^p_{0;\sigma }}
&\leq &
C'(\|\xi \|_{L^p_{1;\sigma }}+\|\xi '\|_{L^p_{;\sigma }})
\|\xi -\xi '\|_{L^p_{1;\sigma }},
\end{eqnarray*}
where $\|\xi \|_{L^p_{1;\sigma }}\leq c$
and $\|\xi '\|_{L^p_{1;\sigma }}\leq c$
and $C'$ depend on $c$, and also for $f_w$ and $f_y$.
We define similar $e^0_{\rho }$ and $e^1_{\rho }$ 
as in the proof ${\bf (ii)}$.
Define $W_{(\overline{h},\rho )}:=
\langle 
-\beta (-\tau -\tau _0)\frac{d}{da}(a*x)|_{a=0}+e^0_{\rho },
\beta (\tau -\tau _0)\frac{d}{da}(a*y)|_{a=0}+e^1_{\rho }
\rangle $
and 
$
H_{(\overline{h},\rho )}:=W_{(\overline{h},\rho )}\oplus 
L^p_{1;\sigma }(x^*TX,x^*TL_0,x^*TL_1)
\oplus 
L^p_{1;\sigma }(w^*TX,w^*TL_0,w^TL_1)
\oplus 
L^p_{1;\sigma }(y^*TX,y^*TL_0,y^*TL_1)
$
and 
the $L^2$-inner product on $H_{(\overline{h},\rho )}$ by
\begin{eqnarray*}
&&\langle \xi ,\xi '\rangle _{H_{(\overline{h},\rho )}}
:=\langle \xi ,\xi '\rangle _{L^2},
\\
&&
\left\langle 
-\beta (-\tau -\tau _0)\frac{d}{da}(a*x)\bigg|_{a=0}+e^0_{\rho },
\xi \right\rangle _{H_{(\overline{h},\rho )}}
:=\left\langle 
\beta (\tau -\tau _0)\frac{d}{da}(a*y)\bigg|_{a=0}+e^1_{\rho },
\xi \right\rangle _{H_{(\overline{h},\rho )}}
:=0,
\\
&&
\left\langle 
-\beta (-\tau -\tau _0)\frac{d}{da}(a*x)\bigg|_{a=0}+e^0_{\rho },
\beta (\tau -\tau _0)\frac{d}{da}(a*y)\bigg|_{a=0}+e^1_{\rho }
\right\rangle _{H_{(\overline{h},\rho )}}
:=0,
\end{eqnarray*}
where $\xi ,\xi '\in
L^p_{1;\sigma }(x^*TX,x^*TL_0,x^*TL_1)
\oplus 
L^p_{1;\sigma }(w^*TX,w^*TL_0,w^TL_1)
\oplus 
L^p_{1;\sigma }(y^*TX,y^*TL_0,y^*TL_1).
$
Let $W_{(\overline{h},\rho )}^{\perp }$ be 
the $L^2$-orthogonal compliment of $W_{(\overline{h},\rho )}$ in
$H_{(\overline{h},\rho )}$.
\begin{pro}\label{pro4.31}%%%%%%%%%%%%%%%%%%%%%%%%%%%%%%%%%%%%%%%%
There are constants $\varepsilon _0>0$ and $\rho _0$ and $C$ 
such that 
for 
$\overline{h}\in \mathcal{M}_{III}^2(\gamma _-,\gamma _+)\cap 
U_{(\varepsilon _0,\rho _0)}
(\overline{u},\overline{w},\overline{v})$ and 
$\xi :=(\xi _x,\xi _w, \xi _y) 
\in W_{(\overline{h},\rho )}^{\perp}$
\[
\|\xi \|_{L^p_{1;\sigma }}
\leq C%_3\rho ^{\frac{3}{2}-\frac{1}{p}}
\|(E_x\xi _x,E_w\xi _w,E_y\xi _y)\|_{L^p_{0;\sigma }}.
\]
\end{pro}%%%%%%%%%%%%%%%%%%%%%%%%%%%%%%%%%%%%%%%%%%%%%%%%%%%%%%%%%
\begin{pro}\label{4.32}%%%%%%%%%%%%%%%%%%%%%%%%%%%%%%%%%%%%%%%%%%%
There exist constants 
$\varepsilon _0>0$ and $\rho _0$ and $C$ 
such that 
for 
$\overline{h}\in \mathcal{M}_{III}^2(\gamma _-,\gamma _+)\cap 
U_{(\varepsilon _0,\rho _0)}
(\overline{u},\overline{w},\overline{v})$
there exists a map 
$G_{(\overline{h},\rho )}
:L^p_{0;\sigma }(x^*TX)\oplus 
L^p_{0;\sigma }(w^*TX)\oplus L^p_{0;\sigma }(y^*TX)
\to W_{(\overline{h},\rho )}^{\perp }$ such that
\begin{eqnarray*}
(E_x\oplus E_w\oplus E_y)G_{(\overline{h},\rho )} 
&=& \mbox{{\rm id}},
\\
\|G_{(\overline{h},\rho )}\xi \|_{L^p_{1;\sigma }} &\leq & 
C
\|\xi \|_{L^p_{0;\sigma }}.
\end{eqnarray*}
\end{pro}%%%%%%%%%%%%%%%%%%%%%%%%%%%%%%%%%%%%%%%%%%%%%%%%%%%%%%%%%
{}From these propositions and the Newton's method,
we can conclude that there are constants $\varepsilon >0$
and $\rho $ and $C$ and a smooth map 
\[
\sharp ':\mathcal{M}_{III}^2(\gamma _-,\gamma _+)
\cap U_{(\varepsilon ,\rho )}
(\overline{u},\overline{w},\overline{v})\to S,\
\overline{h}\mapsto 
(\exp _x(\xi _{\overline{h};x}),
\exp _w(\xi _{\overline{h};w}),
\exp _y(\xi _{\overline{h};y}))
\]
with $\|\xi _{\overline{h};x}\|_{L^p_{1;\sigma }}
\leq C\|\overline{\partial }_{J_t}x\|_{L^p_{0;\sigma }}$,
and also $w$ and $y$.
Divide them by the ${\bf R}$ actions,
then we obtain a map 
$\hat{\sharp }':\hat{\mathcal{M}}_{III}^1(\gamma _-,\gamma _+)
\cap \hat{U}_{(\varepsilon ,\rho )}
(\overline{u},\overline{w},\overline{v})\to \hat{S}$.

{}From the construction of $\sharp $ and $\sharp '$,
if $\rho _0$ is large and $\varepsilon $ is small enough,
then $\sharp \circ \sharp '$ and $\sharp '\circ \sharp $
are diffeomorphisms. We finish proving the gluing
argument ${\bf (vii)}$.
\\

Finally we will prove the gluing argument ${\bf (viii)}$.
(Most of the proof is similar to that of ${\bf (ii)}$,
we will show a sketch.)
Take a lift of  
$\hat{\mathcal{M}}_V^0(\gamma _-,\gamma )
\times \hat{\mathcal{M}}_{III}^0(\gamma ,\gamma _+)$
in 
$\mathcal{M}_V^2(\gamma _-,\gamma )
\times \mathcal{M}_{III}^1(\gamma ,\gamma _+)$
and consider the orbit of the lift by the similar 
${\bf R}^2$-action 
to that of ${\bf (v)}$.
Note that the orbit is diffeomorphic to 
$\hat{\mathcal{M}}_V^0(\gamma _-,\gamma )
\times \hat{\mathcal{M}}_{III}^0(\gamma ,\gamma _+)
\times {\bf R}^2$.
We choose a compact set $S$ in the orbit, and we will construct
a gluing map $\sharp :S\times [\rho _0,\infty )\to 
\mathcal{M}_{III}^2(\gamma _-,\gamma _+)$.
In a similar way to the proof for ${\bf (ii)}$
we construct a strip $w_{\chi }(\tau ,t)$ 
for $\chi :=(\overline{u},\overline{v},\rho )
\in S\times [\rho _0,\infty )$
which satisfies the Lagrangian boundary conditions 
and ${\bf (III)}$
and decay conditions ${\bf (2)}$ and ${\bf (2')}$ and
\begin{eqnarray*}
\|\overline{\partial }_{J_t}w_{\chi }\|_{L^p_{0;\sigma }}
&\leq &Ce^{-d\rho},
\end{eqnarray*}
where $C$ and $d>0$ are constants depending only on $S$
and $\rho _0$.
Moreover, for $\|\xi \|_{L^p_{1;\sigma }}\leq c$ 
and $\|\xi '\|_{L^p_{1;\sigma }}\leq c$,
\begin{eqnarray*}
\|N_{w_{\chi }}(\xi )-N_{w_{\chi }}(\xi ')\|_{L^p_{0;\sigma }}
&\leq &
C(\|\xi \|_{L^p_{1;\sigma }}+\|\xi '\|_{L^p_{1;\sigma }})
\|\xi -\xi '\|_{L^p_{1;\sigma }},
\end{eqnarray*}
where $C$ is a constant depending only on 
$\|\nabla w_{\chi }\|_{L^p_{0;\sigma }}$ and $c$.
We define similar $e^0_{\rho }$ and $e^1_{\rho }$
as in the proof of ${\bf (ii)}$.
Let $W_{w_{\chi }}^{\perp }$ be 
the $L^2$- orthogonal compliment of 
$W_{w_{\chi }}
:=\{e^0_{\rho },
\beta (\tau -2\tau _0-2T\rho )\frac{d}{da}(a*\overline{v})|_{a=0}
(\tau -2T\rho ,t)+e^1_{\rho }\}$
in 
$\langle \frac{d}{da}(a*w_{\chi })|_{a=0}\rangle 
\oplus L^p_{1;\sigma }
(w_{\chi }^*TX,w_{\chi }^*TL_0,w_{\chi }^*TL_1)$.
\begin{pro}\label{pro4.33}%%%%%%%%%%%%%%%%%%%%%%%%%%%%%%%%%%%%%%%%
There exist constants $\rho _0$ and $C$ such that 
for $\chi \in S\times [\rho _0,\infty )$
and $W_{w_{\chi }}^{\perp }$
\[
\|\xi \|_{L^p_{1;\sigma }}
\leq C\rho ^{\frac{3}{2}-\frac{1}{p}}
\|E_{w_{\chi }}\xi \|_{L^p_{0;\sigma }}.
\]
\end{pro}%%%%%%%%%%%%%%%%%%%%%%%%%%%%%%%%%%%%%%%%%%%%%%%%%%%%%%%%%
\begin{pro}\label{pro4.34}%%%%%%%%%%%%%%%%%%%%%%%%%%%%%%%%%%%%%%%%
There exist constants $\rho _0$ and $C$ such that 
for $\chi \in S\times [\rho _0,\infty )$
there exists a map 
$G_{w_{\chi }}:L^p_{0;\sigma }\to W_{w_{\chi}}^{\perp } $
such that 
\begin{eqnarray*}
E_{w_{\chi }}G_{w_{\chi }}
&=&
\mbox{{\rm id}},
\\
\|G_{w_{\chi }}\xi \|_{L^p_{1;\sigma }}
&\leq &
C\rho ^{\frac{3}{2}-\frac{1}{p}}
\|\xi \|_{L^p_{0;\sigma }}.
\end{eqnarray*}
\end{pro}%%%%%%%%%%%%%%%%%%%%%%%%%%%%%%%%%%%%%%%%%%%%%%%%%%%%%%%%%
{}From these propositions and the Newton's method,
we can conclude that 
there are constants $\rho _0$ and $C$ and a smooth map 
\[
\sharp :S\times [\rho _0,\infty )\to 
\mathcal{M}_{III}^2(\gamma _-,\gamma _+),\
\xi \mapsto \exp _{w_{\chi }}(\xi _{\chi })
\]
with $\|\xi _{\chi }\|_{L^p_{1;\sigma }}\leq 
C\|\overline{\partial }_{J_t}w_{\chi }\|_{L^p_{0;\sigma }}$.
Divide them by the ${\bf R}$ actions, we obtain a gluing map
$\hat{\sharp }:\hat{S}\times [\rho _0,\infty )\to 
\hat{\mathcal{M}}_{III}^1(\gamma _-,\gamma _+)$.
\\

The next step is to show the surjectivity of
$\sharp :S\times [\rho _0,\infty )\to 
\mathcal{M}_{III}^2(\gamma _-,\gamma _+)$.
In a similar way of the proof for ${\bf (ii)}$,
for a map $\overline{h}$ which satisfies the Lagrangian 
boundary conditions and ${\bf (III)}$ and the decay conditions 
${\bf (2)}$ and ${\bf (2')}$,
we define
$x:=x_{\rho }$ and $y:=y_{\rho }$
and $U_{(\varepsilon ,\rho _0)}
(\overline{u},\overline{v})$.
If $\overline{h}\in \mathcal{M}_{III}^2(\gamma _-,\gamma _+)\cap 
U_{(\varepsilon ,\rho _0)}
(\overline{u},\overline{v})$,
then for a smooth map $f_x$ there are constants $C$ and $C'$
\begin{eqnarray*}
\|f_x(0)\|_{L^p_{0;\sigma }}
&\leq &Ce^{-\frac{1}{\varepsilon }\rho },
\\
\|N_x(\xi )-N_x(\xi ')\|_{L^p_{0;\sigma }}
&\leq &
C'(\|\xi \|_{L^p_{1;\sigma }}+\|\xi '\|_{L^p_{;\sigma }})
\|\xi -\xi '\|_{L^p_{1;\sigma }},
\end{eqnarray*}
where $\|\xi \|_{L^p_{1;\sigma }}\leq c$
and $\|\xi '\|_{L^p_{1;\sigma }}\leq c$
and $C'$ depend on $c$, and also for $f_y$.
Define $W_{(\overline{h},\rho )}
:=\langle e^0_{\rho },e^1_{\rho }\rangle $
and 
$
H_{(\overline{h},\rho )}:=W_{(\overline{h},\rho )}\oplus 
L^p_{1;\sigma }
(x^*TX,x^*TL_0,x^*TL_1)
\oplus L^p_{1;\sigma }(y^*TX,y^*TL_0,y^*TL_1),
$
and 
the $L^2$-inner product on $H_{(\overline{h},\rho )}$
in a similar way to that of ${\bf (v)}$.
Let $W_{(\overline{h},\rho )}^{\perp }$ be 
the $L^2$-orthogonal compliment of $W_{(\overline{h},\rho )}$ in
$H_{(\overline{h},\rho )}$.
\begin{pro}\label{pro4.35}%%%%%%%%%%%%%%%%%%%%%%%%%%%%%%%%%%%%%%%%
There are constants $\varepsilon _0>0$ and $\rho _0$ and $C$ 
such that 
for 
$\overline{h}\in \mathcal{M}_{III}^2(\gamma _-,\gamma _+)\cap 
U_{(\varepsilon _0,\rho _0)}
(\overline{u},\overline{v})$ and 
$\xi :=(\xi _x, \xi _y) 
\in W_{(\overline{h},\rho )}^{\perp}$
\[
\|\xi \|_{L^p_{1;\sigma }}
\leq C
\|(E_x\xi _x,E_y\xi _y)\|_{L^p_{0;\sigma }}.
\]
\end{pro}%%%%%%%%%%%%%%%%%%%%%%%%%%%%%%%%%%%%%%%%%%%%%%%%%%%%%%%%%
\begin{pro}\label{pro4.36}%%%%%%%%%%%%%%%%%%%%%%%%%%%%%%%%%%%%%%%%
There exist constants 
$\varepsilon _0>0$ and $\rho _0$ and $C$ 
such that 
for 
$\overline{h}\in \mathcal{M}_{III}^2(\gamma _-,\gamma _+)\cap 
U_{(\varepsilon _0,\rho _0)}
(\overline{u},\overline{v})$
there exists a map 
$G_{(\overline{h},\rho )}
:L^p_{0;\sigma }(x^*TX)\oplus L^p_{0;\sigma }(y^*TX)
\to W_{(\overline{h},\rho )}^{\perp }$ such that
\begin{eqnarray*}
(E_x\oplus E_y)G_{(\overline{h},\rho )} 
&=& \mbox{{\rm id}},
\\
\|G_{(\overline{h},\rho )}\xi \|_{L^p_{1;\sigma }} &\leq & 
C
\|\xi \|_{L^p_{0;\sigma }}.
\end{eqnarray*}
\end{pro}%%%%%%%%%%%%%%%%%%%%%%%%%%%%%%%%%%%%%%%%%%%%%%%%%%%%%%%%%
{}From these propositions and the Newton's method,
we can conclude that there are constants $\varepsilon >0$
and $\rho $ and $C$ and a smooth map 
\[
\sharp ':\mathcal{M}_{III}^2(\gamma _-,\gamma _+)
\cap U_{(\varepsilon ,\rho )}
(\overline{u},\overline{v})\to S,\
\overline{h}\mapsto 
(\exp _x(\xi _{\overline{h};x}),
\exp _y(\xi _{\overline{h};y}))
\]
with $\|\xi _{\overline{h};x}\|_{L^p_{1;\sigma }}
\leq C\|\overline{\partial }_{J_t}x\|_{L^p_{0;\sigma }}$,
and also $y$.
Divide them by the ${\bf R}$ actions,
then we obtain a map 
$\hat{\sharp }':\hat{\mathcal{M}}_{III}^1(\gamma _-,\gamma _+)
\cap 
\hat{U}_{(\varepsilon ,\rho )}
(\overline{u},\overline{v})\to \hat{S}$.

{}From the construction of $\sharp $ and $\sharp '$,
if $\rho _0$ is large and $\varepsilon $ is small enough,
then $\sharp \circ \sharp '$ and $\sharp '\circ \sharp $
are diffeomorphisms. We finish proving the gluing
argument ${\bf (viii)}$.
\\

We observe the dimensions of the moduli spaces 
in gluing arguments.
For example, we consider the following case.  
Take a lift of 
\begin{eqnarray*}
\hat{\mathcal{M}}_{II'}^{e_0}(p_-,\gamma _-^1)
\times 
\hat{\mathcal{M}}_{IV}^{d_1}(\gamma _-^1,\gamma _+^1)
\times 
\hat{\mathcal{M}}_{III}^{e_1}(\gamma _+^1,\gamma _-^2)
\times 
\hat{\mathcal{M}}_{IV}^{d_2}(\gamma _-^2,\gamma _+^2)
\times
\cdots 
\\
\times 
\hat{\mathcal{M}}_{III}^{e_{k-1}}(\gamma _+^{k-1},\gamma _-^k)
\times 
\hat{\mathcal{M}}_{IV}^{d_k}(\gamma _-^k,\gamma _+^k)
\times
\hat{\mathcal{M}}_{II}^{e_k}(\gamma _+^k,p_+)
\end{eqnarray*}
in
\begin{eqnarray*}
\mathcal{M}_{II'}^{e_0+1}(p_-,\gamma _-^1)
\times 
\mathcal{M}_{IV}^{d_1+2}(\gamma _-^1,\gamma _+^1)
\times 
\mathcal{M}_{III}^{e_1+1}(\gamma _+^1,\gamma _-^2)
\times 
\mathcal{M}_{IV}^{d_2+2}(\gamma _-^2,\gamma _+^2)
\times
\cdots 
\\
\times 
\mathcal{M}_{III}^{e_{k-1}+1}(\gamma _+^{k-1},\gamma _-^k)
\times 
\mathcal{M}_{IV}^{d_k+2}(\gamma _-^k,\gamma _+^k)
\times
\mathcal{M}_{II}^{e_k+1}(\gamma _+^k,p_+)
\end{eqnarray*}
and consider the orbit
of the lift by the following 
${\bf R}^{2k}$-action:
\\
for $((0,0),\ldots,(a_l,0),\ldots ,(0,0))
\in {\bf R}^{2k}$
\begin{eqnarray*}
&&((0,0)\ldots ,(a_l,0),\ldots ,(0,0))\cdot 
(\overline{u}_0,\overline{w}_1,\overline{u}_1,
\ldots ,\overline{w}_k,\overline{u}_k)
\\
&:=&
(\overline{u}_0,\overline{w}_1,\overline{u}_1,
\ldots ,a*\overline{w}_l,\ldots ,\overline{w}_k,\overline{u}_k),
\end{eqnarray*}
and for $((0,0),\ldots ,(0,b_l),\ldots ,(0,0))
\in {\bf R}^{2k}$
\begin{eqnarray*}
&&
((0,0),\ldots ,(0,b_l),\ldots ,(0,0))
\cdot 
(\overline{u}_0,\overline{w}_1,\overline{u}_1,
\ldots ,\overline{w}_k,\overline{u}_k)
\\
&:=&
(\overline{u}_0,\overline{w}_1,
\overline{u}_1,
\ldots ,
b_l\sharp \overline{w}_l,\ldots 
\overline{w}_k,\overline{u}_k).
\end{eqnarray*}
Then the orbit is diffeomorphic to
\begin{eqnarray*}
\hat{\mathcal{M}}_{II'}^{e_0}(p_-,\gamma _-^1)
\times 
\hat{\mathcal{M}}_{IV}^{d_1}(\gamma _-^1,\gamma _+^1)
\times 
\hat{\mathcal{M}}_{III}^{e_1}(\gamma _+^1,\gamma _-^2)
\times 
\hat{\mathcal{M}}_{IV}^{d_2}(\gamma _-^2,\gamma _+^2)
\times
\cdots 
\\
\times 
\hat{\mathcal{M}}_{III}^{e_{k-1}}(\gamma _+^{k-1},\gamma _-^k)
\times 
\hat{\mathcal{M}}_{IV}^{d_k}(\gamma _-^k,\gamma _+^k)
\times
\hat{\mathcal{M}}_{II}^{e_k}(\gamma _+^k,p_+)
\times {\bf R}^{2k},
\end{eqnarray*}
and we can construct a smooth map
\[
\sharp :S\times [\rho _0,\infty )\to 
\mathcal{M}_I^{e_0+\cdots +d_k+e_k+2k}(p_-,p_+),
\]
where $S$ is a compact set in the orbit
and $\rho _0$ is a constant depending on $S$.
Divide them by the ${\bf R}$-actions, then we obtain 
a gluing map 
$\hat{\sharp }:\hat{S}\to 
\hat{\mathcal{M}}_I^{e_0+\cdots + d_k+e_k+2k-1}(p_-,p_+)$.
This implies that each element of 
$\hat{\mathcal{M}}_I^{d_l}(\gamma _-^l,\gamma _+^l)$
contribute dimension $2$ to
$\sum_{i=0}^ke_i+\sum_{j=1}^kd_l+2k$,
the dimension of 
$\mathcal{M}_I^{e_0+\cdots + d_k+e_k+2k}(p_-,p_+)$.

%%%%%%%%%%%%%%%%%%%%%%%%%%%%%%%%%%%%%%%%%%%%%%%%%%%%%%%%%%%%%%%%%%
\section{Bubbling off phenomena for pseudo-holomorphic curves }%%
%%%%%%%%%%%%%%%%%%%%%%%%%%%%%%%%%%%%%%%%%%%%%%%%%%%%%%%%%%%%%%%%%%

We owe most of this section to \cite{h}.
Let $(R,\infty )\times M_+$ be a convex end of $X$
such that $L_0|_{(R,\infty )\times M_+}$ and 
$L_1|_{(R,\infty )\times M_+}$ 
are isomorphic to the products of 
$(R,\infty )$ and Legendrian submanifolds.
Assume that almost complex structures on the end are
of the form $\overline{I}_t$.
\begin{lem}\label{lem5.1}%%%%%%%%%%%%%%%%%%%%%%%%%%%%%%%%%%%%%%%%%
If $\overline{u}$ is a pseudo-holomorphic strip,
then the image of $\overline{u}$ is contained in 
$X\setminus (R,\infty )\times M_+$.
\end{lem}%%%%%%%%%%%%%%%%%%%%%%%%%%%%%%%%%%%%%%%%%%%%%%%%%%%%%%%%%
{\it Proof}.
We denote $\overline{u}$ on $(R,\infty )\times M_+$
by $\overline{u}:=(\alpha ,u)$.
We can compute 
\[
u^*d\lambda _+=\frac{1}{2}
[
g_{M_+}(\pi _{\xi _+}u_{\tau },\pi _{\xi _+}u_{\tau })^{1/2}
+
g_{M_+}(\pi _{\xi _+}u_t,\pi _{\xi _+}u_t)^{1/2}
]
d\tau dt,
\]
where $u_{\tau }:=u_*(\frac{\partial }{\partial \tau })$
and
$u_t:=u_*(\frac{\partial }{\partial t })$, and
\begin{eqnarray*}
(-\Delta \alpha )d\tau \wedge dt
&=&
d(d\alpha \circ i)
\\
&=& 
d(-u^*\lambda _+)
\\
&=&
-u^*d\lambda _+
\\
&=&
-\frac{1}{2}
[
g_{M_+}(\pi _{\xi _+}u_{\tau },\pi _{\xi _+}u_{\tau })
+
g_{M_+}(\pi _{\xi _+}u_t,\pi _{\xi _+}u_t)
]
d\tau \wedge dt,
\end{eqnarray*}
where 
$\Delta =\partial ^2/\partial \tau ^2+\partial ^2/\partial t^2$.
Then, by the maximum principle, the maximum of $\alpha $
has to be achieved at a boundary point 
$p_0\in {\bf R}\times \partial [0,1]$.
Since $\overline{u}_*(\frac{\partial }{\partial \tau })_{p_0}$
is tangent to $L_0$ or $L_1$,
$\overline{u}_*(\frac{\partial }{\partial \tau })_{p_0}
\in \xi _+$.
We assume that 
almost complex structures are of the form $\overline{I}_t$,
hence $\overline{u}_*(\frac{\partial }{\partial t})_{p_0}
\in \xi _+$.
Therefore the image of $\overline{u}$ is tangent to 
$\alpha (p_0)\times M_+$ at the boundary point $p_0$,
which contradicts to the strong maximum principle.
\qed
\\
\\
{}From this lemma, it is enough for us 
to consider only concave ends.
\\

In the following, we assume that our almost complex structures
on the symplectization of $(M,\lambda )$
are of the form $\overline{I}_t$.
We recall the following important matter \cite{h} and \cite{hwzI}.
Let $\Phi :=\{\varphi :{\bf R}\to [0,1], \varphi '\geq 0\}$.
To $\varphi \in \Phi $ we associate 
$\omega _{\varphi }:=d(\varphi \lambda )$ and define
\[
E_{\Phi }(\overline{u})
:=\sup _{\varphi \in \Phi}
\int_{\Sigma }\overline{u}^*\omega _{\varphi },
\]
for $\overline{u}:\Sigma \to {\bf R}\times M$.
Note that, if $\overline{u}$ is pseudo-holomorphic, then 
\begin{eqnarray*}
\overline{u}^*\omega _{\varphi }
&=&\overline{u}^*(d\varphi \wedge \lambda +\varphi d\lambda )
\\
&=&\frac{1}{2}[(\varphi '(\alpha )(
\alpha _{\tau }^2+\alpha _t^2
+\lambda (\alpha _{\tau })^2+\lambda (\alpha _t)^2)
\\
&&
+\varphi (\alpha )
(g_M(\pi _{\xi }u_{\tau },\pi _{\xi }u_{\tau })+
g_M(\pi _{\xi }u_t,\pi _{\xi }u_t))
]d\tau \wedge dt.
\end{eqnarray*}
Hence $\overline{u}^*\omega _{\varphi }\geq 0$.
\\

Hofer proved the following 
Lemma \ref{lem5.2}, Proposition \ref{pro5.3}, 
Theorem \ref{theo5.4} and Theorem \ref{theo5.5}.
Let $\Sigma $ be ${\bf C}$ or ${\bf R}\times S^1$
and $\overline{u}=(\alpha ,u)
:\Sigma \to {\bf R}\times M$ a pseudo-holomorphic map.
\begin{lem}\label{lem5.2}%%%%%%%%%%%%%%%%%%%%%%%%%%%%%%%%%%%%%%%%%
If $E_{\Phi }(\overline{u})<\infty $ and 
$\int_{\Sigma }u^*\lambda =0$, then 
$\overline{u}$ is a constant map.
\end{lem}%%%%%%%%%%%%%%%%%%%%%%%%%%%%%%%%%%%%%%%%%%%%%%%%%%%%%%%%%
\begin{pro}\label{pro5.3}%%%%%%%%%%%%%%%%%%%%%%%%%%%%%%%%%%%%%%%%%
If $E_{\Phi }(\overline{u})<\infty $, then 
$\sup_{z\in \Sigma }|\nabla \overline{u}|<\infty $.
\end{pro}%%%%%%%%%%%%%%%%%%%%%%%%%%%%%%%%%%%%%%%%%%%%%%%%%%%%%%%%%
\begin{theo}\label{theo5.4}%%%%%%%%%%%%%%%%%%%%%%%%%%%%%%%%%%%%%%%
If there is a constant $c$ such that 
$\sup_{z\in \Sigma }|\nabla \overline{u}|<c $,
then for $\beta :=(\beta _1,\beta _2)$ there exist constants 
$c_{\beta }$ such that
$|D^{\beta }\overline{u}|\leq c_{\beta }$.
\end{theo}%%%%%%%%%%%%%%%%%%%%%%%%%%%%%%%%%%%%%%%%%%%%%%%%%%%%%%%%
Let $\phi :{\bf R}\times S^1\to {\bf C}\setminus \{0\}$ 
be a map $\phi (\tau ,t):=e^{2\pi (\tau +it)}$.
We use $\overline{v}=(\alpha ,v)$ to denote 
$\overline{u}\circ \phi :{\bf R}\times S^1\to {\bf R}\times M$.
\begin{theo}\label{theo5.5}%%%%%%%%%%%%%%%%%%%%%%%%%%%%%%%%%%%%%%%
Let $\overline{u}:{\bf C}\to {\bf R}\times M$
be a non-constant pseudo-holomorphic map 
such that $E_{\Phi }(\overline{u})<\infty $.
Then there is a closed characteristic 
$x:{\bf R}/T{\bf Z}\to M$ and a sequence $s_k\to \infty $ 
such that $x_k(t):=v(s_k,t/T)$ converges to $x(t)$
in the $C^{\infty }$ topology.
\end{theo}%%%%%%%%%%%%%%%%%%%%%%%%%%%%%%%%%%%%%%%%%%%%%%%%%%%%%%%%
Note that the above closed characteristic $x$ is contractible.
We can prove completely parallel arguments 
for $\Sigma =\{z\in {\bf C}|\mbox{Im}z\geq 0\}$
to the above results.
Let $\Lambda $ be a Legendrian submanifold and
$\overline{u}=(\alpha ,u):\Sigma \to {\bf R}\times M$ 
is pseudo-holomorphic map
such that $\overline{u}(\{z\in {\bf C} |\mbox{Im}z=0\})
\subset {\bf R}\times \Lambda $.
\begin{lem}\label{lem5.6}%%%%%%%%%%%%%%%%%%%%%%%%%%%%%%%%%%%%%%%%%
If $E_{\Phi }(\overline{u})<\infty $ and 
$\int_{\Sigma }u^*\lambda =0$, then 
$\overline{u}$ is a constant map.
\end{lem}%%%%%%%%%%%%%%%%%%%%%%%%%%%%%%%%%%%%%%%%%%%%%%%%%%%%%%%%%
\begin{pro}\label{pro5.7}%%%%%%%%%%%%%%%%%%%%%%%%%%%%%%%%%%%%%%%%%
If $E_{\Phi }(\overline{u})<\infty $, then 
$\sup_{z\in \Sigma }|\nabla \overline{u}|<\infty $.
\end{pro}%%%%%%%%%%%%%%%%%%%%%%%%%%%%%%%%%%%%%%%%%%%%%%%%%%%%%%%%%
\begin{theo}\label{theo5.8}%%%%%%%%%%%%%%%%%%%%%%%%%%%%%%%%%%%%%%%
If there is a constant $c$ such that 
$\sup_{z\in \Sigma }|\nabla \overline{u}|<c $,
then for $\beta :=(\beta _1,\beta _2)$ there exist constants 
$c_{\beta }$ such that
$|D^{\beta }\overline{u}|\leq c_{\beta }$.
\end{theo}%%%%%%%%%%%%%%%%%%%%%%%%%%%%%%%%%%%%%%%%%%%%%%%%%%%%%%%%
Let $\phi :{\bf R}\times [0,1]\to 
\{z\in {\bf R}|\mbox{Im}z\geq 0\}$ be a map 
$\phi (\tau ,t):=e^{\pi (\tau +it)}$.
We use $\overline{v}=(\alpha ,v)$ to denote 
$\overline{v}:=\overline{u}\circ \phi 
:{\bf R}\times [0,1]\to {\bf R}\times M$.
\begin{theo}\label{theo5.9}%%%%%%%%%%%%%%%%%%%%%%%%%%%%%%%%%%%%%%%
Let $\overline{u}:\{z\in {\bf C}|\mbox{{\rm Im}}z\geq 0\}
\to {\bf R}\times M$
be a non-constant pseudo-holomorphic map 
such that $E_{\Phi }(\overline{u})<\infty $.
Then there is a Reeb chord 
$x:[0,T]\to M$ from $\Lambda $ to itself
and a sequence $s_k\to \infty $
such that $x_k(t):=v(s_k,t/T)$ converges to $x(t)$
in the $C^{\infty }$ topology.
\end{theo}%%%%%%%%%%%%%%%%%%%%%%%%%%%%%%%%%%%%%%%%%%%%%%%%%%%%%%%%
Note that the above Reeb chord $x$ is contractible.
Moreover, we can similarly prove the following theorem.
Let $\Lambda _0$ and $\Lambda _1$ be Legendrian submanifolds
such that $\Lambda _0\cap \Lambda _1=\emptyset $.
\begin{theo}\label{theo5.10}%%%%%%%%%%%%%%%%%%%%%%%%%%%%%%%%%%%%%%
Let $\overline{u}:\{z\in {\bf C}|\mbox{{\rm Im}}z\geq 0\}
\setminus \{0\}\to 
{\bf R}\times M$ be a non-constant pseudo-holomorphic map
such that $E_{\Phi }(\overline{u})<\infty $ and
$\overline{u}(\{z\in {\bf C}|
\mbox{{\rm Im}}z=0,
\mbox{{\rm Re}}z>0)\subset {\bf R}\times \Lambda _0$
and $\overline{u}(\{z\in {\bf C}|
\mbox{{\rm Im}}z=0,
\mbox{{\rm Re}}z<0)\subset {\bf R}\times \Lambda _1$.
Then $\overline{u}$ satisfies one of the following:
\begin{itemize}
\item 
There are Reeb chords $x:[0,T]\to M$ 
from $\Lambda _0$ to $\Lambda _1$ and 
$x':[0,T']\to M$ 
from $\Lambda _0$ to $\Lambda _1$
and sequences $s_k\to \infty $ and $s'_k\to -\infty $
such that 
$x_k(t):=v(s_k,t/T)$ converges to $x(t)$
and
$x'_k(t):=v(s'_k,t/T')$ converges to $x'(t)$
in the $C^{\infty }$ topology.
($\lim _{k\to \infty }\alpha (s_k,[0,1])=\infty $ and
$\lim _{k\to \infty }\alpha (s'_k,[0,1])=-\infty $.)
\item 
There are Reeb chords $x:[0,T]\to M$ 
from $\Lambda _1$ to $\Lambda _0$ and 
$x':[0,T']\to M$ 
from $\Lambda _1$ to $\Lambda _0$
and sequences $s_k\to -\infty $ and $s'_k\to \infty $
such that 
$x_k(t):=v(s_k,t/T)$ converges to $x(1-t)$
and
$x'_k(t):=v(s'_k,t/T')$ converges to $x'(1-t)$
in the $C^{\infty }$ topology.
($\lim _{k\to \infty }\alpha (s_k,[0,1])=\infty $ and
$\lim _{k\to \infty }\alpha (s'_k,[0,1])=-\infty $.)
\item 
There are Reeb chords $x:[0,T]\to M$ 
from $\Lambda _1$ to $\Lambda _0$ and 
$x':[0,T']\to M$ 
from $\Lambda _0$ to $\Lambda _1$
and sequences $s_k\to -\infty $ and $s'_k\to \infty $
such that 
$x_k(t):=v(s_k,t/T)$ converges to $x(1-t)$
and
$x'_k(t):=v(s'_k,t/T')$ converges to $x'(t)$
in the $C^{\infty }$ topology.
($\lim _{k\to \infty }\alpha (s_k,[0,1])=\infty $ and
$\lim _{k\to \infty }\alpha (s'_k,[0,1])=\infty $.)
\end{itemize}
\end{theo}%%%%%%%%%%%%%%%%%%%%%%%%%%%%%%%%%%%%%%%%%%%%%%%%%%%%%%%%
(Note that the first case is ${\bf (V)}$, 
the second is $({\bf V'})$ and the third is ${\bf (IV)}$.)
\\

{}From here $(-\infty ,R]\times M$ denotes a concave end of $X$.
Let $\overline{u}_i:{\bf R}\times [0,1]\to X$ 
be a sequence of pseudo-holomorphic strips 
with the Lagrangian boundary conditions and $\{p_i\}$ 
a sequence of ${\bf R}\times [0,1]$ 
such that 
$|\nabla \overline{u}_i(p_i)|\to \infty $.
Then there is a sequence $\varepsilon _i\to 0$ such that 
\begin{eqnarray*}
\varepsilon _i|\nabla \overline{u}_i(p_i)|
&\to & \infty,
\\
|\nabla \overline{u}_i(z)|
&\leq & 2|\nabla \overline{u}_i(p_i)| 
\mbox{ for }|z-p_i|\leq \varepsilon _i,
\end{eqnarray*}
see \cite{h} and \cite{ms}.
If $\bigcup _i\{\overline{u}_i(z), |z-p_i|\leq \varepsilon _i\}$
is contained in a compact set,
then we can adopt the usual bubbling off phenomena in closed 
symplectic manifolds.
Put $\overline{u}_i:=(\alpha _i,u_i)$ on the concave end.
Let $q_i$ be a sequence such that 
$|p_i-q_i|\leq \varepsilon _i$ and $\alpha _i(q_i)\to -\infty $.
Take a sequence $\delta _i\to 0$ such that 
$\delta _i\leq \varepsilon _i$ and $\alpha _i(z)< R$ for
$|z-q_i|< \delta _i$ and 
$\overline{u}_i(\{|z-q_i|= \delta _i\})\cap 
(\{R\}\times M)\neq \emptyset $.
Then, from the mean value theorem, 
there is a point $q'_i\in \{|z-q_i|\leq \delta _i\}$
such that $\delta _i|\nabla \overline{u}_i(q'_i)|\to \infty $.
By slightly modifying $\{q_i\}$ we may assume
$\delta _i|\nabla \overline{u}_i(q_i)|\to \infty $. 
First we consider the case when
we can choose a subsequence $\{q_{i_k}\}$ such that 
$\{|z-q_{i_k}|< \delta _{i_k}\}\subset {\bf R}\times (0,1)$.
(In the following we shall use $\{q_i\}$ 
to denote the subsequence.)
Put 
\[
\overline{v}_i(z)=(b_i(z),v_i(z)):=
(\alpha _i(q_i+z/|\nabla \overline{u}_i(q_i)|)-\alpha _i(q_i),
u_i(q_i+z/|\nabla \overline{u}_i(q_i)|)),
\]
then 
\[
b_i(0)=0,\ |\nabla \overline{v}_i(0)|=1,\ 
|\nabla \overline{v}_i(z)|\leq 2 
\mbox{ {\rm for }}|z|\leq \delta _i|\nabla \overline{u}_i(q_i)|.
\]
Then we can conclude that there is a pseudo-holomorphic map
$\overline{v}=(b,v):{\bf C}\to {\bf R}\times M$
and a subsequence of $\{\overline{v}_i\}$ 
which converges to $\overline{v}$
in the $C^{\infty }_{loc}$ topology.
Secondly we consider the case of
$\mbox{Im}q_i\to 0$ and $\{|z-q_i|<\delta _i\}
\cap ({\bf R}\times \{0\})\neq \emptyset $, or
$\mbox{Im}q_i\to 1$ and $\{|z-q_i|<\delta _i\}
\cap ({\bf R}\times \{1\})\neq \emptyset $.
Then we can similarly obtain a pseudo-holomorphic map
$\overline{v}:\{z\in {\bf C}|\mbox{Im}z\geq 0\}
\to {\bf R}\times M$ and a subsequence of 
$\{\overline{v}_i\}$ which converges to 
$\overline{v}$ in the $C^{\infty }_{loc}$ topology.

Let $\overline{u}_i:{\bf R}\times [0,1]\to X$ 
be a sequence of pseudo-holomorphic strips 
with the Lagrangian boundary conditions and 
$\{q_i\}$ a sequence of ${\bf R}\times [0,1]$ 
such that $\alpha _i(q_i)\to -\infty $, 
where $\overline{u}_i:=(\alpha ,u_i)$ on a concave end.
Take a sequence $\delta _i$ such that 
$|\alpha _i(z)|<R$ for $|z-q_i|<\delta _i$ and
$\overline{u}_i(\{|z-q_i|=\delta _i\})\cap (\{R\}\times M)\neq 
\emptyset $.
If $\delta _i$ are bounded, then, from the mean value theorem,
there are points $p_i$ such that $|q_i-p_i|<\delta _i$
and $|\nabla \overline{u}_i(p_i)|\to \infty $.
In this case we can return to the previous one. 
If there is a subsequence $\delta _{i_k}\to \infty $,
then we put
\[
\overline{w}_k:=(\alpha _{i_k}(z-\mbox{Re}q_{i_k})-\alpha (q_i),
u_{i_k}(z-\mbox{Re}q_{i_k})).
\]
If there are points $r_k$ such that 
$|\nabla \overline{w}_k(r_k)|\to \infty$,
then we can also return to the previous bubbling phenomena.
Hence we assume that 
the differential of $\overline{w}_k$ are bounded.
Then there is a pseudo-holomorphic strip
$\overline{v}:{\bf R}\times [0,1]\to {\bf R}\times M$
and a subsequence of $\{\overline{w}_k\}$ which converges to 
$\overline{v}$ in the $C^{\infty }_{loc}$ topology.
\begin{pro}\label{pro5.11}%%%%%%%%%%%%%%%%%%%%%%%%%%%%%%%%%%%%%%%%
If $\overline{u}_i$ are of the type ${\bf (I)}$ or ${\bf (II)}$ or
${\bf (II')}$ or ${\bf (III)}$ or ${\bf (III')}$, and
$\int _{{\bf R}\times [0,1]}\overline{u}_i^*\omega \leq C$,
then $E_{\Phi }(\overline{v})\leq e^{-R}C$.
\end{pro}%%%%%%%%%%%%%%%%%%%%%%%%%%%%%%%%%%%%%%%%%%%%%%%%%%%%%%%%%
{\it Proof}.
Notation: $K_{[a,b]}:=[a,b]\times M$, $M_d:=d\times M$
and $C_d:=\overline{u}^{-1}((-\infty ,d]\times M)$.
Let $d_1\leq d_2\leq R$ and 
$\varphi =\varphi (\alpha ):{\bf R}\to [0,1]$ a function
such that $\varphi '\geq 0$. 
Since $\partial \overline{u}^{-1}(K_{[d_1,d_2]})
=\overline{u}^{-1}(M_{d_1})\cup 
\overline{u}^{-1}(M_{d_2})\cup 
\overline{u}^{-1}([d_1,d_2]\times \Lambda )$,
where $\Lambda $ is Legendrian,
\begin{eqnarray*}
&&
\int_{C_{d_2}}\overline{u}^*d(\varphi \lambda )
-\int_{C_{{d_1}}}\overline{u}^*d(\varphi \lambda )
\\
&=&
\int_{\overline{u}^{-1}(K_{[d_1,d_2]})}\overline{u}^*
d(\varphi \lambda ) 
\\
&=&
\int_{\partial \overline{u}^{-1}(K_{[d_1,d_2]})}
\overline{u}^*(\varphi \lambda )
\\
&=&
\int_{\overline{u}^{-1}(M_{d_2})}\overline{u}^*(\varphi \lambda ) 
-
\int_{\overline{u}^{-1}(M_{d_1})}\overline{u}^*(\varphi \lambda )
\\
&=&
\varphi (d_2)e^{-d_2}\int_{\overline{u}^{-1}(M_{d_2})}
\overline{u}^*(e^{\alpha }\lambda )
-
\varphi (d_1)e^{-d_1}\int_{\overline{u}^{-1}(M_{d_1})}
\overline{u}^*(e^{\alpha }\lambda )
\\
&=& \varphi (d_2)e^{-d_2}\int_{C_{d_2}}\overline{u}^*\omega 
-\varphi (d_1)e^{-d_1}\int_{C_{d_1}}\overline{u}^*\omega .
\end{eqnarray*}
If $\overline{u}$ satisfies the exponential decay conditions, then
$\lim_{d_1\to -\infty }\int_{C_{d_1}}
\overline{u}^*d(\varphi \lambda )
=\lim_{d_1\to -\infty }T(\varphi (d_1)-\varphi (-\infty ))=0$.
Hence 
\[
\int_{C_{d_2}}\overline{u}^*d(\varphi \lambda )
=\varphi (d_2)e^{-d_2}\int_{C_{d_2}}\overline{u}^*\omega 
-\lim_{d_1\to -\infty }\varphi (d_1)
e^{-d_1}\int_{C_{d_1}}\overline{u}^*\omega 
\leq e^{-R}C.
\]
For any $\varepsilon >0$, there are $\varphi \in \Phi$
such that 
\[
E_{\Phi }(\overline{v})-\varepsilon \leq 
\int_{\Sigma }\overline{v}^*d(\varphi \lambda ).
\]
For any compact set $K\subset \Sigma $,
there are $\overline{u}_i$ such that 
\[
\left|\int_K\overline{v}^*d(\varphi \lambda )-
\int_K
(\alpha _i\circ \phi _i-\alpha _i(q_i),u_i\circ \phi _i)^*
d(\varphi \lambda )\right|
\leq \varepsilon ,
\]
where $\phi _i:K\to C_R$ is a suitable map.
Since $\overline{u}_i^*d(\varphi \lambda )\geq 0$, 
\[
\int_K(\alpha _i\circ \phi _i-\alpha _i(q_i),
u_i\circ \phi _i)^*d(\varphi \lambda)
=\int_{\phi _i(K)}\overline{u}_i^*d(\varphi _i\lambda )
\leq \int_{C_R}\overline{u}_i^*d(\varphi _i\lambda ),
\]
where $\varphi _i(\alpha ):=\varphi (\alpha -\alpha (q_i))$.
Hence
\[
E_{\Phi }(\overline{v})\leq e^{-R}C.
\]
\qed
\\
\\
{}From this proposition 
we can apply Theorem \ref{theo5.5} to $\overline{v}$
which is the limit of the sequence $\{\overline{v}_i\}$
and also Theorem \ref{theo5.10} to $\overline{v}$,
the limit of the sequence $\{\overline{w}_k\}$.
Moreover, we can conclude the exponential decay conditions 
for $\overline{v}$ \cite{hwzI}.
\begin{pro}\label{pro5.12}%%%%%%%%%%%%%%%%%%%%%%%%%%%%%%%%%%%%%%%%
If there is an open neighborhood $U\subset S^1\times {\bf R}^{2n}$
of $S^1\times \{0\}$ and an open neighborhood $V\subset M$ of 
a closed characteristic $x$ with the minimal period $\tau _0$ and
a diffeomorphism $\varphi :U\to V$ mapping $S^1\times \{0\}$
to $x$ such that $\varphi ^*\lambda =f\lambda _0$ 
with $\lambda _0:=d\theta +\sum_{i=1}^nx_ndy_n$ and a positive
smooth function $f:U\to {\bf R}$ satisfying 
$f(\theta ,0,0)=\tau _0$ and $df(\theta ,0,0)=0$
for all $\theta \in S^1$,
then $x_k$ as in Theorem \ref{theo5.5} converges to $x$ with
the exponential decay conditions.
\end{pro}%%%%%%%%%%%%%%%%%%%%%%%%%%%%%%%%%%%%%%%%%%%%%%%%%%%%%%%%%
Similarly we can prove 
\begin{pro}\label{pro5.13}%%%%%%%%%%%%%%%%%%%%%%%%%%%%%%%%%%%%%%%%
Under Assumption \ref{ass3.3}, 
$x_k$ and $x'_k$ as in Theorem \ref{theo5.9} and \ref{theo5.10}
converge to $x$ and $x'$ with the exponential decay conditions.
\end{pro}%%%%%%%%%%%%%%%%%%%%%%%%%%%%%%%%%%%%%%%%%%%%%%%%%%%%%%%%%

Under Assumption \ref{ass3.2} there are no point $q_i$ 
such that $|\nabla \overline{u}_i(q_i)|\to \infty $ and
$\alpha _i(q_i)\to -\infty $,
and under Assumption \ref{ass3.5} $\overline{v}$, 
the limit of $\overline{w}_k$, is of the type
${\bf (IV)}$ or the trivial ones of the type 
${\bf (V)}$ and ${\bf (V')}$.
Let $\overline{u}_i:{\bf R}\times [0,1]\to X$ 
be a sequence of pseudo-holomorphic strips 
of $\mathcal{M}^2_{I}(p_-,p_+)$
and let $\{q^s_i\}_{i=1,2,\ldots }$, $s=1,2,\ldots $, be 
sequences of ${\bf R}\times [0,1]$ 
such that $\alpha _i(q^s_i)\to -\infty $.
Take a sequence $\delta ^s_i$ such that 
$|\alpha _i(z)|<R$ for $|z-q^s_i|<\delta ^s_i$ and
$\{|z-q^{s_1}_i|<\delta ^{s_1}_i\}\cap 
\{|z-q^{s_2}_i|<\delta ^{s_2}_i\}=\emptyset $ for $s_1\neq s_2$.
We assume $\delta ^s_i\to \infty $,
and put
\[
\overline{w}^s_i:=(\alpha _i(z-\mbox{Re}q^s_i)
-\alpha (q^s_i), u_i(z-\mbox{Re}q^s_i)).
\]
Then there are pseudo-holomorphic strips
$\overline{v}^s:{\bf R}\times [0,1]\to {\bf R}\times M$
of $\mathcal{M}_{IV}^{d_s}(\gamma ^s_-,\gamma ^s_+)$ 
or the trivial ones of ${\bf (V)}$ or ${\bf (V')}$
and subsequences of $\{\overline{w}^s_i\}_{i=1,2,\ldots }$ 
which converges to $\overline{v}^s$ 
in the $C^{\infty }_{loc}$ topology.
Assume that $\overline{v}^1$ and $\overline{v}^2$ 
are {\it not} the trivial ones, i.e., of the type ${\bf (IV)}$.
Then we obtain 
$\overline{x}\in \hat{\mathcal{M}}_{II'}^{d_x}(p_-,\gamma ^1_-)$,
$\overline{y}\in \hat{\mathcal{M}}_{III}^{d_y}
(\gamma ^1_+,\gamma ^2_-)$ and
$\overline{w}\in \hat{\mathcal{M}}_{II}^{d_w}(\gamma ^2_+,p_+)$
such that we can glue them 
with $\overline{v}^1$ and $\overline{v}^2$ 
to reconstruct $\overline{u}_i$.
{}From Assumption \ref{ass3.1}
we can calculate 
$d_x+d_1+d_y+d_2+d_w=2$, which contradicts to 
$d_x,d_y,d_w\geq 0$ and $d_1,d_2\geq 2$. 
Then there is at most one $\overline{v}$ of the type ${\bf (IV)}$
which appears at the limit of the sequence $\{\overline{u}_i\}$.
Similarly also for sequences of pseudo-holomorphic strips 
of $\mathcal{M}_{II}^2(\gamma _-,p_+)$ or 
$\mathcal{M}_{II'}^2(p_-,\gamma _+)$ or
$\mathcal{M}_{III}^2(\gamma _-,\gamma _+)$.
\\

Finally we completely finish proving Theorem \ref{theo3.4}.

%%%%%%%%%%%%%%%%%%%%%%%%%%%%%%%%%%%%%%%%%%%%%%%%%%%%%%%%%%%%%%%%%%
\appendix %%%%%%%%%%%%%%%%%%%%%%%%%%%%%%%%%%%%%%%%%%%%%%%%%%%%%%%%
%%%%%%%%%%%%%%%%%%%%%%%%%%%%%%%%%%%%%%%%%%%%%%%%%%%%%%%%%%%%%%%%%%

%%%%%%%%%%%%%%%%%%%%%%%%%%%%%%%%%%%%%%%%%%%%%%%%%%%%%%%%%%%%%%%%%%
\section{Newton's method}%%%%%%%%%%%%%%%%%%%%%%%%%%%%%%%%%%%%%%%%%
%%%%%%%%%%%%%%%%%%%%%%%%%%%%%%%%%%%%%%%%%%%%%%%%%%%%%%%%%%%%%%%%%%

In this appendix we adopt Newton's method with proof \cite{f1} 
and \cite{f4}.
\begin{pro}{\bf (Newton's method)}\label{app1}%%%%%%%%%%%%%%%%%%%%
Let $(E, \|\ \|_E)$ and $(F, \|\ \|_F)$ be Banach spaces and
$f:E\to F$ a smooth map.
We denote the Taylor expansion of $f$ by 
\[
f(\xi )=f(0)+Df(0)\xi +N(\xi ),
\]
and assume that $Df(0)$ has a right inverse $G:F\to E$,
$Df(0)G=\mbox{{\rm id}}_F$,
such that
\begin{equation}
\|GN(\xi )-GN(\zeta )\|_E
\leq C_N(\|\xi \|_E+\|\zeta \|_E)\|\xi -\zeta \|_E,
\label{a1}
\end{equation}
for some constant $C_N$.
Then the zero-set of $f$ in 
$B:=\{\xi \in E| \|\xi \|_E<(4C_N)^{-1} \}$
is a smooth manifold,
whose dimension is equal to that of $\mbox{{\rm Ker}} Df(0)$.
In fact,
if we put
\begin{eqnarray*}
K&=&\{\xi \in \mbox{{\rm Ker}} 
Df(0)| \| \xi \|_E<(4C_N)^{-1} \},
\\
K^{\perp }&=&\{\xi \in GF| \| \xi \|_E <(4C_N)^{-1} \},
\end{eqnarray*}
then there is a smooth map
$\phi :K\to K^{\perp }$ such that 
$f(\xi +\phi (\xi ))=0$,
and
all zeroes of $f$ in $B$ are of the form 
$\xi +\phi (\xi )$.
Moreover, we have the estimate
\begin{equation}
\| \phi (\xi )\|_E  \leq 2\| Gf(0)\|_E. 
\label{a2}
\end{equation}
On the other hand, 
if we have the inequality
\begin{equation}
\|Gf(0)\|_E\leq (8C_N)^{-1},
\label{a3}
\end{equation}
then there must exist the zeros of $f$ in $B$. 
\end{pro}
{\it Proof}.
If we put $v:=Df(0)u$ for $u\in E$, then 
$Df(0)Gv=v=Df(0)u$, 
and $Df(0)(Gv-u)=0$.
Moreover,
if we have $Gb\in \mbox{Ker} Df(0)\cap GF$, then 
$0=Df(0)Gb=b$,
and $\mbox{Ker} Df(0)\cap GF=\{0\}$.
Hence we obtain the direct decomposition
$E=\mbox{Ker} Df(0)\oplus GF$. 
Denote $E_1:=\mbox{Ker} Df(0)$,
$E_2:=GF$ and  $\mathcal{F}:=G\circ f:E_1\oplus E_2\to E_2$.
For the natural projection $\pi :E_1\oplus E_2\to E_2$,
we obtain $D\mathcal{F}(0)=GDf(0)=\pi$, and 
the Taylor expansion of $\mathcal{F}$ is
\[
\mathcal{F}(\xi )
=
\mathcal{F}(0)+\pi (\xi )+\mathcal{N}(\xi ),
\]
where $\mathcal{N}:=GN$.
The inequalities $(\ref{a1})$ $(\ref{a3})$ become
\begin{eqnarray}
\|\mathcal{N}(\xi )-\mathcal{N}(\zeta )\|_E
& \leq &
C_N(\| \xi \|_E +\| \zeta \|_E )
\| \xi -\zeta \|_E ,
\label{a4}
\\
\| \mathcal{F}(0) \|_E 
& \leq & 
(8C_N)^{-1}.
\label{a5}
\end{eqnarray}
Because $G:F\to GF$ is an isomorphism,
we will prove the proposition for $\mathcal{F}:E\to E$
instead of $f:E\to F$.

For simplicity we denote the norm $\|\ \|_E$ by $\|\ \|$.
{}From $(\ref{a4})$ we have
\begin{eqnarray*}
\| D\mathcal{N}(\xi)\zeta \|
& \leq &
\| \mathcal{N}(\xi +\zeta )-\mathcal{N}(\xi )\|
+o(\zeta )\| \zeta \|
\\
& \leq &
C_N(\| \xi +\zeta \| +\| \xi \| )
\| \zeta \| +o(\zeta )\| \zeta \|,
\end{eqnarray*}
where $o(\zeta )$ is a function such that
$\lim_{\zeta \to 0}o(\zeta )=0$, and then
\[
\|D\mathcal{N}(\xi )
(\varepsilon \zeta/\|\varepsilon \zeta \|)\|
\leq C_N(\|\xi \|+\|\xi +\varepsilon \zeta \|)
+o(\varepsilon \zeta ).
\]
The limit of the above estimate as $\varepsilon \to 0$ is
\[
\|D\mathcal{N}(\xi )(\zeta /\|\zeta \|)\|\leq 2C_N\|\xi \|,
\]
and hence the operator norm of $D\mathcal{N}(\xi )$ satisfies
\[
\|D\mathcal{N}(\xi )\|<1/2
\]
for $\xi \in B$.
Moreover,
from the differential of the Taylor expansion of $\mathcal{F}$
we obtain
\[
D\mathcal{F}(\xi )=\pi +D\mathcal{N}(\xi ),
\]
and then
\[
D\mathcal{F}(\xi )|_{E_2}=
\mbox{id}_{E_2}+D\mathcal{N}(\xi )|_{E_2}.
\]
Hence, the restriction $D\mathcal{F}(\xi )|_{E_2}:E_2\to E_2$
is an isomorphism for $\xi \in B$.
Now we use:
\begin{theo}(Implicit function theorem)%%%%%%%%%%%%%%%%%%%%%%%%%%%
Let $\mathcal{F}:E_1\oplus E_2\to E_2$ be a smooth map,
where $E_1$ and $E_2$ are Banach spaces, such that
the differential 
$D\mathcal{F}(\xi _1,\xi _2 )|_{E_2}:E_2\to E_2$ 
is an isomorphism at zeros 
of $\mathcal{F}$, i.e., $\mathcal{F}(\xi _1,\xi _2 )=0$.
Then there is a neighborhood $W_{\xi _1}\subset E_1$ of $\xi _1$
and a smooth map 
$\phi _{\xi _1}:W_{\xi _1}\to E_2$
such that, for any $\zeta \in W_{\xi _1}$, 
$\mathcal{F}(\zeta ,\phi _{\xi _1}(\zeta ))=0$.
\end{theo}%%%%%%%%%%%%%%%%%%%%%%%%%%%%%%%%%%%%%%%%%%%%%%%%%%%%%%%%
Let $\xi :=(\xi _1,\xi _2)$ and $\xi ':=(\xi _1,\xi _2')$
be in $B$ such that $\mathcal{F}(\xi )=\mathcal{F}(\xi ')=0$.
{}From the Taylor expansion of $\mathcal{F}$
\[
\mathcal{F}(0)+\xi _2+\mathcal{N}(\xi )
=
\mathcal{F}(0)+\xi _2'+\mathcal{N}(\xi ')
=0,
\]
and then 
\begin{eqnarray*}
\|\xi _2-\xi _2'\|
&=&
\|\mathcal{N}(\xi )-\mathcal{N}(\xi ')\|
\\
&\leq &
C_N(\|\xi \|+\|\xi '\|)\|\xi -\xi '\|
\\
&\leq &
\frac{1}{2}\|\xi _2-\xi _2'\|.
\end{eqnarray*}
Hence we can write the zeros of $\mathcal{F}$ in $B$ of the form
$(\zeta , \phi (\zeta )),\zeta \in E_1$, and then 
$\{\xi \in B|\mathcal{F}(\xi )=0\}$ is a smooth manifold
whose tangent spaces are isomorphic to $E_1$.

{}From the equations
\begin{eqnarray*}
\mathcal{F}(\xi _1+\phi (\xi _1))
&=&
\mathcal{F}(0)+\phi (\xi _1)+\mathcal{N}(\xi _1+\phi (\xi _1))=0
\\
\mathcal{F}(\xi _1)
&=&
\mathcal{F}(0)+0+\mathcal{N}(\xi _1)
\end{eqnarray*}
for $\xi _1\in E_1$, we obtain 
\[
\phi (\xi _1)=-\mathcal{F}(\xi _1)
-\mathcal{N}(\xi _1+\phi (\xi _1))
+\mathcal{N}(\xi _1),
\]
and then
\begin{eqnarray*}
\|\phi (\xi _1)\|
&\leq &
\|\mathcal{F}(\xi _1)\|
+\|\mathcal{N}(\xi _1+\phi (\xi _1))-\mathcal{N}(\xi _1)\|
\\
&\leq &
\|\mathcal{F}(\xi _1)\|
+C_N(\|\xi _1+\phi (\xi _1)\|+\|\xi _1\|)\|\phi (\xi _1)\|
\\
&\leq &
\|\mathcal{F}(\xi _1)\|+\frac{1}{2}\|\phi (\xi _1)\|.
\end{eqnarray*}
This is the inequality (\ref{a2}).

Finally 
we prove the existence of the zeros of $\mathcal{F}$ in $B$
when the inequality $(\ref{a5})$ holds.
Define a map $g:E\to E_2$ by
\[
g(\xi ):=-\mathcal{F}(0)-\mathcal{N}(\xi ).
\]
Because the inequality (\ref{a4})
with $\xi \in B$ and $\zeta =0$
becomes
\[
\|\mathcal{N}(\xi )\|\leq C_N\|\xi \|^2<(16C_N)^{-1}
\]
and the inequality $(\ref{a5})$ holds,
the image $g(B)$ is contained in 
$H:=\{\xi \in GF|\|\xi \|\leq 3/(16C_N)\}$.
Moreover, for 
$\xi $ and $\zeta \in \overline{B}$,
\begin{eqnarray*}
\|g(\xi )-g(\zeta )\|
&\leq &
\|\mathcal{N}(\xi )-\mathcal{N}(\zeta )\|
\\
&\leq &
C_N(\|\xi \|+\|\zeta \|)\|\xi -\zeta \|
\\
&\leq &\frac{1}{2}\|\xi -\zeta \|.
\end{eqnarray*}
Now we use:
\begin{theo}(Fixed points theorem)%%%%%%%%%%%%%%%%%%%%%%%%%%%%%%%%
Let $(X,d)$ be a complete metric space,
$g:X\to X$ a map such that
\[
d(g(x),g(y))\leq Cd(x,y),\ x,y\in X
\]
for some constants $C<1$.
Then there uniquely exists the pont $y_0$ such that
$g(y_0)=y_0$.
\end{theo}%%%%%%%%%%%%%%%%%%%%%%%%%%%%%%%%%%%%%%%%%%%%%%%%%%%%%%%%
We can find the fixed point $y_0\in H$ with respect to 
a map $g|_H:H\to H$.
Then, from the Taylor expansion of $\mathcal{F}$ 
and the definition of $g$,
\[
\mathcal{F}(y_0)=\mathcal{F}(0)+\pi (y_0)+\mathcal{N}(y_0)=0,
\]
i.e., $y_0$ is a zeros of $\mathcal{F}$.
Note that $y_0=\phi (0)$.
\qed 
%%%%%%%%%%%%%%%%%%%%%%%%%%%%%%%%%%%%%%%%%%%%%%%%%%%%%%%%%%%%%%%%%%

%%%%%%%%%%%%%%%%%%%%%%%%%%%%%%%%%%%%%%%%%%%%%%%%%%%%%%%%%%%%%%%%%%
\end{document}